\documentclass[a4paper,reqno,11pt]{amsart}
\usepackage[rgb,dvipsnames]{xcolor}
\usepackage[margin=1.18in]{geometry}
\usepackage{amssymb,amsmath,cite}
\usepackage{graphicx,capt-of,needspace}
\usepackage{tikz}

\usepackage[colorlinks=true,linkcolor=MidnightBlue,citecolor=BrickRed,urlcolor=MidnightBlue]{hyperref}
\hypersetup{bookmarksnumbered=true,pdftitle={Gradient concentration in elasticity with nearly touching cavities or inclusions of finite stiffness},pdfauthor={Haigang Li and Junhua Zhang}}
\allowdisplaybreaks[4]
\numberwithin{equation}{section}
\newtheorem{theorem}{Theorem}[section]
\newtheorem{corollary}[theorem]{Corollary}
\newtheorem{lemma}[theorem]{Lemma}
\newtheorem{prop}[theorem]{Proposition}

\theoremstyle{definition}
\newtheorem{remark}[theorem]{Remark}

\newcommand{\RN}[1]{\textup{\uppercase\expandafter{\romannumeral#1}}}

\newcommand{\norm}[1]{\lVert#1\rVert}

\newcommand{\abs}[1]{\left\lvert#1\right\rvert}

\newcommand{\Div}{\operatorname{div}}

\def\bC{\mathbb{C}}

\def\bR{\mathbb{R}} 
\def\bS{\mathbb{S}} \def\cS{\mathcal{S}}
 \def\cT{\mathcal{T}}

\def\bZ{\mathbb{Z}} 
\def\rd{\mathrm{d}}\def\ri{\mathrm{i}}
\def\Re{\operatorname{Re}}\def\Im{\operatorname{Im}}

\newcommand{\be}{\begin{equation}}
\newcommand{\ee}{\end{equation}}
\newcommand\lb{\label}\newcommand\eq{\eqref}
\newcommand{\pa}{\partial}

\newcommand{\nn}{\nonumber}

\begin{document}

\title[Gradient concentration in elasticity]{Gradient concentration in elasticity with nearly touching cavities or inclusions of finite stiffness}

\author[H.~Li]{Haigang Li}
\address[H.~Li]{School of Mathematical Sciences, Beijing Normal University, Laboratory of Mathematics and Complex Systems, Ministry of Education, Beijing 100875, China}
\email{hgli@bnu.edu.cn}

\author[J.~Zhang]{Junhua Zhang}
\address[J.~Zhang]{School of Mathematics and Computer Sciences, Nanchang University, Nanchang 330031, China}
\email{zhangjunhua@ncu.edu.cn}

\begin{abstract}
We determine the sharp rate of gradient and strain concentration between
two nearly touching traction-free cavities in a three-dimensional
isotropic elastic solid. The cavities are strictly convex bodies of
revolution interchanged by midplane reflection and approaching quadratic
contact. For gap $\varepsilon$ and Poisson's ratio $\nu$, the full
displacement gradient is $O(\varepsilon^{-\beta_\nu})$ up to both cavity
faces. A fixed smooth outer displacement attains this rate in strain for
every sufficiently small gap, even as the elastic energy in the central
ligament tends to zero. The sharp exponent is determined by a membrane
mode of angular order two for $\nu<1/5$ and by radial stretching for
$\nu>1/5$. At $\nu=1/5$ both modes attain the same rate.

For inclusions of finite stiffness, we bound the pointwise gradient error
in the ligament interior relative to the cavity and rigid solutions at
the same gap. The bounds give explicit conditions on the shear and bulk
moduli under which soft inclusions retain the sharp cavity strain rate
or stiff inclusions retain the leading term of the rigid midpoint
gradient for odd data with nonzero limiting traction resultant. For two
planar disks, the stiff-inclusion estimate extends to both interfaces,
and only the shear modulus ratio carries a logarithmic factor in the
soft-inclusion energy estimate.
\end{abstract}

\maketitle

\section{Introduction and main results}

Consider two traction-free cavities in a three-dimensional elastic body with a prescribed displacement on its outer boundary. As the gap $\varepsilon$ closes, the peak stress in the ligament can diverge even though the energy stored within a lateral distance $O(\sqrt\varepsilon)$ of the axis tends to zero. In 2007, Cardone, Nazarov, Soko{\l}owski, and Taskinen \cite[Section~6, p.~767]{CNST07} identified the asymptotic behavior of three-dimensional elastic fields in thin ligaments as an open problem. For two strictly convex cavities of revolution related by midplane reflection and approaching quadratic contact, Theorem~\ref{sph:thm:general-cavity} determines the sharp rate. Its upper bound controls the full displacement gradient up to the free faces. One fixed smooth outer displacement produces a strain of the same order for every sufficiently small gap. Even for spheres, the leading mode changes at Poisson's ratio $\nu=1/5$.

For inclusions of finite stiffness, two independent material parameters enter the problem. At fixed moduli, the computations of Babu\v{s}ka et al.~\cite{BASL99} suggested that the strain between approaching fibers remains bounded. Li and Vogelius \cite[p.~94, problem~(b)]{LV00} established a separation-independent gradient estimate for scalar equations and asked how its constant depends on ellipticity. For two conducting disks, Ammari et al.~\cite{AKLLL07} and Dong and Li \cite{DL19} obtained bounds explicit in the gap and the conductivities. For elliptic systems, Li and Nirenberg \cite{LN03} proved a separation-independent estimate at fixed ellipticity. In the rigid limit the three-dimensional gradient can grow as $(\varepsilon|\log\varepsilon|)^{-1}$ \cite{BLL17,L21}. Energy convergence to a limiting displacement \cite{AKKL13,FJ24} allows large pointwise errors in the ligament. Theorems~\ref{thm:spherical-contrast} and~\ref{thm:disks} estimate these errors in terms of the gap and the shear and bulk moduli.

\subsection{Background and motivation}

The stress between two nearby cavities depends on the load applied far
from the ligament. In planar examples, Markenscoff and Dundurs
\cite{MD92} showed how a distant boundary changes the force or moment
transmitted through the ligament. For two circular holes under uniform
normal loading at infinity, Lim and Yu \cite{LY17} determined the
singular part of the stress explicitly. In three dimensions, Tsuchida,
Nakahara, and Kodama \cite{TNK76} computed the field for two spherical
cavities under transverse uniaxial tension. To determine the sharp rate
for general outer displacements, one must identify the fastest-growing
deformation mode and find a fixed load that excites it.

The geometry at contact gives some indication of the possible singular
modes. Movchan and Nazarov \cite{MN92} studied elastic fields near sharp
defects. Nazarov and Polyakova analyzed the ligament between a cavity
and the outer boundary \cite{NP92} and the stress and strain near a
three-dimensional beak tip \cite{NP93}. Cardone, Nazarov, and
Soko{\l}owski \cite[Section~1.2]{CNS09} describe the three-dimensional
asymptotics of Nazarov and Polyakova as formal. A specific candidate
for the radial mode appears in the plate model of Syngellakis
\cite[Eq.~(37)]{Syn12}. For a spherical cavity in contact with a plane,
the radial displacement contains a term of order $r^{-1+\sqrt{2-2\nu}}$,
where $r$ is the distance to the axis. This is the radial power in our
membrane equation. As Syngellakis observed, the contact calculation
leaves open how rapidly the peak stress grows as a positive gap closes.

Maz'ya, Nazarov, and Plamenevskij \cite[Vol.~II, pp.~103--130]{MNP00}
studied the scalar Dirichlet problem in domains with thin ligaments.
The scalar analogue of a traction-free cavity is an insulating inclusion,
with zero normal flux on its boundary. For two insulating disks in the
plane, Ammari et al.~\cite{AKLLL07} proved that the gradient blow-up rate
$\varepsilon^{-1/2}$ is sharp. In dimensions $n\ge3$, Li and Yang
\cite{LY23sys} improved the upper bound to order
$\varepsilon^{-1/2+\sigma}$ for some $\sigma>0$. Weinkove
\cite{W23} used the maximum principle to make this gain explicit for
two unit balls in dimensions $n\ge4$. By analyzing a degenerate equation
for the potential averaged across the gap, Dong, Li, and Yang \cite{DLY25}
determined the sharp exponent for two balls in every dimension $n\ge3$.
Li and Zhao \cite{LZ24} proved the estimate at the optimal exponent for
more general convex inclusions. For a survey of the scalar theory,
see \cite[Section~3.1]{DY26}.

Li and Yang asked whether the scalar gain over $\varepsilon^{-1/2}$
extends to elliptic systems. For the Lam\'e system with traction-free
cavity faces, their conormal estimate
\cite[Eq.~(4.2) and Remark~4.2]{LY23sys} gives
$|\nabla u(x)|\le C(\varepsilon+|x'|^2)^{-1/2}
\|u\|_{L^\infty(G_\varepsilon(R_0))}$ in $G_\varepsilon(r_0)$,
where $0<r_0<R_0$ are fixed lateral radii. Even for two spherical
cavities in three dimensions, the sharp rate of gradient concentration
as the gap closes remained open.

For the reflected cavities of Theorem~\ref{sph:thm:general-cavity},
the full displacement gradient is bounded up to both cavity faces by
$C\|\varphi\|_{C^{2,\gamma}}\varepsilon^{-\beta_\nu}$, with optimal
exponent
\[
  \beta_\nu
  = 1-\frac12\min\left\{
      \sqrt{2-2\nu},\,
      \sqrt{7-\nu-\sqrt{\nu^2-30\nu+33}}
    \right\}.
\]
The first index comes from radial stretching and the second from a
membrane mode of angular order two. Radial stretching gives the stronger
singularity for $\nu>1/5$, the order-two mode for $\nu<1/5$, and both
attain the same rate at $\nu=1/5$. For each $\nu$, a single smooth outer
displacement produces strain of order $\varepsilon^{-\beta_\nu}$ within
a lateral distance $O(\sqrt{\varepsilon})$ of the axis for every
sufficiently small gap.

At the rigid limit, the pointwise gradient estimates of Bao, Li, and Li
\cite{BLL15,BLL17} give peak rates of order $\varepsilon^{-1/2}$ in
two dimensions and $(\varepsilon|\log\varepsilon|)^{-1}$ in three for
inclusions approaching quadratic contact. Li \cite{L21} proved matching
lower bounds for outer displacements with nonzero blow-up factor.

For inclusions of finite stiffness, Ammari, Kang, Kim, and Lee
\cite{AKKL13} proved convergence in energy to the cavity and rigid
limits, with bounds that depend on the separation when several
inclusions are present. Fu and Jing \cite{FJ24} obtained convergence
rates uniform in the period for arrays whose inclusion diameters and
spacing scale together. Here the inclusions retain their size as the
gap closes, and a small energy error can coexist with a large gradient
error in the ligament. Theorem~\ref{thm:spherical-contrast} bounds the
pointwise gradient error relative to the cavity or rigid solution at
the same gap. The bounds account for the local ligament width, the
minimum gap, and the shear and bulk modulus ratios.

\subsection{The three-dimensional cavity theorem}\label{sec:cavity-setting}

In the ball $B_L\subset\bR^3$, let $\mathcal D_\pm^0\Subset B_L$ be strictly convex bodies of revolution about the $z$ axis with $C^{2,\gamma}$ boundaries ($0<\gamma<1$). Assume $\mathcal D_-^0=\mathsf R\mathcal D_+^0$, where $\mathsf R=\operatorname{diag}(1,1,-1)$. The closures of $\mathcal D_+^0$ and $\mathcal D_-^0$ lie in $\{z\ge0\}$ and $\{z\le0\}$, respectively, and meet the midplane only at the origin. For a gap $\varepsilon>0$, set
\[
 \mathcal D_\pm^\varepsilon=\mathcal D_\pm^0\pm\tfrac12\varepsilon e_3,\qquad
 M_\varepsilon=B_L\setminus\overline{\mathcal D_+^\varepsilon\cup\mathcal D_-^\varepsilon}.
\]
Write $x=(x',z)$ and $r=|x'|$. For $r<4R$ with $R>0$ fixed and small, the facing surfaces in Figure~\ref{fig:neck} are $z=f_\pm(r)$, where $2f_+(r)=-2f_-(r)=\varepsilon+g(x')$. Assume quadratic contact, with $D^2g(0)=2\kappa I_2$ for some $\kappa>0$. Set
\[
 \delta_\varepsilon(r)=f_+(r)-f_-(r),\qquad G_\varepsilon(t)=\{r<t,\ f_-(r)<z<f_+(r)\}.
\]

\noindent\begin{minipage}{\linewidth}
\centering
\begin{tikzpicture}[font=\small,line cap=round,line join=round,
  interface/.style={draw=black!85,line width=.68pt},
  guide/.style={draw=black!48,densely dashed,line width=.36pt},
  window/.style={draw=black!52,densely dotted,line width=.52pt},
  dim/.style={<->,>=stealth,line width=.52pt},
  declare function={
    face(\r)=.08+.58*(1-sqrt(1-(\r/1.18)^2))-.065*(\r/1.18)^2;
  }]
  \def\outerRadius{2.8}
  \draw[interface] (0,0) circle (\outerRadius);
  \foreach \sgn in {1,-1}{
    \begin{scope}[yscale=\sgn]
      \filldraw[interface,fill=black!9]
        plot[domain=0:360,samples=100,smooth,variable=\ang]
        ({1.18*sin(\ang)},{.08+.58*(1-cos(\ang))-.065*sin(\ang)^2}) -- cycle;
    \end{scope}
  }
  \draw[guide] (-\outerRadius,0)--(\outerRadius,0);
  \draw[guide,->,>=stealth] (0,{-\outerRadius+.12})--(0,{\outerRadius-.15}) node[below right,text=black] {$z$};
  \node[anchor=east] at ({-\outerRadius-.04},0) {$z=0$};
  \node[anchor=south east] at (135:\outerRadius) {$\partial B_L$};
  \node at (-.63,-1.63) {$M_\varepsilon$};
  \foreach \sgn/\side in {1/+,-1/-}
    \node at (.43,{\sgn*.76}) {$\mathcal D_{\side}^{\varepsilon}$};
  \draw[window] (-.84,-.47) rectangle (.84,.47);
\hspace{1cm}
  \begin{scope}[xshift=5.63cm,scale=2.85]
    \foreach \sgn in {1,-1}{
      \begin{scope}[yscale=\sgn]
        \fill[black!9] plot[domain=-.84:.84,samples=61] (\x,{face(\x)})
          -- (.84,.47)--(-.84,.47)--cycle;
        \draw[interface] plot[domain=-.84:.84,samples=61] (\x,{face(\x)});
      \end{scope}
    }
    \draw[window] (-.84,-.47) rectangle (.84,.47);
    \draw[guide] (-.84,0)--(.84,0);
    \draw[guide] (0,-.47)--(0,.47);
    \draw[dim] (0,-.08)--(0,.08);
    \node[anchor=east,fill=white,inner sep=1pt] at (-.025,0) {$\varepsilon$};
    \draw[dim] (.56,{-face(.56)})--(.56,{face(.56)});
    \node[anchor=east,fill=white,inner sep=1pt] at (.54,0) {$\delta_\varepsilon(r)$};
    \draw[guide] (.56,{-face(.56)})--(.56,-.39);
    \draw[dim] (0,-.35)--(.56,-.35)
      node[midway,fill=black!9,inner sep=1pt] {$r$};
    \node[anchor=west] at (.86,{face(.84)}) {$z=f_+(r)$};
    \node[anchor=west] at (.86,{-face(.84)}) {$z=f_-(r)$};
  \end{scope}
\end{tikzpicture}

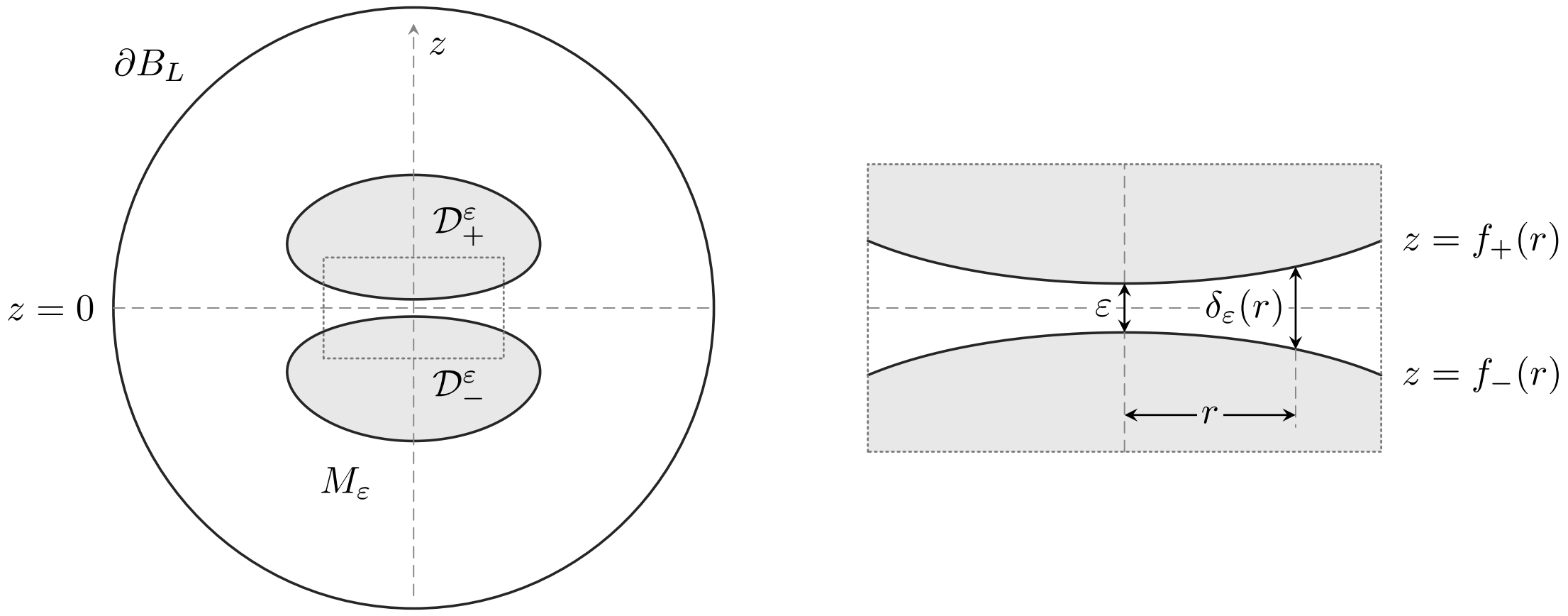
\captionof{figure}{Meridional section of the reflected cavities and an enlargement of the ligament. The minimum gap is $\varepsilon$, and $\delta_\varepsilon(r)$ is its width at distance $r$ from the axis.}\label{fig:neck}
\end{minipage}\par\medskip

For each gap $\varepsilon>0$, let $u_{\rm cav}$ solve the Lam\'e system in $M_\varepsilon$ with zero traction on the cavity faces and prescribed outer displacement $\varphi\in C^{2,\gamma}(\partial B_L;\bR^3)$. The matrix tensor is 
$$\bC_0 E=2\mu E+\lambda\operatorname{tr}(E)I_3,$$ with $\mu>0$ and $\lambda+2\mu/3>0$. We use $(\nabla u)_{ij}=\partial_j u_i$, write $e(u)=(\nabla u+\nabla u^{\mathsf T})/2$, and denote the reflection-odd part by
\[
 u_{\rm cav}^-(x)=\tfrac12\bigl(u_{\rm cav}(x)-\mathsf Ru_{\rm cav}(\mathsf Rx)\bigr).
\]
For Poisson's ratio $\nu=\lambda/[2(\lambda+\mu)]\in(-1,1/2)$, define the membrane indices and the gradient exponent by
\be\lb{main:eq:exponent}
 q_0=\sqrt{2-2\nu},\quad q_2=\sqrt{7-\nu-\sqrt{\nu^2-30\nu+33}},\quad
 \beta_\nu=1-\tfrac12\min\{q_0,q_2\}.
\ee

\begin{theorem}[Gradient and strain between traction-free cavities]\lb{sph:thm:general-cavity}
Suppose in addition that the gap profile $g$ is of class $C^5$ in a neighborhood of the origin.

\emph{(i) Upper bound.} For a fixed sufficiently small $r_0>0$ and all sufficiently small $\varepsilon>0$,
\be\lb{sph:eq:general-cavity-upper}
|\nabla u_{\rm cav}(x)|+\frac{|e(u_{\rm cav}^-)(x)|}{\sqrt{\varepsilon/\kappa}+r}\le C\|\varphi\|_{C^{2,\gamma}}(\varepsilon/\kappa+r^2)^{-\beta_\nu}
\ee
in $G_\varepsilon(r_0)$, including both matrix-side face traces. Outside the neck $|\nabla u_{\rm cav}|\le C\|\varphi\|_{C^{2,\gamma}}$. The constant is independent of $\varepsilon$ and $\varphi$.

\smallskip
\emph{(ii) Fixed-load lower bound.} There are a fixed $\varphi_*\in C^\infty(\partial B_L;\bR^3)$ and constants $c,C_*>0$ such that the corresponding solution satisfies
\be\lb{main:eq:sharp-lower}
\|e(u_{\rm cav})\|_{L^\infty(G_\varepsilon(C_*\sqrt{\varepsilon/\kappa}))}\ge c\,\varepsilon^{-\beta_\nu}
\ee
for every sufficiently small gap. The load may be chosen in the mode of angular order two for $\nu<1/5$ and in the radial sector for $\nu>1/5$. At $\nu=1/5$ either mode attains the bound. In particular, $\beta_\nu\in[1-\sqrt{2/5},1)$.
\end{theorem}

Since $|e(u_{\rm cav})|\le|\nabla u_{\rm cav}|$, the fixed-load
lower bound \eq{main:eq:sharp-lower} shows that the exponent $\beta_\nu$
in \eq{sph:eq:general-cavity-upper} is optimal. In the radial sector,
Proposition~\ref{sph:prop:scalar-strain} gives, for each $\nu$, a fixed
smooth outer displacement for which
\[
  c\,\varepsilon^{-1+q_0/2}
  \le |e_{11}(u_{\rm cav})(0)|
  = |e_{22}(u_{\rm cav})(0)|
  \le C\,\varepsilon^{-1+q_0/2}.
\]
For $\nu\ge 1/5$, these two midpoint strain components attain the full
rate $\varepsilon^{-\beta_\nu}$.

At $\nu=0$, the first Cartesian component of the radial membrane
displacement satisfies the same reduced equation as the thickness
average for two insulating balls. The midpoint estimate becomes
\[
  c\,\varepsilon^{-(2-\sqrt{2})/2}
  \le |e_{11}(u_{\rm cav})(0)|
  = |e_{22}(u_{\rm cav})(0)|
  \le C\,\varepsilon^{-(2-\sqrt{2})/2}.
\]
Thus the radial strain recovers the sharp scalar exponent for two
insulating balls \cite{DLY25}. The larger elastic exponent
$\beta_0=1-q_2/2$ is attained by a membrane mode of angular order two.

For the reflection-even load in part~\textup{(ii)}, the strain supremum grows while the elastic energy in the same central neck tends to zero. This energy is comparable to $\varepsilon^{2-2\beta_\nu}$ by Proposition~\ref{sph:thm:cav-lower}, \eq{sph:eq:general-cavity-upper} and the $O(\varepsilon^2)$ volume of the neck. Part~\textup{(i)} gains a factor of $\sqrt{\varepsilon/\kappa}+r$ for the strain of the reflection-odd bending displacement.

\subsection{Finite stiffness in three dimensions}\label{sec:main-contrast}

We compare the finite-stiffness gradient with the cavity and rigid gradients at the same gap and with the same outer displacement $\varphi$. We use the matrix tensor and reflected geometry of Section~\ref{sec:cavity-setting}. The transmission solution $u_{\rm fin}$ has continuous displacement and traction across each interface. The rigid solution $u_\infty$ minimizes the energy subject to $e(u_\infty)=0$ in both inclusions.

Both inclusions have the same isotropic tensor. For $n=2,3$, write
\[
 K=\lambda+2\mu/n,\qquad K_1=\lambda_1+2\mu_1/n,\qquad
 k_{\rm sh}=\mu_1/\mu,\qquad k_{\rm vol}=K_1/K.
\]
We assume $\mu,\mu_1,K,K_1>0$ and set
\[
 d_{\rm s}=\max\{k_{\rm sh},k_{\rm vol}\},\qquad
 d_{\rm h}=\max\{k_{\rm sh}^{-1},k_{\rm vol}^{-1}\}.
\]
For $0<\theta<1/2$, define the interior comparison region by
\[
 N_\varepsilon^\theta=\{r<r_0,\ f_-(r)+\theta\delta_\varepsilon(r)\le z\le f_+(r)-\theta\delta_\varepsilon(r)\}.
\]

\begin{theorem}[Gap and material contrast for reflected inclusions]\lb{thm:spherical-contrast}
Fix the inclusion shapes and $0<\theta<1/2$. For all sufficiently small gaps, the following bounds are uniform in $\varepsilon$, $\mu_1$, and $K_1$ in the stated ranges.

\smallskip
\emph{(i) Comparison with cavities.} Suppose $g$ is $C^5$ near the origin and $d_{\rm s}\le1$. For $x\in N_\varepsilon^\theta$,
\be\lb{eq:soft-balls-weighted}
 (\varepsilon+r^2)^{\beta_\nu}|\nabla (u_{\rm fin}-u_{\rm cav})(x)|\le C_\theta \|\varphi\|_{C^{2,\gamma}}\,\frac{\sqrt{d_{\rm s}}}{\varepsilon^2}.
\ee
The full gradient satisfies
\be\lb{eq:soft-field-explicit}
 |\nabla u_{\rm fin}(x)|\le C_\theta\|\varphi\|_{C^{2,\gamma}}
 \left(\frac{1}{(\varepsilon+r^2)^{\beta_\nu}}+
 \frac{\sqrt{d_{\rm s}}}{\varepsilon^{\beta_\nu}(\varepsilon+r^2)^2}\right).
\ee
If $d_{\rm s}\le M\varepsilon^4$ for fixed $M>0$, then
$|\nabla u_{\rm fin}(x)|\le C_{\theta,M}\|\varphi\|_{C^{2,\gamma}}(\varepsilon+r^2)^{-\beta_\nu}$.
For the same contrast range, take the fixed load $\varphi_*$ constructed in Proposition~\ref{sph:thm:cav-lower} for a leading cavity mode. There are $\theta_*\in(0,1/2)$ and $A,c>0$ such that, for every sufficiently small gap,
\[
 \|e(u_{\rm fin})\|_{L^\infty(N_\varepsilon^{\theta_*}\cap\{r<A\sqrt\varepsilon\})}
 \ge c\,\varepsilon^{-\beta_\nu}
\]

\smallskip
\emph{(ii) Comparison with rigid inclusions.} If $d_{\rm h}\le1$, then for $x\in N_\varepsilon^\theta$,
\be\lb{eq:main-balls-rigid}
 \nabla u_{\rm fin}(x)=\nabla u_\infty(x)
 +O_\theta\!\left(\frac{\|\varphi\|_{C^{2,\gamma}}\sqrt{d_{\rm h}}}
 {\varepsilon^{1/4}(\varepsilon+r^2)^{3/2}|\log\varepsilon|}\right).
\ee
In particular, for arbitrary outer data,
\be\lb{eq:hard-field-explicit}
 |\nabla u_{\rm fin}(x)|\le
 C_\theta\|\varphi\|_{C^{2,\gamma}}
 \left(\frac{1}{(\varepsilon+r^2)|\log\varepsilon|}
 +\frac{r}{\varepsilon+r^2}
 +\frac{\sqrt{d_{\rm h}}}
 {\varepsilon^{1/4}(\varepsilon+r^2)^{3/2}|\log\varepsilon|}\right).
\ee
For fixed odd data $\varphi(-x)=-\varphi(x)$ with nonzero limiting traction resultant as in Proposition~\ref{sph:thm:hard}, the midpoint gradients satisfy
\[
 \lim_{\varepsilon\to0}\varepsilon |\log\varepsilon|\nabla u_{\rm fin}(0)
 =\lim_{\varepsilon\to0}\varepsilon |\log\varepsilon|\nabla u_\infty(0)\ne0
 \qquad\text{if }d_{\rm h}=o(\varepsilon^{3/2}).
\]
\end{theorem}

\begin{remark}[The rigid-inclusion bound]
For rigid inclusions, Bao, Li, and Li \cite[Remark~1.1, (1.12)]{BLL17} proved
\be\lb{eq:rigid-BLL-exact}
 |\nabla u_\infty(x)|\le C\|\varphi\|_{C^2}
 \left(\frac{1}{(\varepsilon+r^2)|\log\varepsilon|}
 +\frac{r}{\varepsilon+r^2}\right),
 \qquad x\in G_\varepsilon(r_0).
\ee
In \eq{eq:hard-field-explicit} the contrast term is at most $\sqrt{d_{\rm h}/\varepsilon^{3/2}}$ times the first term. Hence $d_{\rm h}=O(\varepsilon^{3/2})$ suffices for the rigid spatial bound in $N_\varepsilon^\theta$. Proposition~\ref{sph:thm:hard} proves both stiff-inclusion upper bounds for general strictly convex $C^{2,\gamma}$ inclusions with nondegenerate quadratic contact. Under the coordinate symmetries and odd loading of that proposition, Li and Xu \cite[Theorem~2.4(i)]{LX25} give the coefficient in the midpoint limit.
\end{remark}

\begin{remark}[Energy estimates for the comparison errors]\lb{rem:comparison-energy}
Under the respective hypotheses of Theorem~\ref{thm:spherical-contrast}, the variational comparisons of Section~\ref{sec:common-variational} give
\[
\begin{aligned}
 \|e(u_{\rm fin}-u_{\rm cav})\|_{L^2(M_\varepsilon)}
 &\le C\|\varphi\|_{C^{2,\gamma}}\sqrt{d_{\rm s}}\,\varepsilon^{-\beta_\nu}
 &&(d_{\rm s}\le1),\\
 \|\nabla(u_{\rm fin}-u_\infty)\|_{L^2(B_L)}
 &\le C\|\varphi\|_{C^{2,\gamma}}
 \frac{\sqrt{d_{\rm h}}}{\varepsilon^{1/4}|\log\varepsilon|}
 &&(d_{\rm h}\le1).
\end{aligned}
\]
For $d_{\rm s}=O(\varepsilon^4)$, the first error is $O(\varepsilon^{2-\beta_\nu})$ under the fixed load $\varphi_*$. The cavity strain in the central neck has $L^2$ norm at least $c\,\varepsilon^{1-\beta_\nu}$, so its lower bound persists for soft inclusions.
\end{remark}

\subsection{Proof strategy for the three-dimensional estimates}

The membrane indices $q_0$ and $q_2$ come from the reflection-even part of the cavity displacement. Normalize the contact curvature to one and set $H=\varepsilon+r^2$. Let $v$ be the horizontal component of the even displacement averaged across the gap, and $M$ its in-plane stress integrated across the gap. The traction-free conditions give $\operatorname{div}'M=0$. With $\bC_{\rm ps}$ denoting the plane stress tensor and $e'$ the lateral symmetric gradient, the averaged displacement satisfies
\[
 \operatorname{div}'(H\bC_{\rm ps}e'(v))=\operatorname{div}'\mathcal E,
 \qquad \mathcal E=H\bC_{\rm ps}e'(v)-M .
\]

Set $d(r)=\sqrt\varepsilon+r$. For $0<q<2$, a bound $|\nabla u|\le C d^{q-2}$ on a larger neck gives $M=O(d^q)$ and $\mathcal E=O(d^{q+1})$ on a smaller one. The same local estimate recovers the full gradient from $\nabla'v$ up to both faces (Section~\ref{sec:membrane-error}). At contact the gradient of a mode of degree $q-1$ has order $r^{q-2}$. The extra factor of $d$ in the error estimate allows a radial Wronskian identity and a weighted estimate in angular order two to recover these powers down to the scale $\sqrt\varepsilon$. Compactness excludes the remaining modes at leading order.

For the reflection-odd displacement, the strain leaves a rotation undetermined. On a layer of thickness $h$, the estimates of Section~\ref{sec:bending} and Appendix~\ref{sec:layer-proof} bound the pointwise strain by $C h^{-1/2}\|e(u)\|_{L^2}$ and the gradient modulo a rigid motion by $C h^{-3/2}\|e(u)\|_{L^2}$. After rescaling, the limiting rotation is determined by the slope of the averaged normal displacement $w$, which solves
\[
 \operatorname{div}'\operatorname{div}'\!\bigl(H_a^3\bC_{\rm ps}D'^2w\bigr)=0,
 \qquad H_a(y)=a^2+|y|^2,
\]
where $a$ is the limiting ratio of $\sqrt\varepsilon$ to the lateral scale. The indicial roots at contact give a Liouville theorem that forces the limiting slope to vanish at the leading membrane rate.

To produce a fixed load that attains a membrane rate, we start with the dual mode of degree $-1-q$, where $q=q_0$ or $q_2$ (Section~\ref{sec:lower}). We lift it to a three-dimensional stress in equilibrium with zero traction on both faces. A correction of bounded strain energy preserves the contact singularity, whose strain energy is infinite. The resulting displacement vanishes on the outer boundary, so unique continuation forces its outer traction to be nonzero. By reciprocity, a fixed smooth outer displacement with nonzero pairing against this traction gives, for every sufficiently small gap,
\[
 \|e(u_{\rm cav})\|_{L^2(G_\varepsilon(C\sqrt\varepsilon))}\ge c\,\varepsilon^{q/2}.
\]
Since the volume of this neck is $O(\varepsilon^2)$, its strain supremum is at least $c\,\varepsilon^{-1+q/2}$.

For soft inclusions, the face estimate bounds the energy of extending $u_{\rm cav}$ into the inclusions and hence the squared $L^2$ strain error in the matrix. Lemma~\ref{lem:neck-line} controls the rotation from a fixed cross-section and propagates the full gradient bound along the center line, also for soft disks. For torsional loading on balls, a scalar reduction in five dimensions gives a bounded cavity gradient and a soft comparison even for unequal radii (Proposition~\ref{sph:thm:torsion}).

For stiff inclusions, the rigid traction $t_{\infty,i}$ on $\partial\mathcal D_i^\varepsilon$ annihilates rigid motions. Korn's inequality on each inclusion gives
\[
 \|\nabla(u_{\rm fin}-u_\infty)\|_{L^2(B_L)}
 \le C\sqrt{d_{\rm h}}\sum_{i=\pm}\|t_{\infty,i}\|_{H^{-1/2}(\partial\mathcal D_i^\varepsilon)} .
\]
On the facing caps the $L^{4/3}$ norm of the traction is bounded by $C\|\varphi\|_{C^{2,\gamma}}\varepsilon^{-1/4}/|\log\varepsilon|$. Interior estimates on balls of radius comparable to the local gap width then give \eq{eq:main-balls-rigid} (Section~\ref{sec:spherical-comparisons}).

\subsection{Two planar disks}

For two disks, the pointwise comparisons extend to both interfaces. The stiff-inclusion bound depends on the sum of the inverse shear and bulk modulus ratios. For the corresponding scalar conductivity estimates, see Ammari et al.~\cite{AKLLL07} and Dong and Li \cite{DL19}.

Let $D_1=\mathcal D_+^\varepsilon$ and $D_2=\mathcal D_-^\varepsilon$ be the unit disks centered at $(1+\varepsilon/2,0)$ and $(-1-\varepsilon/2,0)$, respectively. Let $\Omega_1\Subset\Omega$ be fixed bounded domains containing both disks, with $\partial\Omega$ of class $C^{2,\gamma}$ and the disks a fixed positive distance from $\partial\Omega_1$. Set $M_\varepsilon=\Omega\setminus\overline{D_1\cup D_2}$. With $t=x_2$ the signed tangential coordinate, the gap width and the neck are
\[
 \delta_\varepsilon(t)=\varepsilon+2(1-\sqrt{1-t^2}),\qquad
 G_\varepsilon(R)=\{(x_1,t): |t|<R,\ |x_1|<\delta_\varepsilon(t)/2\}.
\]
For $|t|<r_0$, with $r_0$ fixed and sufficiently small, $\delta_\varepsilon(t)$ is comparable to $H(t)=\varepsilon+t^2$.

On symmetric $2\times2$ matrices let $\bC_j E=2\mu_j E+\lambda_j(\operatorname{tr}E)I_2$, with $(\lambda_0,\mu_0)=(\lambda,\mu)$. We use the Frobenius norm. We use the positive shear and bulk moduli of Section~\ref{sec:main-contrast}, with $K=\lambda+\mu$ and $K_1=\lambda_1+\mu_1$ in the plane. For prescribed outer displacement $\varphi\in C^{1,\gamma}(\partial\Omega;\bR^2)$, the transmission problem is
\[
 -\Div(\bC_\varepsilon e(u_{\rm fin}))=0\quad\text{in }\Omega,\qquad
 u_{\rm fin}=\varphi\quad\text{on }\partial\Omega,
 \qquad
 \bC_\varepsilon=\bC_0\chi_{M_\varepsilon}+\bC_1(\chi_{D_1}+\chi_{D_2}).
\]
The cavity and rigid solutions have the same trace. When comparing gaps, we write $u_{\rm cav}(\cdot\,;\varepsilon)$ and $u_{\rm fin}(\cdot\,;\varepsilon)$. For soft disks set
\be\lb{soft:eq:order}
 \mathcal A_\varepsilon=\sqrt{(1+|\log\varepsilon|)k_{\rm sh}+k_{\rm vol}},\qquad
 0<k_{\rm sh},k_{\rm vol}<1.
\ee

\begin{theorem}[Gap and contrast estimates for two disks]\label{thm:disks}
The following estimates hold for all sufficiently small gaps, with constants independent of $\varepsilon$, $k_{\rm sh}$, and $k_{\rm vol}$ in the indicated ranges.

\smallskip
\emph{(i) Stiff inclusions.} There is a constant $k_*>1$, depending only on the matrix material and the geometry, such that for $k_{\rm sh},k_{\rm vol}\ge k_*$,
\be\lb{finite:eq:mat}
 |\nabla u_{\rm fin}(x_1,t)|\le
 \frac{C\|\varphi\|_{C^{1,\gamma}}}{\sqrt\varepsilon+|t|+k_{\rm sh}^{-1}+k_{\rm vol}^{-1}}\qquad ((x_1,t)\in G_\varepsilon(r_0)).
\ee
On $(\Omega_1\cap M_\varepsilon)\setminus G_\varepsilon(r_0)$ the gradient is bounded by $C\|\varphi\|_{C^{1,\gamma}}$. In $G_\varepsilon(r_0)$,
\be\lb{finite:eq:diff}
 |\nabla(u_{\rm fin}-u_\infty)(x_1,t)|
 \le \frac{C\|\varphi\|_{C^{1,\gamma}}}{\sqrt\varepsilon+|t|}
           \min\!\left\{1,\frac{k_{\rm sh}^{-1}+k_{\rm vol}^{-1}}{\sqrt\varepsilon}\right\}.
\ee
Both neck estimates include the matrix-side traces on the circles.

\smallskip
\emph{(ii) Soft inclusions.} If $0<k_{\rm sh},k_{\rm vol}<1$, then
\be\lb{soft:eq:energycompare}
 \|e(u_{\rm fin}-u_{\rm cav})\|_{L^2(M_\varepsilon)}
 \le C\|\varphi\|_{C^{1,\gamma}}\min\{1,\mathcal A_\varepsilon\}.
\ee
If also $d_{\rm s}\le c_0\sqrt\varepsilon$ for a sufficiently small fixed $c_0>0$, then, up to both matrix-side traces,
\be\lb{soft:eq:grad-difference}
 \|\nabla(u_{\rm fin}-u_{\rm cav})\|_{L^\infty(\Omega_1\cap M_\varepsilon)}
 \le C\|\varphi\|_{C^{1,\gamma}}\varepsilon^{-5/4}\mathcal A_\varepsilon .
\ee
\end{theorem}

The disk traction map separates the moduli in Fourier variables. Its negative modes depend only on shear. The real part of the first displacement coefficient gives the mean normal displacement and has traction multiplier $2K_1$ in the inclusion. For stiff disks, the interaction between the two circles has a $4\times4$ principal symbol in bipolar coordinates. A weighted bound for its inverse yields \eq{finite:eq:mat}.

For soft disks, pairing the cavity trace with its inclusion traction gives the extension energy. The mean normal displacement remains bounded as the gap closes, so only the shear term carries the logarithm in \eq{soft:eq:energycompare}. By Corollary~\ref{cor:disk-soft-transfer}, the soft disks retain the cavity gradient bound if $\mathcal A_\varepsilon^2=O(\varepsilon^{3/2})$. If $\mathcal A_\varepsilon^2=o(\varepsilon^{3/2})$, the error relative to the cavity profile is $o(\varepsilon^{-1/2})$ on $|t|\le L\sqrt\varepsilon$ in the neck, for each fixed $L>0$ and fixed outer data.

\section{Energy comparisons and gradient propagation}\label{sec:common-variational}
\subsection{Variational comparisons}
The two comparisons use different boundary quantities. For a rigid inclusion the relevant quantity is its traction acting on displacements modulo rigid motions. For a cavity it is the energy needed to extend the traction-free displacement across the interface.

Let $\Omega\subset\bR^n$, $n=2,3$, be a fixed bounded Lipschitz domain. The inclusions $B_1,B_2$ are translates of fixed bounded Lipschitz domains, with disjoint closures. Both stay a fixed positive distance from $\partial\Omega$. Write $M_\varepsilon=\Omega\setminus\overline{B_1\cup B_2}$ and $\Gamma_i=\partial B_i$. For $j=0,1$, assume $\mu_j>0$ and $K_j=\lambda_j+2\mu_j/n>0$, with $(\lambda_0,\mu_0)=(\lambda,\mu)$. If $E$ is symmetric, then
\be\lb{sph:eq:energy-density}
\bC_jE:E=2\mu_j\left|E-\frac{\operatorname{tr}E}{n}I\right|^2+K_j(\operatorname{tr}E)^2.
\ee
The tensor order is therefore equivalent to the componentwise order of the shear and bulk moduli. In particular, $\bC_1\ge d_{\rm h}^{-1}\bC_0$ and $\bC_1\le d_{\rm s}\bC_0$, with the dimension-dependent bulk moduli above. Write $u_{\rm fin},u_\infty,u_{\rm cav}$ for the finite-stiffness, rigid and traction-free solutions with common outer trace $\varphi\in H^{1/2}(\partial\Omega;\bR^n)$.

A lift of the outer trace supported in a fixed collar has zero strain in the inclusions. Minimality bounds the $L^2(M_\varepsilon)$ norm of the strain of each solution by $C\|\varphi\|_{H^{1/2}(\partial\Omega)}$, uniformly in the separation and inclusion moduli. When $\bC_1\ge\bC_0$, Korn's inequality applied after subtracting the lift gives a uniform $H^1$ bound for $u_{\rm fin}$ on the entire body.

For a finite or rigid inclusion, let $t_i$ denote its matrix traction with normal directed out of $B_i$. The traction $t_i$ annihilates every rigid displacement on $\Gamma_i$. For finite stiffness this follows from traction continuity and the weak equation in $B_i$. For rigid inclusions it follows from stationarity of the matrix energy under independent rigid variations, which the positive gap permits.

\begin{lemma}[Comparison with rigidity and zero traction]
Let $v=u_{\rm fin}-u_\infty$. If $L$ is a matrix displacement of zero outer trace and $v-L$ has rigid traces on both inclusions, then
\be\lb{eq:common-rigid-lift}
\int_{M_\varepsilon}\bC_0e(v):e(v)\,\rd x\le\int_{M_\varepsilon}\bC_0e(L):e(L)\,\rd x.
\ee
If $d_{\rm h}\le1$, then
\be\lb{eq:common-rigid-traction}
\|\nabla v\|_{L^2(\Omega)}\le C\sqrt{d_{\rm h}}\,\sum_i\|t_{\infty,i}\|_{H^{-1/2}(\Gamma_i)}.
\ee
Here $t_{\infty,i}$ is the rigid solution's matrix traction. The constant depends on the reference inclusions and the outer domain and is uniform under translations of the inclusions.

If $U_{\rm cav}$ is any $H^1$ extension of $u_{\rm cav}$ across the cavities, then
\be\lb{eq:common-soft-lift}
\begin{aligned}
 &\int_{M_\varepsilon}\bC_0e(u_{\rm fin}-u_{\rm cav}):e(u_{\rm fin}-u_{\rm cav})\,\rd x
 +\sum_i\int_{B_i}\bC_1e(u_{\rm fin}):e(u_{\rm fin})\,\rd x\\
 &\qquad\le\sum_i\int_{B_i}\bC_1e(U_{\rm cav}):e(U_{\rm cav})\,\rd x.
\end{aligned}
\ee
\end{lemma}

\begin{proof}
Test the matrix equation for $v$ against $v-L$. Rigid-motion balance cancels the boundary pairing, giving
\[
\int_{M_\varepsilon}\bC_0e(v):e(v)\,\rd x=\int_{M_\varepsilon}\bC_0e(v):e(L)\,\rd x.
\]
Cauchy--Schwarz gives \eq{eq:common-rigid-lift}.

Let $Q^2$ be the sum of the matrix energy of $v$ and its inclusion energies, with their respective tensors. Since $e(u_\infty)=0$ in the inclusions, the finite variational equation tested against $v$ gives
\[
Q^2=\sum_i\langle t_{\infty,i},v\rangle_{\Gamma_i}.
\]
The sign uses the normal pointing out of the inclusion. Subtract a rigid displacement on each $\Gamma_i$ from the second argument. Korn's inequality and the trace theorem on each inclusion bound the resulting $H^{1/2}$ norm by $C\|e(v)\|_{L^2(B_i)}$. Translation preserves the Korn and trace constants of each inclusion. By comparison of the tensors, $\|e(v)\|_{L^2(B_i)}\le C\sqrt{d_{\rm h}}\,Q$. The fixed-domain Korn inequality for $v$ with zero outer trace then yields \eq{eq:common-rigid-traction}.

Finally, use $U_{\rm cav}$ as a competitor for $u_{\rm fin}$. When the matrix energy of $u_{\rm fin}$ is expanded about $u_{\rm cav}$, the cross term vanishes by the traction-free variational equation, since $u_{\rm fin}-u_{\rm cav}$ has zero outer trace. Canceling the matrix energy of $u_{\rm cav}$ proves \eq{eq:common-soft-lift}.
\end{proof}

\subsection{Gradient control along the center line}
An $L^2$ bound for the strain leaves the rigid rotation undetermined. A bound for the full gradient at one fixed section supplies the missing normalization.

\begin{lemma}[Gradient propagation along the center line]\label{lem:neck-line}
Let $n\in\{2,3\}$, $x(t)=x_0+t\xi$ with $|\xi|=1$, and $\ell(t)=\varepsilon+t^2$, where $0<\varepsilon\le1$. Suppose $B_{c\ell(t)}(x(t))\subset G$ for $0\le t\le r_0$, with fixed $c,r_0>0$. If $v$ solves the homogeneous matrix Lam\'e system in $G$ and $E=\|e(v)\|_{L^2(G)}$, then
\be\label{eq:neck-line-gradient}
 |\nabla v(x(r))|\le|\nabla v(x(r_0))|+CE(\sqrt\varepsilon+r)^{-n-1/2},
 \qquad 0\le r\le r_0 .
\ee
The constant depends only on $c$, $r_0$, $n$, and the matrix tensor.
\end{lemma}

\begin{proof}
We may reduce $c$ by a fixed factor. Subtract a rigid displacement on $B_{c\ell(t)}(x(t))$. Scaled Korn and interior estimates give
\[
 |\nabla^2v(x(t))|\le C\ell(t)^{-1-n/2}E_t,\qquad
 E_t=\|e(v)\|_{L^2(B_{c\ell(t)}(x(t)))}.
\]
For $y$ in this ball, its projection $a=(y-x_0)\cdot\xi$ satisfies $|a-t|\le c\ell(t)$. Since $\ell$ has uniformly bounded derivative on a fixed interval containing all such $a$, the quantities $\ell(a)$ and $\ell(t)$ are comparable. Thus the admissible centers lie in an interval of length at most $C\ell(a)$. Fubini's theorem yields
\[
 \int_0^{r_0} E_t^2\frac{\rd t}{\ell(t)}\le CE^2 .
\]
Integrating $\nabla^2 v$ along the line and applying Cauchy--Schwarz with this measure gives a factor
\[
 \left(\int_r^{r_0}\ell(t)^{-n-1}\,\rd t\right)^{1/2}
 \le C(\sqrt\varepsilon+r)^{-n-1/2},
\]
since $2\ell(t)\ge(\sqrt\varepsilon+t)^2$. This proves \eq{eq:neck-line-gradient}.
\end{proof}

\section{Sectional forces and the membrane model}\label{sec:membrane}
Integration of the equilibrium equation across the ligament gives exact balance laws for the sectional forces. The constitutive law for the in-plane force then identifies the membrane operator. We first solve its radial equation, whose regular solution will enter both the gradient estimate and the midpoint strain profile.

\subsection{The gap scale and estimates away from contact}
For the cavity problem we use the reflected $C^{2,\gamma}$ geometry, outer trace and matrix tensor of Section~\ref{sec:cavity-setting}. The dilation $y=\kappa x$ normalizes the quadratic coefficient to one and replaces $\varepsilon$ by $\kappa\varepsilon$, preserving the Lam\'e system and $\nu$. Through Section~\ref{sec:lower} we take $\kappa=1$, write $G_\varepsilon=G_\varepsilon(r_0)$, and suppress the subscript on $\delta_\varepsilon$. The reflection-even part $u_{\rm cav}^+=u_{\rm cav}-u_{\rm cav}^-$ has even horizontal and odd vertical components in $z$. Set $s_\varepsilon=\sqrt\varepsilon$, $H(r)=\varepsilon+r^2$, and $d(r)=s_\varepsilon+r$. The $C^{2,\gamma}$ graphs and the quadratic contact give $\delta(r)=H(r)+O(r^{2+\gamma})$ and $\delta'(r)=2r+O(r^{1+\gamma})$. In particular the gap is increasing for small $r>0$.
Write $A_t=\{x'\colon t<|x'|<2t\}$ for the lateral annulus. The subscript zero denotes the contact configuration.

A lift $\widehat\varphi$ supported in a fixed outer collar is an admissible competitor for the cavity energy. Hence
\be\nn
\int_{M_\varepsilon}\bC_0e(u_{\rm cav}):e(u_{\rm cav})\le\int_{M_\varepsilon}\bC_0e(\widehat\varphi):e(\widehat\varphi)\le C\|\varphi\|_{C^{2,\gamma}}^2.
\ee

Outside a fixed small ball about the contact point, the two cavity faces are uniformly separated and have uniform $C^{2,\gamma}$ charts. The truncated matrix retains the outer Dirichlet boundary. The truncating sphere may be chosen transverse to the cavity faces. These domains are related by uniformly small $C^1$ deformations, so Korn anchored at the outer boundary has a uniform constant. On each such truncated domain, it bounds the $H^1$ norm of every displacement of zero outer trace by its strain norm, uniformly in $\varepsilon$. Applying it to $u_{\rm cav}-\widehat\varphi$ gives a local $H^1$ bound of $C\|\varphi\|_{C^{2,\gamma}}$. Interior and traction boundary regularity \cite{ADN64} then give, for every fixed sufficiently small $t>0$,
\be\lb{geom:eq:collar}
\sup_{M_\varepsilon\setminus G_\varepsilon(t)}(|u_{\rm cav}|+|\nabla u_{\rm cav}|)\le C_t\|\varphi\|_{C^{2,\gamma}}.
\ee
The bound includes the separated cavity faces.

In the reflected geometry, rotations about the axis and reflection in the midplane preserve the matrix, the equation and the zero cavity tractions. In particular $x\mapsto\mathsf R u_{\rm cav}(\mathsf R x)$ solves the same traction problem with the reflected outer trace. The two parity components therefore satisfy separate boundary value problems. All angular and reflection projections used below have data norms controlled by $C\|\varphi\|_{C^{2,\gamma}}$.

\subsection{Sectional forces and the membrane equation}
The sectional balance laws below hold for arbitrary outer data. Write $I_\varepsilon(x')=(-\delta(r)/2,\delta(r)/2)$ for the thickness interval, $\sigma=\bC_0e(u_{\rm cav})$ for the stress, and $\nabla'$ and $\operatorname{div}'$ for the derivatives in $x'$. We write $e^\prime$ for the lateral symmetric gradient. With $\alpha,\beta\in\{1,2\}$, define
\[
M_{\alpha\beta}=\int_{I_\varepsilon}\sigma_{\alpha\beta}\,\mathrm dz,\quad\mathcal M_{\alpha\beta}=\int_{I_\varepsilon}z\sigma_{\alpha\beta}\,\mathrm dz,\quad Q_\alpha=\int_{I_\varepsilon}\sigma_{\alpha3}\,\mathrm dz,\quad P_\alpha=\int_{I_\varepsilon}z\sigma_{\alpha3}\,\mathrm dz.
\]
The plane stress tensor associated with the matrix is
\[
\bC_{\rm ps}E=2\mu E+\frac{2\lambda\mu}{\lambda+2\mu}(\operatorname{tr}E)I_2\qquad(E\in\bR_{\rm sym}^{2\times2}).
\]

Put $\eta_{\rm ps}=\lambda/(\lambda+2\mu)$. For the in-plane displacement $u_{\rm cav}'=(u_{{\rm cav},1},u_{{\rm cav},2})$, set
\[
v=\frac{1}{\delta}\int_{I_\varepsilon}u_{\rm cav}'\,\mathrm dz,\qquad b=\frac{u_{\rm cav}'(x',\delta/2)+u_{\rm cav}'(x',-\delta/2)}{2}-v,
\]
where both traces are taken from the matrix. We write $a\odot c=(a\otimes c+c\otimes a)/2$.

\begin{prop}[Balance and constitutive law for the sectional forces]
In the disk $r<r_0$,
\be\lb{sph:eq:resultant-equilibrium}
\operatorname{div}'M=0,\qquad\operatorname{div}'\mathcal M=Q,\qquad\operatorname{div}'Q=0,\qquad\int_{I_\varepsilon}\sigma_{33}\,\mathrm dz=\operatorname{div}'P.
\ee
The in-plane resultant is
\be\lb{sph:eq:resultant-constitutive}
M=\delta\,\bC_{\rm ps}e'(v)-\bC_{\rm ps}(b\odot\nabla'\delta)+\frac{\lambda}{\lambda+2\mu}(\operatorname{div}'P)I_2.
\ee
\end{prop}

\begin{proof}
In the weak equilibrium equation use test displacements $(\psi,0)$, $(z\psi,0)$, $(0,0,\zeta)$ and $(0,0,z\zeta)$, with the in-plane test functions supported in $r<r_0$. They may have nonzero trace on the two cavity faces since the tractions there vanish. Integrating in $z$ yields \eq{sph:eq:resultant-equilibrium}, hence $\operatorname{div}'\operatorname{div}'\mathcal M=0$. This test uses only the free tractions and works for any two graphs $f_-,f_+$. In particular the moment about the middle height $c=(f_++f_-)/2$ satisfies
\[
\operatorname{div}'(\mathcal M-cM)=Q-M\nabla'c.
\]
Thus a varying middle height couples the in-plane force to the bending moment through $-M\nabla'c$. Here reflection gives $c=0$.

Next,
\[
(\sigma_{\alpha\beta})_{\alpha,\beta\le2}=\bC_{\rm ps}\bigl((e(u_{\rm cav})_{\alpha\beta})_{\alpha,\beta\le2}\bigr)+\frac{\lambda}{\lambda+2\mu}\sigma_{33}I_2.
\]
Differentiating the integral of $u_{\rm cav}'$ at the moving endpoints yields
\[
\int_{I_\varepsilon}(e(u_{\rm cav})_{\alpha\beta})_{\alpha,\beta\le2}\,\rd z=\delta e'(v)-b\odot\nabla'\delta.
\]
Substitution gives \eq{sph:eq:resultant-constitutive}.
\end{proof}

\begin{remark}
For related thin-ligament asymptotics, see Nazarov and Polyakova \cite{NP92}. In a symmetric layer of thickness $H$, the first thickness moment vanishes. The zeroth moment $H$ and the second moment $H^3/12$ give the weights in the membrane equation \eq{sph:eq:membrane-model} and the bending equation \eq{layer:eq:general-bending-limit}. The constant $1/12$ cancels in the homogeneous bending equation.
\end{remark}

\subsection{The radial membrane model}
Rotational averaging and reflection in a vertical plane select the radial meridional sector. Here $u_{\rm cav}=u_r e_r+u_3 e_3$, where $u_r(r,z)$ is even and $u_3(r,z)$ is odd in $z$.

The leading term $\delta\bC_{\rm ps}e'(v)$ in \eq{sph:eq:resultant-constitutive}, with $\delta$ replaced by $H(r)=\varepsilon+r^2$, gives the homogeneous membrane equation
\be\lb{sph:eq:membrane-model}
\operatorname{div}'\bigl(H\bC_{\rm ps}e'(v)\bigr)=0.
\ee
With $q_0$ as in \eq{main:eq:exponent}, put $p=q_0-1$ and $a=(1-p)/2$. Since $-1<\nu<1/2$, we have $0<p<1$ and $0<a<1/2$.

\begin{prop}[Sharp radial estimate for the membrane model]\lb{sph:prop:radial-model}
Fix $0<r_1<r_0$. For $\varepsilon>0$, let $v(x')=V_\varepsilon(r)e_r$ solve \eq{sph:eq:membrane-model} in $r<r_1$, with $V_\varepsilon(r_1)=1$ and $V_\varepsilon(r)/r$ bounded at the origin. There are $c,C>0$, independent of small $\varepsilon>0$, such that
\be\lb{sph:eq:radial-model-bound}
c\,H(r)^{-a}\le V_\varepsilon'(r)\le\frac{V_\varepsilon(r)}{r}\le C\,H(r)^{-a},\qquad0\le r\le r_1,
\ee
where the quotient at $r=0$ is defined by its limit. For $\varepsilon=0$, the finite membrane energy solution with the same outer value is $V_0(r)=(r/r_1)^p$.
\end{prop}

\begin{proof}
The radial and azimuthal stresses in $\bC_{\rm ps}e'(V_\varepsilon e_r)$ are a common positive factor times $V_\varepsilon'+\nu V_\varepsilon/r$ and $V_\varepsilon/r+\nu V_\varepsilon'$, respectively. Radial equilibrium gives
\be\lb{sph:eq:radial-model-ode}
\bigl(rH(V_\varepsilon'+\nu V_\varepsilon/r)\bigr)'=H(V_\varepsilon/r+\nu V_\varepsilon').
\ee
At $\varepsilon=0$, substituting $r^s$ gives $s^2+2s+2\nu-1=0$. The root $p$ has finite energy $\int_0^{r_1}r^3(|V_0'|^2+|V_0/r|^2)\,\mathrm dr$, while the other root, $-1-\sqrt{2-2\nu}$, has infinite energy.

For $\varepsilon>0$ use $t=r/\sqrt\varepsilon$ and $V_\varepsilon(r)=Y_*(t)/Y_*(r_1/\sqrt\varepsilon)$, where $Y_*(t)=t\,F(t^2)$ and
\be\lb{sph:eq:radial-model-integral}
F(s)=\frac{\int_0^1q^{1-a}(1-q)^{a-1}(1+sq)^{-a}\,\mathrm dq}{\int_0^1q^{1-a}(1-q)^{a-1}\,\mathrm dq}.
\ee
Indeed, the differential equation for $F$ is
\[
s(1+s)F''+(2+3s)F'+a(2-a)F=0.
\]
On the numerator integrand this operator gives $a\,\partial_q[q^{2-a}(1-q)^a(1+sq)^{-a-1}]$. Its integral vanishes, and the normalization gives $F(0)=1$. The second independent solution of the equation for $Y_*$ is singular like $t^{-1}$ at zero, as follows from reduction of order using $Y_*(t)/t\to1$. The regularity condition therefore selects the stated $Y_*$.

For every $s\ge0$, the integral gives
\[
(1+s)^{-a}\le F(s)\le C_a(1+s)^{-a},\qquad(1-2a)F(s)\le F(s)+2sF'(s)\le F(s).
\]
The upper estimate for $F$ follows from $(1+sq)^{-a}\le s^{-a}q^{-a}$ for $s\ge1$. The derivative inequalities follow by differentiating inside the integral, where the multiplier is $1-2asq/(1+sq)$. Since $V_\varepsilon(r)/r=F(r^2/\varepsilon)/(r_1F(r_1^2/\varepsilon))$, these estimates prove \eq{sph:eq:radial-model-bound}.
\end{proof}

\section{From the membrane model to the gradient}
The membrane equation controls the thickness average of the displacement. To recover the three-dimensional gradient we first estimate the error of the plane stress approximation on a rescaled layer. The same estimate controls the forcing in the membrane equation, so that a weighted bound for its solution yields a bound for the full gradient up to both faces.

\subsection{Local estimates on a symmetric thin layer}
We use the flattened coordinates
\[
 \begin{aligned}
 \Phi_\tau(y,t)&=(y,\tau a(y)t), & W&=U\circ\Phi_\tau,\\
 n&=(-\tau t\nabla_y a,1), & G_y&=(\partial_{y_1}W,\partial_{y_2}W,0).
 \end{aligned}
\]
For a symmetric elasticity tensor $\bC$, write $\sigma=\bC e(U)$ and $\mathcal T=(\sigma\circ\Phi_\tau)n$. In the flattened cylinder the lateral fluxes are $\tau a\,\sigma_{\cdot\alpha}\circ\Phi_\tau$, and the transverse flux is $\mathcal T$. Equilibrium and the free-traction conditions give
\be\label{layer:eq:conormal-balance}
 \partial_t\mathcal T=-\tau\sum_{\alpha=1}^2\partial_{y_\alpha}(a\,\sigma_{\cdot\alpha}\circ\Phi_\tau),
 \qquad \mathcal T|_{t=\pm1}=0 .
\ee
Since $\nabla U\circ\Phi_\tau=G_y+(\partial_Z U\circ\Phi_\tau)\otimes n$, we have
\be\label{layer:eq:conormal-recovery}
 \mathcal T=\bC\operatorname{sym}(G_y)n+\mathsf A_{\bC}(n)(\partial_ZU\circ\Phi_\tau),
 \qquad \mathsf A_{\bC}(n)p=(\bC\operatorname{sym}(p\otimes n))n .
\ee
The tensor is evaluated at $\Phi_\tau(y,t)$. Ellipticity gives $p\cdot\mathsf A_{\bC}(n)p\ge c|n|^2|p|^2/2$. Thus \eq{layer:eq:conormal-recovery} determines $\partial_ZU$ uniformly in the layer thickness.

\begin{lemma}[Membrane approximation through the free faces]\lb{sph:lem:even-local}
Let $D_1\Subset D$ be disks or concentric annuli. Let the half-thickness $a$ have fixed positive lower and upper bounds and bounded derivatives through order five. On $\mathcal P_\tau=\{y\in D,\ |Z|<\tau a(y)\}$, let $U$ solve the isotropic Lam\'e equation with zero traction on both faces, with $U'$ even and $U_3$ odd in $Z$. Put $W(y,t)=U(y,\tau a(y)t)$, $v=\overline W'=\tfrac12\int_{-1}^1W'\,\rd t$, and $\sigma=\bC_0e(U)(y,\tau a(y)t)$. On $D_1\times[-1,1]$,
\be\lb{sph:eq:layer-approx}
\left\|\nabla U(y,\tau a(y)t)-\begin{pmatrix}
\nabla_yv&0\\
0&-\eta_{\rm ps}\operatorname{div}_yv
\end{pmatrix}\right\|_{C_y^1L_t^\infty}+\|\sigma_{33}\|_{C_y^1L_t^\infty}\le C\sqrt\tau\,m.
\ee
Here $m=\|\nabla U\|_{L^2(\mathcal P_\tau)}$, and $\|\nabla_yv\|_{C^1(D_1)}\le C\tau^{-1/2}m$. For meridional $U$ one may take $m=\|e(U)\|_{L^2(\mathcal P_\tau)}$ in both estimates. Tangential derivatives are taken at fixed $t$, and the constant is uniform for sufficiently small $\tau$.
\end{lemma}

\begin{proof}
Denote the physical derivatives after flattening by
\[
L_\alpha=\partial_{y_\alpha}-t(\partial_{y_\alpha}a/a)\partial_t,\qquad L_3=(\tau a)^{-1}\partial_t.
\]
Write $\nabla_\tau$ for the gradient with columns $L_j$ and $e_\tau$ for its symmetric part. Throughout this proof the norms in flattened coordinates use $\mathrm dy\,\mathrm dt$. Set $E=\tau^{-1/2}m$.

\smallskip
\noindent\emph{Korn estimates in the even class.} Cover any fixed inner lateral subdomain by circular columns with lateral radius comparable to $\tau$ and with the full thickness. Use larger columns in $D$ with bounded overlap. After isotropic rescaling by $\tau$, each column is a uniformly smooth graph deformation of a circular cylinder of height comparable to one. The graph changes by $O(\tau)$ in $C^1$ on the rescaled column. To see uniformity of Korn's inequality, subtract the $L^2$ projection onto the physical rigid motions. If uniformity failed, there would be a sequence of such remainders of unit $H^1$ norm and with strain tending to zero. Pullback to a limiting cylinder preserves the vanishing-strain limit and the orthogonality to rigid motions. Rellich compactness gives strong $L^2$ convergence to zero. Korn's second inequality on that cylinder then gives strong $H^1$ convergence to zero, a contradiction.

Averaging the rigid displacement with its midplane reflection preserves the estimate for fields whose horizontal components are even and whose vertical component is odd in $Z$. This leaves only horizontal translations and rotation about the vertical axis. The vertical derivatives $\partial_Z U$ and the two horizontal derivatives of $U_3$ of such a rigid displacement vanish. The column estimates therefore control these derivatives of the field by its strain. This conclusion applies to all fields of the stated parity, with no rotational symmetry assumption on these fields.

For meridional $U$, the in-plane physical gradient of $U'$ is symmetric on each horizontal slice. In this case the column bound, and in the general case the definition of $m$, give $\|\nabla_\tau W\|_{L^2}\le CE$ on any fixed inner lateral subdomain. For use with test fields, let $\Phi(y,t)$ have the same midplane parity and compact lateral support in a fixed inner subdomain. The column argument applied to its physical representative controls $\|L_3\Phi\|_{L^2}+\sum_\alpha\|L_\alpha\Phi_3\|_{L^2}$ by $C\|e_\tau(\Phi)\|_{L^2}$. At each $t$, extend $\Phi'$ by zero laterally and apply planar Korn on $\bR^2$ to control $\|\nabla_y\Phi'\|_{L^2}$ by $\|e_y(\Phi')\|_{L^2}$. The latter symmetric gradient differs from the in-plane block of $e_\tau(\Phi)$ by a term bounded by $C\tau|L_3\Phi|$, because $\partial_{y_\alpha}=L_\alpha+\tau t(\partial_{y_\alpha}a)L_3$. We obtain the uniform estimate
\be\lb{sph:eq:layer-korn}
\|\nabla_\tau \Phi\|_{L^2}\le C\|e_\tau(\Phi)\|_{L^2}
\ee
for compactly supported lateral test fields of this parity. Their values at $t=\pm1$ are unrestricted.

\smallskip
\noindent\emph{Tangential derivatives.} The flattened weak equation has bilinear form
\[
A_\tau(\Phi,G)=\int_{D\times(-1,1)} a(y)\,\bC_0e_\tau(\Phi):e_\tau(G)\,\mathrm dy\,\mathrm dt.
\]
The common factor $\tau$ of the volume element is omitted. Tangential derivatives of the coefficients of $L_j$ yield multiples of $L_3$ with uniformly bounded coefficients. Through order four these commutators differentiate $(\partial_\alpha a)/a$ at most four times, so derivatives of $a$ through order five suffice. On nested inner lateral subdomains we consequently have
\be\lb{sph:eq:layer-tangential}
\sum_{j=0}^4\|\nabla_\tau D_y^jW\|_{L^2}\le CE.
\ee
Differentiating the weak equation $j$ times tangentially gives an equation for $\Phi=D_y^jW$, with the absolute value of its right-hand side on a test $\psi$ bounded by a constant times $\sum_{k<j}\|\nabla_\tau D_y^kW\|_{L^2}\|\nabla_\tau\psi\|_{L^2}$. This follows by the product rule on both factors of the bilinear form. Derivatives may also be implemented by lateral difference quotients on the fixed cylinder.

Test with $\psi=\zeta^2\Phi$, where $\zeta$ is a lateral cutoff. Tangential differentiation preserves the midplane parity, so \eq{sph:eq:layer-korn} applies to $\zeta \Phi$. The identity
\[
A_\tau(\Phi,\zeta^2\Phi)=A_\tau(\zeta \Phi,\zeta \Phi)-\int_{D\times(-1,1)}a\,\bC_0(\Phi\odot\nabla_y\zeta):(\Phi\odot\nabla_y\zeta)\,\rd y\,\rd t
\]
together with Cauchy--Schwarz and absorption gives
\[
\|\nabla_\tau(\zeta \Phi)\|_{L^2}\le C\left(\|\Phi\|_{L^2}+\sum_{k<j}\|\nabla_\tau D_y^kW\|_{L^2}\right).
\]
The $L^2$ norm of $\Phi$ is controlled by the preceding tangential gradient norm. This closes the induction, starting from the estimate for $\nabla_\tau W$ just proved.

\smallskip
\noindent\emph{Recovery on the faces.} Apply \eq{layer:eq:conormal-balance} to the flattened stress $\bC_0e_\tau(W)$. Integrate in $t$, and use the planar $H^2$ embedding and \eq{sph:eq:layer-tangential}. This yields, on a smaller cylinder,
\[
\|\mathcal T\|_{L^\infty}+\|\nabla_y\mathcal T\|_{L^\infty}\le C\tau E.
\]
The same embedding in $y$, together with the one-dimensional $H^1$ bound in $t$, gives $\|\nabla_yW\|_{L^\infty}\le CE$. The same argument bounds $\|\nabla_y^2W\|_{L^\infty}$ by $CE$. For both estimates, the $L^2_tH^2_y$ norm of the required mixed derivative $\partial_t\nabla_y^jW$, $j=1,2$, is at most $C\tau E$ by \eq{sph:eq:layer-tangential}. For $j=1,2$, integrate $\partial_tD_y^jW=\tau aL_3D_y^jW$ from $t=0$ for the odd component, and subtract the thickness average for the even components. The same bounds give
\[
\|D_y^jW_3\|_\infty+\|D_y^j(W'-v)\|_\infty\le C\tau E.
\]
Averaging $D_yW'$ and $D_y^2W'$ also gives $\|\nabla_yv\|_{C^1}\le CE$.

For the isotropic tensor the acoustic matrix in \eq{layer:eq:conormal-recovery} is $\mathsf A_{\bC_0}(n)=\mu|n|^2I+(\lambda+\mu)n\otimes n$, and $\partial_ZU\circ\Phi_\tau=L_3W$. Differentiating this formula in $y$ and using the bounds for $\mathcal T$, $\nabla_yW$, and $\nabla_y^2W$ also bounds $L_3W$ and $\nabla_yL_3W$ by $CE$. In particular, $\sigma$ and its first $y$-derivatives are bounded by $CE$.

As $n=e_3+O(\tau)$, its acoustic matrix differs by $O(\tau)$ from $\operatorname{diag}(\mu,\mu,\lambda+2\mu)$. The term $\bC_0\operatorname{sym}(G_y)n$ is $\lambda\operatorname{div}_yW'\,e_3+O(\tau E)$, since $\nabla_yW_3=O(\tau E)$. Consequently
\[
L_3W=-\eta_{\rm ps}\operatorname{div}_y\overline W'\,e_3+O(\tau E).
\]
The expansion holds after one tangential derivative as well. Indeed $n-e_3=O(\tau)$ in $C_y^1$, and the preceding mixed derivative estimates control the $C_y^1$ errors in $\nabla_yW_3$ and $\nabla_y(W'-v)$. Substituting into $G_y+L_3W\otimes n$ gives the gradient assertion. Taking the $33$ component of the stress gives the second bound since the approximating block has zero normal stress.
\end{proof}

\subsection{The membrane error on gap-scale patches}\label{sec:membrane-error}
We now use the local $C^5$ hypothesis of Theorem~\ref{sph:thm:general-cavity}. Radial symmetry gives $D^3g(0)=0$, so $|D^j(g-|x'|^2)|\le Cr^{4-j}$ for $0\le j\le4$. For the in-plane thickness average $v$ of a reflection-even displacement $u$, define
\be\label{sph:eq:even-membrane-force}
 \mathcal E=H\bC_{\rm ps}e'(v)-M,\qquad
 \operatorname{div}'(H\bC_{\rm ps}e'(v))=\operatorname{div}'\mathcal E .
\ee
The equation follows from $\operatorname{div}'M=0$. Fix a small radius $\bar r$ with $s_\varepsilon=\sqrt\varepsilon<\bar r/8$. For $0<q<2$, set $K=\sup_{G_\varepsilon(\bar r)}d(r)^{2-q}|\nabla u|$.

For $r<\bar r/8$, let $R=\max\{s_\varepsilon,r\}$ and rescale by $(x',z)=(Ry,RZ)$. The rescaled half-thickness is $R a(y)$, where $a(y)=\delta(R|y|)/(2R^2)$. If $r<s_\varepsilon$, use nested disks in $|y|<4$ with inner disk containing $\overline{B'_1}$. Otherwise use nested annuli in $1/4<|y|<4$ with inner annulus containing $|y|=1$. The enlarged patches lie in $G_\varepsilon(\bar r)$. On either patch $a$ is bounded below and has a uniform $C^5$ bound. Indeed, $g(0)=Dg(0)=0$, and the derivatives of $g$ of orders $j\ge2$ occur with the factor $R^{j-2}$.

The enlarged physical patch has volume $O(R^4)$, so the $L^2$ norm of $\nabla u$ there is at most $CKR^q$. For $U(y,Z)=u(Ry,RZ)$ the corresponding norm of $\nabla U$ is $R^{-1/2}$ times as large. Returning to the physical coordinates in Lemma~\ref{sph:lem:even-local} gives
\be\label{sph:eq:even-point-transfer}
 \left|\nabla u-\begin{pmatrix}\nabla'v&0\\0&-\eta_{\rm ps}\operatorname{div}'v\end{pmatrix}\right|
 \le CKd(r)^{q-1}
\ee
throughout the section. The same bound holds after one scaled lateral derivative at fixed flattened height. The lemma also yields $|\nabla'v|+R|D'^2v|\le CKR^{q-2}$.

Let $\overline\sigma'$ be the thickness average of the in-plane stress. The approximating matrix in \eq{sph:eq:even-point-transfer} has in-plane stress $\bC_{\rm ps}e'(v)$. Thus $\overline\sigma'-\bC_{\rm ps}e'(v)$ and its first scaled lateral derivatives are bounded by $CKR^{q-1}$. In the identity
\[
 \mathcal E=(H-\delta)\bC_{\rm ps}e'(v)+\delta(\bC_{\rm ps}e'(v)-\overline\sigma'),
\]
the coefficients and their scaled derivatives have orders $R^4$ and $R^2$, respectively. Hence
\be\label{sph:eq:even-forcing-size}
 |\mathcal E|+d(r)|\nabla'\mathcal E|\le CKd(r)^{q+1},\qquad r<\bar r/8 .
\ee
The constants in these two estimates are independent of $\bar r$. On $\bar r/8\le r\le\bar r$, the error, its first derivative, and the trace of $v$ at $r=\bar r$ are bounded by $C_{\bar r}\|\varphi\|_{C^{2,\gamma}}$ by the anchored Korn inequality and traction regularity on the separated faces used in \eq{geom:eq:collar}.

\subsection{The radial upper bound}
\begin{prop}[Radial concentration]\lb{sph:thm:radial-sharp}
Assume the local regularity of $g$ in Theorem~\ref{sph:thm:general-cavity}. For axisymmetric meridional data $\varphi=\varphi_r e_r+\varphi_3 e_3$ with $\varphi_r$ even and $\varphi_3$ odd in $z$,
\be\nn
|\nabla u_{\rm cav}(x)|\le C\|\varphi\|_{C^{2,\gamma}}d(r)^{-2+q_0}\qquad(x\in G_\varepsilon(r_0)).
\ee
The estimate includes both free faces.
\end{prop}

\begin{proof}
Put $p=q_0-1\in(0,1)$ and use the weighted supremum $K$ from Section~\ref{sec:membrane-error} with $q=q_0$. At each positive gap $K$ is finite by regularity. The average is $v=V(r)e_r$. With $\Lambda=2\mu+2\lambda\mu/(\lambda+2\mu)$, the radial component of \eq{sph:eq:even-membrane-force} reads
\[
 (rHV')'+(\nu H'-H/r)V=F_r,\qquad
 F_r=(r/\Lambda)(\operatorname{div}'\mathcal E)\cdot e_r .
\]
Estimate \eq{sph:eq:even-forcing-size} gives $|F_r|\le CKr d(r)^{p+1}$ in $r<\bar r/8$. The fixed-collar estimates give $|F_r|\le C_{\bar r}\|\varphi\|_{C^{2,\gamma}}$ on $\bar r/8\le r\le\bar r$.

Choose the regular homogeneous solution $Y(r)=r\varepsilon^{(p-1)/2}F(r^2/\varepsilon)$, with $F$ from \eq{sph:eq:radial-model-integral}. Proposition~\ref{sph:prop:radial-model} and its proof give
\[
 c r d(r)^{p-1}\le Y(r)\le C r d(r)^{p-1},\qquad |Y'(r)|\le C d(r)^{p-1}.
\]
Regularity at the axis fixes the integration constant in the Wronskian identity for $f=V/Y$:
\[
 rHY^2 f'(r)=\int_0^r Y(t)F_r(t)\,\rd t .
\]
In the central disk the integrand satisfies
\[
 |Y(t)F_r(t)|\le CKt^2d(t)^{2p}.
\]
As $p>0$, integration over $0<t<\min\{r,\bar r/8\}$ gives at most $CKr^3d(r)^{2p}$. The remaining contribution occurs only where the denominator has a positive lower bound depending on $\bar r$. Using $rHY^2\ge c r^3d(r)^{2p}$ therefore gives $|f'|\le CK+C_{\bar r}\|\varphi\|_{C^{2,\gamma}}$. The collar estimate bounds $|f(\bar r)|$ by $C_{\bar r}\|\varphi\|_{C^{2,\gamma}}$. Integration and the estimates for $Y,Y'$ give
\[
 |V'|+|V/r|\le \big(C_{\bar r}\|\varphi\|_{C^{2,\gamma}}+C\bar rK\big)d(r)^{p-1}.
\]
Substitution in \eq{sph:eq:even-point-transfer} bounds the weighted gradient on $r<\bar r/8$. Combining this with the collar bound gives $K\le C_{\bar r}\|\varphi\|_{C^{2,\gamma}}+C\bar rK$. Choose $\bar r$ so that $C\bar r<1/2$. This proves the assertion up to both free faces.
\end{proof}

\begin{remark}[Equal spherical cavities]
For two balls of radius $a$, in the original coordinates $\kappa=1/a$ and the gap has the exact remainder
\[
 \delta_\varepsilon(r)-\varepsilon-r^2/a
 =\frac{r^4}{a\bigl(a+\sqrt{a^2-r^2}\bigr)^2}.
\]
Differentiation gives the geometric bounds of Section~\ref{sec:membrane-error}.
\end{remark}

\section{Angular competition in the membrane equation}

Angular projections rotate both the argument and the vector components. In the polar frame they are the usual Fourier projections, and they commute with the three-dimensional equation and with midplane reflection. We use the common error estimate \eq{sph:eq:even-forcing-size} first in angular order two and then on the complement of the two leading modes.

\subsection{Angular exponents and the order-two mode}
\begin{lemma}[Membrane spectrum]\lb{sph:lem:membrane-spectrum}
The angular modes of $\operatorname{div}'(r^2\bC_{\rm ps}e'(v))$ have powers $r^{-1+Q}$. For angular order $m\ge0$, the indicial equation is
\be\lb{sph:eq:membrane-all-roots}
Q^4+(-2m^2+2\nu-6)Q^2+m^4+(6\nu-2)m^2+8-8\nu=0.
\ee
The least positive root corresponding to a nonrigid solution over all angular modes is $q_*=\min\{q_0,q_2\}$, where
\[
q_0=\sqrt{2-2\nu},\qquad q_2=\sqrt{7-\nu-\sqrt{\nu^2-30\nu+33}}.
\]
Among the modes remaining after removing the radial component of order zero and all of order two, all positive nonrigid roots are strictly larger than $q_*$.
\end{lemma}

\begin{proof}
For a displacement $r^se^{\mathrm im\vartheta}(Ae_r+Be_\vartheta)$, divide the plane stress components by $\Lambda r^{s-1}e^{\mathrm im\vartheta}$, where $\Lambda=2\mu+2\lambda\mu/(\lambda+2\mu)$. Their coefficients are
\[
S_{rr}=(s+\nu)A+\mathrm im\nu B,\qquad S_{\vartheta\vartheta}=(1+\nu s)A+\mathrm imB,\qquad S_{r\vartheta}=\tfrac12(1-\nu)(\mathrm imA+(s-1)B).
\]
Force balance gives the two equations
\[
(s+2)S_{rr}-S_{\vartheta\vartheta}+\mathrm imS_{r\vartheta}=0,\qquad(s+3)S_{r\vartheta}+\mathrm imS_{\vartheta\vartheta}=0.
\]
Substitution with $Q=s+1$ gives the determinant $(1-\nu)/2$ times the polynomial in \eq{sph:eq:membrane-all-roots}. For $m=0$ the positive roots are $q_0$ in the radial component and $2$ in the azimuthal component. The latter gives a rigid rotation. For $m=1$ the positive root $1$ gives a translation, while the other is $\sqrt{7-2\nu}$.

For $m\ge2$ the smaller squared root is
\[
m^2+3-\nu-\sqrt{8(1-\nu)m^2+(1+\nu)^2}.
\]
It is positive and strictly increasing as a function of $m^2$ on $[4,\infty)$. Indeed, at $4$ it equals $q_2^2>0$, and the derivative is $1-4(1-\nu)/\sqrt{8(1-\nu)m^2+(1+\nu)^2}>0$ for $-1<\nu<1/2$.
\end{proof}

\begin{lemma}[The finite-gap membrane equation in angular order two]\lb{sph:lem:mode-two-inverse}
On the disk $r<b$, let $v\in H^1$ be a weak solution of \eq{sph:eq:even-membrane-force} of angular order two, with $s=\sqrt\varepsilon$ in $(0,b/32)$. Suppose $\mathcal E$ has the same order, in the polar tensor frame, and
\[
|\mathcal E|+d|\nabla'\mathcal E|\le Ld^{q_2+1}\quad(r<b/8),\qquad\|\mathcal E\|_{L^\infty(b/8<r<b)}+\|v\|_{C^1(r=b)}\le A.
\]
Then for $r<b/16$,
\[
|\nabla'v(x')|\le(C_bA+CbL)d(r)^{-2+q_2}.
\]
The coefficient of $bL$ is independent of $b$ and $s$.
\end{lemma}

\begin{proof}
First consider a regular homogeneous solution $Y$ of order two on a disk $r<t$, with $s\le t$. Put $\mathsf E_t(Y)^2=\int_{r<t}H|e(Y)|^2\,\mathrm dx'$. We claim that
\be\lb{sph:eq:mode-two-propagation}
|\nabla'Y(x')|\le C\mathsf E_t(Y)t^{-q_2}d(r)^{q_2-2}\qquad(r<t/2).
\ee
To see this, first let $s=1$ by the horizontal rescaling. The coefficient then becomes $1+r^2$. The complex regular solution space for the angular frequency $2$ has dimension two. Its trace on the circle $r=1$ determines it on the disk by coercivity, and every trace in that frequency has a smooth solution on the disk which is continued by the ordinary differential equation.

At infinity this differential equation, in $\log r$, is a constant system plus an $O(r^{-2})$ perturbation. Its four distinct exponents are $-1\pm q_2$ and $-1\pm\widehat q_2$, where $\widehat q_2>q_2$ is the larger positive root of \eq{sph:eq:membrane-all-roots}. In the variable $1/r$ this system is regular singular, with coefficients analytic in $r^{-2}$ after the Euler derivative $r\partial_r$ is used. The four indicial roots are simple. The Frobenius solutions therefore have leading terms $r^\sigma$ times the corresponding eigenvector, also after one Euler derivative. At a resonance a logarithm can multiply a lower power without changing these leading terms. Equivalently one may use asymptotic integration \cite{Eastham85}.

A regular solution whose two positive-root coefficients vanish is zero. Indeed, its traction work at radius $R$ is $O(R^{-2q_2})$ as $R\to\infty$. Testing with the solution on that disk gives zero strain on the plane, and the angular order excludes rigid motions. Thus the coefficients of the two positive roots give an isomorphism on the regular solution space. Here the sign refers to $Q$ in $r^{-1+Q}$. Even when $q_2<1$, the annular strain energy of that branch grows like $R^{2q_2}$ at infinity.

Choose regular solutions $Y_1,Y_2$ with positive-root coefficients $(1,0)$ and $(0,1)$, respectively. Their gradients satisfy $|\nabla Y_j(r)|\le C(1+r)^{q_j'-2}$, where $q_1'=q_2$ and $q_2'=\widehat q_2$. On any annulus $T/2<r<T$, $T\ge1$, the weighted strain norm of $c_1Y_1+c_2Y_2$ bounds $\sum_{j=1}^2|c_j|T^{q_j'}$ from below up to a fixed factor. For large $T$ this follows by scaling the two strains to $1/2<r<1$, where the limiting strain profiles are linearly independent, since their powers differ and neither is zero. For $T$ in a compact subinterval of $[1,\infty)$ it follows from independence on every annulus and compactness. In the physical coordinates the basis is $Y_j(x'/s)$. The energy norm on the rescaled annulus is multiplied by $s$, and the gradient by $s^{-1}$. With $T=t/s$ this gives
\[
|\nabla'Y(x')|\le C\mathsf E_t(Y)\sum_{j=1}^2t^{-q_j'}d(r)^{q_j'-2}.
\]
Since $d(r)\le3t/2$ on $r<t/2$, this proves \eq{sph:eq:mode-two-propagation}.

Separate the solution with the given lateral trace from the zero-trace response to $\mathcal E$. Choose a smooth radial cutoff $\xi$ equal to one on $r\le3b/32$ and zero on $r\ge b/8$, with $|\nabla\xi|\le C/b$. Split the source of the zero-trace problem as
\[
\xi\mathcal E,\qquad(1-\xi)\mathbf 1_{\{r<b/8\}}\mathcal E,\qquad\mathbf 1_{\{r\ge b/8\}}\mathcal E.
\]
All three tensors are admissible divergence data in the energy formulation. Only the first will be differentiated. Partition it smoothly into a central disk $r<4s$ and dyadic annuli $t/2<r<2t$, $s\le t\le b$, with bounded overlap and cutoff derivatives bounded by $C/t$.

For an inner tensor piece at scale $t$, the energy identity gives
\[
\mathsf E_b(v_t)\le C\|\mathcal E_t/\sqrt H\|_{L^2}\le CLt^{q_2+1}.
\]
The central piece has the same bound with $t=s$.

Consider a point with $d=d(r)$ and $r<b/16$. If $t$ is larger than a fixed multiple of $d$, apply \eq{sph:eq:mode-two-propagation} on a concentric disk of radius comparable to $t$ disjoint from the support. The resulting gradient bound is $CLtd^{q_2-2}$. If $t$ is much smaller than $d$, a ball or annulus at scale $d$ separates the point from the forcing. Rescaled interior estimates bound the gradient by $Cd^{-2}\mathsf E_b(v_t)$. At $t$ comparable to $d$, the same local estimate with forcing adds
\[
Cd^{-2}\bigl(\|\mathcal E_t\|_{L^\infty}+d\|\nabla'\mathcal E_t\|_{L^\infty}\bigr).
\]
Here one rescales by $d$ on a fixed annulus or on the central disk, where $H/d^2$ is uniformly elliptic. The $L^2$ term in the interior gradient estimate is reduced to the strain norm by planar Korn after subtracting a rigid displacement. Equivalently one can apply Korn on the annulus or disk to the order-two field, which is orthogonal to all planar rigid motions. The forcing term follows from the interior $C^{1,\alpha}$ estimate for a divergence datum with a bounded first derivative. The sum of $t$ for the outer dyadic scales is at most $Cb$. The inner scales give $CLd^{-2}\sum_{t\lesssim d}t^{q_2+1}\le CLd^{q_2-1}$. These bounds also include the finitely many comparable scales.

It remains to bound the last two source pieces and the lateral trace. Denote their zero-trace responses by $v_{\rm mid}$ and $v_{\rm out}$, and the homogeneous solution with the prescribed trace by $v_{\partial}$. Testing the first two equations and extending the lateral trace in a fixed collar gives
\[
\mathsf E_b(v_{\rm mid})\le CLb^{q_2+1},\qquad\mathsf E_b(v_{\rm out})+\mathsf E_b(v_{\partial})\le C_bA.
\]
Indeed $H\ge cb^2$ on both source supports, and the middle tensor has size at most $CLb^{q_2+1}$ on a set of area $O(b^2)$. All three fields are homogeneous on $r<3b/32$. For $r<b/32$ apply \eq{sph:eq:mode-two-propagation} with $t=3b/32$. For $b/32\le r<b/16$, rescale an annulus separated from the sources by a fixed multiple of $b$. Interior regularity and Korn's inequality in order two bound each gradient there by $Cb^{-2}$ times its energy norm. Since $d(r)$ is comparable to $b$ on this annulus, the middle piece contributes at most $CbLd(r)^{q_2-2}$ in both regions. The lateral and outer terms contribute at most $C_bAd(r)^{q_2-2}$. Together with the dyadic sum this proves the claim.
\end{proof}

\subsection{The full membrane estimate}
\begin{prop}[The even endpoint estimate]\lb{sph:thm:even-sharp}
For arbitrary outer data, the reflection-even part of the cavity solution satisfies
\[
|\nabla u_{\rm cav}^+(x)|\le C\|\varphi\|_{C^{2,\gamma}}d(r)^{q_*-2}\qquad(x\in G_\varepsilon(r_0)),
\]
including both free faces.
\end{prop}

\begin{proof}
The radial component has the stated estimate with $q=q_0$ by Proposition~\ref{sph:thm:radial-sharp}. For the component of order two, apply \eq{sph:eq:even-point-transfer} and \eq{sph:eq:even-forcing-size} with $q=q_2$. Lemma~\ref{sph:lem:mode-two-inverse} and the fixed-collar bound give
\[
K\le C_b\|\varphi\|_{C^{2,\gamma}}+CbK
\]
for a fixed sufficiently small $b$, the radius denoted by $\bar r$ in Section~\ref{sec:membrane-error}, with $C$ independent of $b$. Take $Cb<1/2$. This bounds the order-two component at its endpoint, including both traces.

It remains to control the projection with these components removed. We prove the asserted estimate with $q=q_*$ by contradiction, normalizing the norm of the outer data to be at most one. If $K$ were unbounded, we could choose gaps tending to zero and points with $r_n\to0$ at which $d(r_n)^{2-q_*}|\nabla u|\ge K/2$, by the fixed-collar estimate. Rescale isotropically by $R_n=\max\{s_{\varepsilon_n},r_n\}$ and normalize the rescaled gradient by $R_n^{2-q_*}/K$. After passage to a subsequence, $s_{\varepsilon_n}/R_n\to a\in[0,1]$.

On each compact planar subdomain, avoiding $0$ when $a=0$, the local membrane approximation bounds the gradients of the in-plane averages in $C^1$. Subtract horizontal translations, if necessary, to obtain a limiting average $V$; in the case $a=0$ an exhaustion starting on a fixed annulus gives the same normalization on every connected compact subdomain. A constant horizontal translation has zero strain and affects only the rigid coefficient of order one. By the vanishing aspect ratio, the local error in the gradient approximation tends to zero. The selected points stay in a fixed inner disk if $a>0$, or on $r=1$ if $a=0$. Thus $\nabla'V\ne0$ at a limit point. Passing to the limit in the resultant equation shows
\[
\operatorname{div}'((a^2+r^2)\bC_{\rm ps}e(V))=0,\qquad|\nabla'V|\le C(a+r)^{q_*-2}.
\]
The excluded angular components are still absent. The bound and the equation hold on the entire plane if $a>0$ and on the punctured plane if $a=0$.

To rule out this limit, project the membrane equation onto one angular mode at a time. When $a=0$ the equations are Euler equations. The bound at infinity rules out all nonrigid positive-root powers by Lemma~\ref{sph:lem:membrane-spectrum}. The bound at zero rules out all negative-root powers, since their nonzero gradients are proportional to $r^{-2-|Q|}$. A translation has zero gradient, and rigid rotation is excluded by the bound at infinity.

When $a>0$ each mode is regular at zero. Its asymptotic basis at infinity consists of the same powers, by asymptotic integration as in Lemma~\ref{sph:lem:mode-two-inverse}. The nonrigid positive-root coefficients again vanish. A possible translation in order one is an exact constant solution and may be subtracted. All remaining terms decay as negative-root powers, so the boundary traction work for that mode tends to zero. The energy identity on disks then shows that its strain is zero. Consequently all the Fourier components have zero gradient after removal of the translation. This contradicts the selected nonzero limit. The argument uses only the even parity, and works at the two matrix-side traces by Lemma~\ref{sph:lem:even-local}.

\end{proof}

\section{Rotation and the bending spectrum}\label{sec:bending}
The reflection-odd problem has a bending limit on rescaled layers. The local estimates allow laterally varying tensors, while the comparison of bending and membrane exponents uses isotropy and radial thickness.

\subsection{Strain, rotation and the bending limit}

Let $D_1\Subset D$ be bounded connected $C^2$ planar domains. Let $a$ be positive on $\overline D$, and let $\bC(y)$ be a self-adjoint uniformly positive elasticity tensor on symmetric $3\times3$ matrices. Assume a fixed positive lower bound for $a$ and fixed ellipticity bounds for $\bC$. The tensor is independent of the thickness coordinate and may be anisotropic.

\begin{theorem}[Uniform strain and rotation estimates]\lb{layer:thm:uniform-bending}
Suppose $a$ and $\bC$ have fixed $C^5$ bounds. For sufficiently small $h>0$, let $U$ solve $\operatorname{div}(\bC(y)e(U))=0$ in
\[
\mathcal P_h(D)=\{(y,Z):y\in D,\ |Z|<ha(y)\},
\]
with zero traction on both faces. Suppose $U'$ is odd and $U_3$ is even in $Z$. There is a rigid displacement $R=(-Zc,c\cdot y+c_0)$ such that
\be\lb{layer:eq:uniform-bounds}
h^{1/2}\|e(U)\|_{L^\infty(\mathcal P_h(D_1))}
+h^{3/2}\|\nabla(U-R)\|_{L^\infty(\mathcal P_h(D_1))}
\le C\|e(U)\|_{L^2(\mathcal P_h(D))}.
\ee
The estimates hold through both free faces. Their constants depend only on the lateral domains and the indicated geometric and material bounds. Both powers are optimal, already for a flat layer with a constant isotropic tensor. No lateral boundary condition is required.
\end{theorem}

For reflection-invariant materials the parity condition follows by projection. Primes denote lateral components or derivatives, and $e'$ denotes the lateral symmetric gradient.

For $p\in\bR^2$ define the skew matrix
\[
\mathcal R(p)=\begin{pmatrix}0_{2\times2}&-p\\p^{\mathsf T}&0\end{pmatrix}.
\]

\begin{lemma}[Local bending regularity]\lb{sph:lem:odd-local}
In the setting of Theorem~\ref{layer:thm:uniform-bending}, assume $a,\bC$ have fixed $C^{N+2}$ bounds for some integer $N\ge3$. Set
\[
K_h=h^{-1/2}\|\nabla U\|_{L^2(\mathcal P_h(D))},\qquad W(y,t)=U(y,ha(y)t),\qquad w_h(y)=\tfrac12\int_{-1}^1W_3(y,t)\,\mathrm dt.
\]
After subtracting a vertical constant, the following estimates hold on the closed layer over $D_1$:
\be\lb{sph:eq:odd-local-pointwise}
|e(U)|+|\nabla U(y,ha(y)t)-\mathcal R(\nabla_yw_h(y))|\le ChK_h,\qquad\|w_h\|_{H^{N+1}(D_1)}\le C_NK_h.
\ee
Thus $\|w_h\|_{C^{2,\alpha}(D_1)}\le C_\alpha K_h$ for $N=3$ and $0<\alpha<1$, and $\|w_h\|_{C^3(D_1)}\le CK_h$ for $N\ge4$. The constants in \eq{sph:eq:odd-local-pointwise} depend only on $N$, the lateral domains and the indicated geometric and material bounds. If $\|\nabla U\|_{L^\infty}\le K$, all these bounds hold with $CK$ in place of $CK_h$.
\end{lemma}

\begin{proof}[Proof of Theorem~\ref{layer:thm:uniform-bending}]
Use the thin-layer Korn inequality \eq{sph:eq:odd-free-korn} in Appendix~\ref{sec:layer-proof} to choose $R$, and apply Lemma~\ref{sph:lem:odd-local} to $U-R$. Its normalized gradient norm is at most $Ch^{-3/2}\|e(U)\|_{L^2}$. The strain bound in~\eq{sph:eq:odd-local-pointwise} gives the first estimate in~\eq{layer:eq:uniform-bounds}. The gradient approximation and the bound for $\nabla_y w_h$ give the second.

For sharpness take $a=1$ and
\[
U(y,Z)=\left(-Zy_1,\ 0,\ \tfrac12y_1^2+\tfrac{\eta_{\rm ps}}{2}Z^2\right),\qquad\eta_{\rm ps}=\frac{\lambda}{\lambda+2\mu}.
\]
For the constant isotropic tensor $\bC_0$, this is an exact homogeneous Lam\'e solution with zero traction on both flat faces. Its only nonzero strains are $e_{11}=-Z$ and $e_{33}=\eta_{\rm ps} Z$. Thus its strain supremum is bounded above and below by positive multiples of $h$, and its strain $L^2$ norm by multiples of $h^{3/2}$. On a fixed inner lateral domain its gradient has distance at least $c>0$ from every constant skew matrix, because $\partial_ZU_1=-y_1$ varies over a fixed interval. These facts prove optimality of both powers.
\end{proof}

\begin{corollary}[Strong rotation compactness and the bending equation]\lb{layer:cor:compactness}
Let $h\to0$, and let $a_h,\bC_h$ satisfy the bounds in Lemma~\ref{sph:lem:odd-local} uniformly, with $N=3$. Suppose $a_h\to a$ in $C^1_{\rm loc}(D)$ and $\bC_h\to\bC$ locally uniformly. Let $U_h$ have the parity and boundary conditions of that lemma, solve the corresponding equation, and satisfy $h^{-1/2}\|\nabla U_h\|_{L^2(\mathcal P_h(D))}\le C$. After subtraction of vertical constants and passage to a subsequence, the averaged normal displacements converge in $C^2_{\rm loc}(D)$ to $w$, and
\[
\nabla U_h(y,ha_h(y)t)\longrightarrow\mathcal R(\nabla_yw(y))\quad\text{locally uniformly on }D\times[-1,1].
\]
Moreover,
\be\lb{layer:eq:general-bending-limit}
\operatorname{div}'\operatorname{div}'\bigl(a^3\bC_{\rm ps}(y)D_y^2w\bigr)=0\quad\text{in }D.
\ee
Here $\bC_{\rm ps}$ is the plane stress Schur complement of $\bC$, characterized, for symmetric $2\times2$ matrices $E$, by
\[
\bC_{\rm ps}(y)E:E=\min_{p\in\bR^3}\bC(y)\bigl(\iota E+\operatorname{sym}(p\otimes e_3)\bigr):\bigl(\iota E+\operatorname{sym}(p\otimes e_3)\bigr),
\qquad\iota E=\begin{pmatrix}E&0\\0&0\end{pmatrix}.
\]
\end{corollary}

\begin{proof}
The normalized gradient bound and Poincar\'e's inequality fix a vertical constant on each connected lateral subdomain. Lemma~\ref{sph:lem:odd-local} bounds the averaged normal displacement in $H^4$ locally on the two-dimensional base, hence in $C^{2,\alpha}$ for each $0<\alpha<1$, and approximates the full gradient with error $O(h)$, uniformly through both faces. Compactness on nested subdomains and a diagonal extraction give the asserted convergences.

Write $W_h(y,t)=U_h(y,ha_h(y)t)$. The lemma also bounds $e(U_h)/h$ locally in the flattened $L^2$ norm. Planar Korn on each $t$-slice, applied after a lateral cutoff, bounds $W_h'/h$ in $H^1$. Indeed, its planar strain is the physical planar strain divided by $h$, plus $t\nabla a_h\odot\partial_ZU_h'$; its $t$-derivative is $a_h\partial_ZU_h'$. All these terms are bounded locally. Oddness and the transverse shear identity give
\[
W_h'/h\rightharpoonup-at\nabla_yw\quad\text{in }H^1_{\rm loc},\qquad h^{-1}(e(U_h)_{\alpha\beta})\rightharpoonup-at\partial_{\alpha\beta}w.
\]
In the second limit the derivatives of $a$ cancel when the flattened derivatives are converted to physical derivatives.

Set $\sigma_h=\bC_h(y)e(U_h)$. For an arbitrary smooth vector $f(y,t)$, compactly supported laterally and smooth up to $t=\pm1$, test equilibrium with $h^2f(y,Z/(ha_h(y)))$. Divide by $h^3$. The terms containing lateral derivatives are $O(h)$, while the remaining term is $\int(\sigma_{h,\cdot3}/h)\cdot\partial_tf$. Since $\partial_tf$ is arbitrary, the weak limit of $\sigma_{h,\cdot3}/h$ is zero. Minimization in the definition of $\bC_{\rm ps}$ is equivalent to zero third traction. It is unique by ellipticity. Thus the planar block of $\sigma_h/h$ converges weakly to $-at\bC_{\rm ps}(y)D_y^2w$.

Finally test with $(-Z\nabla_y\psi,\psi)$, where $\psi\in C_c^\infty(D)$. Its only nonzero strains are the planar ones $-ZD_y^2\psi$. Flattening, dividing by $h^3$, and canceling $\int_{-1}^1t^2\,\mathrm dt$ proves \eq{layer:eq:general-bending-limit}.
\end{proof}

\subsection{The bending exponents and the limiting equation}

Corollary~\ref{layer:cor:compactness} identifies the bending operator of the rescaled layers. For the constant isotropic matrix its plane stress tensor is the tensor defined in Section~\ref{sec:membrane}. In rescalings of the quadratic neck the limiting $a(y)$ is a constant multiple of $a_0^2+|y|^2$, where $a_0\ge0$ is the limiting ratio of the gap scale to the rescaling radius. In the contact case $a_0=0$ this gives
\be\lb{sph:eq:bending-contact}
\operatorname{div}'\operatorname{div}'\bigl(r^6\bC_{\rm ps}D'^2w\bigr)=0.
\ee
Comparing its homogeneous solutions with the membrane displacements reduces the question of bending concentration to the locations of the roots of two polynomials.

\begin{prop}[Separation of membrane and bending exponents]\lb{sph:prop:bending-spectrum}
For an angular mode $w=r^{Q-2}e^{\mathrm im\vartheta}$ in \eq{sph:eq:bending-contact}, the indicial polynomial is
\be\lb{sph:eq:bending-polynomial}
Q^4+(-2m^2+6\nu-14)Q^2+m^4+(22-30\nu)m^2+40-24\nu.
\ee
If $\operatorname{Re}Q>0$ and the corresponding mode is nonaffine with finite bending energy, then $\operatorname{Re}Q>1+q_*$. This also holds for logarithmic modes associated with repeated positive roots.
\end{prop}

\begin{proof}
Write $s=Q-2$. Up to a common factor, the polar Hessian components of $r^se^{\mathrm im\vartheta}$ are $a=s(s-1)$, $b=s-m^2$, $c=\mathrm im(s-1)$. After applying $\bC_{\rm ps}$, dividing by $2\mu+2\lambda\mu/(\lambda+2\mu)$, and multiplying by $r^6$, the three components of the moment tensor are $a+\nu b$, $b+\nu a$, $(1-\nu)c$, with radial power $s+4$. For a moment tensor with this power and components denoted by $M_{rr},M_{\vartheta\vartheta},M_{r\vartheta}$, two polar divergences give the coefficient
\[
(s+4)(s+5)M_{rr}-(s+4+m^2)M_{\vartheta\vartheta}+\mathrm im(2s+10)M_{r\vartheta}.
\]
Substitution gives \eq{sph:eq:bending-polynomial}. An individual nonaffine mode has finite energy near the origin if $\operatorname{Re}Q>0$, since its radial bending energy has density comparable, after angular integration, to $r^{2\operatorname{Re}Q-1}$.

For $m=0$, the positive roots are $2$ and $\sqrt{10-6\nu}$. For $m=1$, they are $3$ and $\sqrt{7-6\nu}$. The roots $2$ at $m=0$ and $3$ at $m=1$ give affine functions. The nonaffine roots in these two orders are larger than $1+q_0$, since
\[
(7-6\nu)-(1+q_0)^2=2q_0(q_0-1)>0.
\]
At the repeated root $3$ for $m=1,\nu=-1/3$, the logarithmic branch also has $3>1+q_0$.

For $m\ge2$ put $t=m^2\ge4$. The roots in $Q^2$ are
\[
t+7-3\nu\ \pm\sqrt{8(3\nu-1)t+9(1-\nu)^2}.
\]
If the discriminant is negative, the real part of each positive-real-part square root is at least $\sqrt{t+7-3\nu}$. If the discriminant is nonnegative, it is at most $4t+9/4$, because it is increasing in $\nu\in(-1,1/2)$ for $t\ge4$. In both cases the squared real part is at least
\[
t+\frac{11}{2}-\sqrt{4t+\frac94}\ \ge\ \frac{19-\sqrt{73}}2.
\]
The last inequality follows by monotonicity in $t\ge4$. On the other hand $q_*\le\sqrt{8/5}$, since $q_0$ is decreasing in $\nu$, $q_2$ is increasing, and they agree at $\nu=1/5$. The inequality $(19-\sqrt{73})/2>(1+\sqrt{8/5})^2$ proves the claim, including repeated roots. For example, $\sqrt{73}<43/5$ and $\sqrt{8/5}<19/15$ bound the difference from below by $14/225$.
\end{proof}

The growing slope of a nonaffine bending mode has radial power $\operatorname{Re}Q-3>q_*-2$. The following Liouville statement supplies the growth exclusion needed for the three-dimensional blow-up argument.

\begin{lemma}[A bending Liouville theorem]\lb{sph:lem:bending-liouville}
Let $a_0\ge0$ and $0<q<2$. Suppose every nonaffine positive-real-part root in \eq{sph:eq:bending-polynomial} has real part greater than $1+q$. Let $w$ be a weak $H^2_{\rm loc}$ solution of
\[
\operatorname{div}'\operatorname{div}'\bigl((a_0^2+r^2)^3\bC_{\rm ps}D'^2w\bigr)=0.
\]
The domain is $\bR^2$ if $a_0>0$, and $\bR^2\setminus\{0\}$ if $a_0=0$. Suppose
\[
|\nabla'w|\le C(a_0+r)^{q-2}.
\]
Then $w$ is constant.
\end{lemma}

\begin{proof}
Project onto angular frequencies. For $a_0=0$ the equation at each frequency is the Euler equation with polynomial \eq{sph:eq:bending-polynomial}. Its nonaffine positive-real-part roots have real part greater than $1+q$, by hypothesis. Their gradient powers therefore exceed $q-2$ and are excluded by the bound at infinity. All negative-real-part roots have gradient powers less than $-3$, and are excluded by the bound at zero, since $q>0$. This argument also excludes the corresponding logarithmic terms at repeated roots. The affine modes have constant gradients, which vanish by the bound at infinity. Only the constant remains.

For $a_0>0$, scaling reduces to $a_0=1$. Each angular equation is a regular fourth-order equation away from zero, with a regular solution at zero. In the variable $\log r$ its first-order system tends to the contact system with an $O(r^{-2})$ coefficient error. Asymptotic integration \cite{Eastham85}, or its exponential dichotomy formulation on a line separated from the characteristic roots, gives the same positive and negative exponents at infinity. More explicitly, the coefficient matrix is analytic in $r^{-2}$ at infinity in Euler variables. In a Frobenius basis a leading term is a power of $r$ times a polynomial in $\log r$; resonant corrections of lower real order preserve the leading growth.

Equal real parts of conjugate roots cause no cancellation in the growth test: after division by the largest real power and logarithmic degree, the leading angular gradient is a finite sum of exponentials in $\log r$. Its mean square on long intervals is positive if any distinct-frequency coefficient is nonzero. This also applies to the logarithmic companion of the affine mode at $m=1$, $\nu=-1/3$. The constant mode is treated separately below. Thus all nonaffine positive-root coefficients vanish by the gradient bound. The exact constant solution in order zero can be subtracted. The exact affine solutions in order one have nonzero constant gradients and are excluded. Thus, after subtraction of the constant, each regular angular solution satisfying the growth bound lies in the subspace of negative roots at infinity.

For one angular frequency, let $B=(1+r^2)^3\bC_{\rm ps}D'^2w$. After the constant is removed if necessary, integration by parts on the disk gives
\[
\int_{r<R}(1+r^2)^3\bC_{\rm ps}D'^2w:D'^2\overline w=\int_{r=R}\bigl((Bn)\cdot\nabla'\overline w-(\operatorname{div}'B\cdot n)\overline w\bigr).
\]
The right-hand side tends to zero. Indeed a negative root $-Q$, $\operatorname{Re}Q>0$, gives boundary work of order $R^{-2\operatorname{Re}Q}$, with a fixed logarithmic factor when needed. The solution is smooth at zero. Positivity of the plate tensor implies $D'^2w=0$ in that frequency. Its gradient vanishes by the bound at infinity. Applying the argument to every Fourier coefficient proves the lemma.
\end{proof}

\subsection{From bending compactness to the full gradient bound}
For the reflection-odd displacement, the quadratic gap profile enters the gradient bound through a Liouville property of the limiting bending equation.

\begin{prop}[A Liouville criterion for quadratic necks]\lb{sph:prop:odd-neck-criterion}
Let $g\in C^5(B'_{r_0})$, with $g(0)=0$, $D'g(0)=0$, and positive definite quadratic part $Q(x')=\tfrac12D'^2g(0)[x',x']$. Decrease $r_0$, if necessary, so that $g$ is positive away from zero, and set
\[
\delta_\varepsilon(x')=\varepsilon+g(x'),\qquad\mathcal N_\varepsilon(r)=\{(x',z):|x'|<r,\ |z|<\delta_\varepsilon(x')/2\}.
\]
Let $\bC(x')$ be a self-adjoint uniformly positive elasticity tensor with fixed $C^5$ bounds. Write $\bC_{\rm ps}(0)$ for the plane stress Schur complement of $\bC(0)$. Fix $1<q<2$. Assume the following Liouville property for every $a\ge0$: a weak $H^2_{\rm loc}$ solution of
\[
\operatorname{div}'\operatorname{div}'\bigl((a^2+Q)^3\bC_{\rm ps}(0)D'^2w\bigr)=0,\qquad|\nabla'w(y)|\le C_w(a+|y|)^{q-2},
\]
on $\bR^2$ when $a>0$, and on $\bR^2\setminus\{0\}$ when $a=0$, is constant. Then every weak solution $u_\varepsilon$ of $\operatorname{div}(\bC(x')e(u_\varepsilon))=0$ in $\mathcal N_\varepsilon(r_0)$, with zero face traction, $u_\varepsilon'$ odd and $u_{\varepsilon,3}$ even in $z$, and collar bound
\[
\sup_{\mathcal N_\varepsilon(r_0)\setminus\mathcal N_\varepsilon(r_0/2)}|\nabla u_\varepsilon|\le 1,
\]
satisfies, for all sufficiently small positive gaps,
\[
|\nabla u_\varepsilon(x)|\le C(\sqrt\varepsilon+|x'|)^{q-2},\qquad|e(u_\varepsilon)(x)|\le C(\sqrt\varepsilon+|x'|)^{q-1}\quad(x\in\mathcal N_\varepsilon(r_0/4)).
\]
For an odd solution with collar gradient supremum $M$, the right-hand sides are multiplied by $M$. The estimates include both faces. The constant is independent of $\varepsilon,u_\varepsilon$.
\end{prop}

\begin{proof}
At each positive gap, regularity up to a free face makes the weighted gradient supremum finite on a smaller closed neck. On the fixed collar, subtract a vertical constant. Parity bounds $u_\varepsilon'$ by the collar gradient times the thickness, and integration along paths in the collar bounds $u_{\varepsilon,3}$ by $C$. A cutoff of this field supported in the collar supplies a lateral Dirichlet lift with energy at most $C$. The variational principle therefore gives the same strain-energy bound in the whole neck.

For gaps bounded away from zero, ordinary Korn and boundary regularity on this nondegenerate family give the required gradient bound. On every truncated neck $|x'|>r_1>0$, Korn anchored at the outer collar has a constant uniform as $\varepsilon\to0$: the truncated domains form a smooth nondegenerate family and retain a fixed lateral Dirichlet portion. The energy bound and local regularity therefore give $|\nabla u_\varepsilon|\le C_{r_1}$ away from the contact point, including the free faces.

Suppose the gradient estimate fails along gaps tending to zero. Put
\[
K_\varepsilon=\sup_{\mathcal N_\varepsilon(3r_0/4)}(\sqrt\varepsilon+|x'|)^{2-q}|\nabla u_\varepsilon(x)|.
\]
The collar bound forces any sequence of points where the weighted gradient is at least $K_\varepsilon/2$, with $K_\varepsilon\to\infty$, to approach the contact point. Write these points as $x_n=(x_n',z_n)$, and set
\[
R_n=\max\{\sqrt{\varepsilon_n},|x_n'|\},\quad a_n=\sqrt{\varepsilon_n}/R_n,\quad U_n(y,Z)=\frac{R_n^{1-q}}{K_{\varepsilon_n}}u_{\varepsilon_n}(R_ny,R_nZ).
\]
After a subsequence, $a_n\to a\in[0,1]$. The half-thickness of the rescaled layer is $R_nA_n(y)$, where
\[
A_n(y)=\frac{\varepsilon_n+g(R_ny)}{2R_n^2}\longrightarrow\tfrac12(a^2+Q(y)).
\]
On compact lateral subdomains, excluding zero if $a=0$, these functions have uniform positive lower bounds and uniform $C^5$ bounds. Their convergence is locally $C^1$. The tensors $\bC(R_ny)$ have uniform $C^5$ bounds and converge locally uniformly to $\bC(0)$. Moreover
\[
|\nabla U_n(y,Z)|\le C(a_n+|y|)^{q-2}.
\]

Consequently the normalized gradient norms in Corollary~\ref{layer:cor:compactness} are locally bounded. Choose vertical constants on one reference subdomain and use overlapping subdomains to keep this normalization consistent. That corollary and a diagonal extraction give a limit $w$ on the plane if $a>0$, or on the punctured plane if $a=0$, with
\[
\nabla U_n(y,R_nA_n(y)t)\longrightarrow\mathcal R(\nabla_yw(y))
\]
locally uniformly through the full thickness. Its equation is the one in the Liouville hypothesis; the factor $1/8$ from $A_n^3$ cancels. It inherits the stated gradient growth bound.

The selected points lie in a compact lateral subdomain: when $a=0$, their lateral radius is one; when $a>0$, a fixed disk suffices. Their normalized gradients are bounded below by a positive constant because $(\sqrt{\varepsilon_n}+|x_n'|)/R_n\in[1,2]$. The convergence through $t\in[-1,1]$ therefore makes $\nabla'w$ nonzero at a limit point, contradicting the Liouville property. This proves the gradient estimate.

Finally rescale at $R=\max\{\sqrt\varepsilon,|x'|\}$, using nested disks in the central region and nested annuli elsewhere, with rescaled displacement $R^{-1}u_\varepsilon(Ry,RZ)$. Its gradient on the larger patch is bounded by $CR^{q-2}$. Lemma~\ref{sph:lem:odd-local}, whose rescaled aspect ratio is $R$, gives $|e(u_\varepsilon)|\le CR^{q-1}$. The uniform comparisons between $R$ and $\sqrt\varepsilon+|x'|$ finish the proof.
\end{proof}

\begin{proof}[Proof of the upper bound in Theorem~\ref{sph:thm:general-cavity}]
The even part is controlled by Proposition~\ref{sph:thm:even-sharp}. For the odd part, Proposition~\ref{sph:prop:bending-spectrum} gives
\[
b_\nu=\min\left\{\sqrt{7-6\nu},\,\sqrt{\frac{19-\sqrt{73}}2}\right\}>\max\{2,1+q_*\}.
\]
Choose $\max\{1,q_*\}<q<\min\{2,b_\nu-1\}$. Lemma~\ref{sph:lem:bending-liouville} verifies the Liouville hypothesis of Proposition~\ref{sph:prop:odd-neck-criterion} for $Q(y)=|y|^2$, the constant tensor $\bC_0$, and this value of $q$.

The fixed-collar bound \eq{geom:eq:collar} applies to the odd projection, including both faces there. The gap has the form $\varepsilon+g(x')$, with $g(x')=|x'|^2+O(|x'|^4)$. Proposition~\ref{sph:prop:odd-neck-criterion} now gives on a fixed smaller neck
\[
|\nabla u_{\rm cav}^-|\le C\|\varphi\|_{C^{2,\gamma}}d(r)^{q-2},\qquad|e(u_{\rm cav}^-)|\le C\|\varphi\|_{C^{2,\gamma}}d(r)^{q-1}.
\]
Combine these bounds up to the faces with the even estimate, using $q>q_*=2-2\beta_\nu$ and $\varepsilon+r^2\le d(r)^2\le2(\varepsilon+r^2)$. The dilation of Section~\ref{sec:membrane} restores $\kappa$ in \eq{sph:eq:general-cavity-upper}, and \eq{geom:eq:collar} controls the gradient off the neck.
\end{proof}

\begin{remark}[Nonradial quadratic contact]
Proposition~\ref{sph:prop:odd-neck-criterion} gives the odd estimate for symmetric quadratic necks under its Liouville hypothesis. A nonradial quadratic form $Q$ couples angular frequencies. A linear change of coordinates that makes $Q$ radial also transforms the elastic operator.
\end{remark}

\section{Exterior loads attaining the membrane rates}\lb{sec:lower}
Throughout this section the cavities satisfy the geometric hypotheses of Theorem~\ref{sph:thm:general-cavity}, including the local $C^5$ assumption on $g$. For each leading membrane mode we construct a dual stress that detects it from the fixed outer sphere. A finite-energy correction preserves its contact singularity. Unique continuation then gives a nonzero outer traction, whose pairing with one fixed smooth outer displacement bounds the central strain energy from below. The volume of the central neck gives the $L^\infty$ lower bound in Theorem~\ref{sph:thm:general-cavity}.

\subsection{A dual stress for each leading mode}

\begin{prop}[Realization of the two leading modes]\lb{sph:thm:cav-lower}
For each $j\in\{0,2\}$ there is a fixed $\varphi_j\in C^\infty(\partial B_L;\bR^3)$ such that the corresponding solution satisfies, for all sufficiently small gaps,
\be\lb{sph:eq:3d-lower}
\|e(u_{\rm cav})\|_{L^2(G_\varepsilon(C\sqrt\varepsilon))}\ge c\,\varepsilon^{q_j/2},\qquad\|e(u_{\rm cav})\|_{L^\infty(G_\varepsilon(C\sqrt\varepsilon))}\ge c\,\varepsilon^{-1+q_j/2}.
\ee
Each $\varphi_j$ is reflection-even and has angular order $j$.
\end{prop}

\begin{proof}
Put $q=q_m$ and $\beta=1-q/2$. For $m=0$ use the axisymmetric radial mode, with no azimuthal displacement. For $m=2$ the polar components have angular frequency $2$. All three-dimensional fields constructed below are even in their first two components and odd in their third component under $z\mapsto-z$. The angular conditions are imposed in addition to this parity. We take real fields, using complex notation for angular modes when computing exponents.

\smallskip
\noindent\emph{A dual solution of the membrane equation.} For $m=0$, \eq{sph:eq:radial-model-ode} at $\varepsilon=0$ has the roots $-1\pm q_0$. For $m=2$, Lemma~\ref{sph:lem:membrane-spectrum} gives the four distinct roots $-1\pm q_2$ and $-1\pm\widehat q_2$, with $\widehat q_2>q_2$. Both selected numbers $q_0,q_2$ belong to $(0,2)$.

Set $s=-1-q$ and $s_\varepsilon=\sqrt\varepsilon$. On $r\ge s_\varepsilon$ we use a solution $v_\varepsilon$ of $\operatorname{div}'(H\bC_{\rm ps}e'(v_\varepsilon))=0$ with
\be\lb{sph:eq:dual-memb}
|D'^jv_\varepsilon|\le C_jr^{s-j},\qquad v_\varepsilon\longrightarrow r^s\Phi(\vartheta)\quad\text{locally in }C^j\text{ on }r>0,\quad0\le j\le4,
\ee
where $\Phi\ne0$ belongs to the chosen angular mode. To see this, rescale $r=s_\varepsilon R$. The radial system for the coefficient $1+R^2$, written in the variable $\log R$, is a constant matrix plus a perturbation $O(R^{-2})$. Its roots are distinct. Asymptotic integration \cite{Eastham85} gives a solution whose radial coefficients and first logarithmic derivatives, divided by $R^s$, converge to $\Phi$ and $s\Phi$, respectively. Differentiating the system gives the corresponding limits and bounds for the required higher derivatives. Extend down to $R=1$ by the ordinary differential equation and multiply the displacement by $s_\varepsilon^s$.

\smallskip
\noindent\emph{A balanced three-dimensional stress.} For $|z|<\delta(r)/2$ and $r\ge s_\varepsilon$, define
\[
\begin{gathered}
S=(H/\delta)\bC_{\rm ps}e'(v_\varepsilon),\qquad T=-z\operatorname{div}'S,\qquad\Sigma=\begin{pmatrix}S&T\\T^{\mathsf T}&n\end{pmatrix},\\
n=\left(\frac{z^2}{2}-\frac{\delta^2}{8}\right)\operatorname{div}'\operatorname{div}'S+\frac12\nabla'\delta\cdot T(x',\delta/2).
\end{gathered}
\]
Then $\operatorname{div}\Sigma=0$. Indeed, $\partial_zT=-\operatorname{div}'S$ and $\partial_zn=-\operatorname{div}'T$. Since $\operatorname{div}'(\delta S)=0$, $T(x',\pm\delta/2)=\pm S\nabla'\delta/2$. The normal stress $n$ has the same value on both faces. Multiplication by the unnormalized normals $(-\nabla'\delta/2,\pm1)$ gives zero traction.

Use the displacement
\be\lb{sph:eq:dual-displ}
U'=v_\varepsilon+\frac12z^2\left(\eta_{\rm ps}\nabla'\operatorname{div}'v_\varepsilon-\mu^{-1}\operatorname{div}'S\right),\qquad U_3=-\eta_{\rm ps} z\operatorname{div}'v_\varepsilon.
\ee
Its transverse shear stress is exactly $T$. The leading in-plane stress is $\bC_{\rm ps}e'(v_\varepsilon)$, and the term of degree zero in $z$ in its normal stress vanishes.

On $A_t$ with $t\ge s_\varepsilon$, $\delta=H+O(t^4)$, $\delta\le Ct^2$, and $|\nabla'\delta|\le Ct$. The in-plane discrepancy is bounded by $Ct^2|D'v_\varepsilon|+Ct^4t^{s-3}$. The normal stress and $n$ are both bounded by $Ct^4t^{s-3}$. Using \eq{sph:eq:dual-memb} and the $O(t^4)$ volume of the neck over $A_t$ gives
\be\lb{sph:eq:dual-error}
\|\bC_0e(U)-\Sigma\|_{L^2(G_\varepsilon\cap\{x'\in A_t\})}\le Ct^{2-q}.
\ee

\smallskip
\noindent\emph{The finite-energy correction.} Choose radial cutoffs $\chi_{\rm i}$ and $\chi_{\rm o}$. The first increases from zero to one on $[s_\varepsilon,2s_\varepsilon]$, with $|\nabla'\chi_{\rm i}|\le C/s_\varepsilon$. The second decreases from one to zero in a fixed annulus strictly inside $r<r_0$. Write $\chi=\chi_{\rm i}\chi_{\rm o}$ and extend $\chi U$ by zero off the neck. For $a_\varepsilon(v,w)=\int_{M_\varepsilon}\bC_0e(v):e(w)\,\rd x$, integration by parts using $\Sigma$ gives
\be\lb{sph:eq:dual-functional}
\begin{aligned}
 a_\varepsilon(\chi U,w)&=\int_{M_\varepsilon}R_\varepsilon:e(w)\,\rd x-
 \int_{M_\varepsilon}f_{\rm o}\cdot w\,\rd x+\mathcal L_\varepsilon(w),\\
 \mathcal L_\varepsilon(w)&=\int_{M_\varepsilon}I_\varepsilon^{\rm in}:e(w)\,\rd x-
 \int_{M_\varepsilon}f_{\rm i}\cdot w\,\rd x,
\end{aligned}
\ee
where
\[
\begin{aligned}
R_\varepsilon&=\chi(\bC_0e(U)-\Sigma)+\chi_{\rm i}\bC_0(U\odot\nabla\chi_{\rm o}),&f_{\rm o}&=\chi_{\rm i}\Sigma\nabla\chi_{\rm o},\\
I_\varepsilon^{\rm in}&=\chi_{\rm o}\bC_0(U\odot\nabla\chi_{\rm i}),&f_{\rm i}&=\chi_{\rm o}\Sigma\nabla\chi_{\rm i}.
\end{aligned}
\]
We extend these fields by zero off their supports. The dyadic sum in~\eq{sph:eq:dual-error} gives $\|R_\varepsilon\|_2\le C$ since $q<2$. The $L^2$ norm of $f_{\rm o}$ is bounded on a fixed annulus away from contact. There the anchored Korn estimate preceding~\eq{geom:eq:collar} bounds $\|w\|_2$ by $C\|e(w)\|_{L^2(M_\varepsilon)}$ for every $w$ of zero outer trace.

For every positive gap solve the zero outer trace variational problem for $V_\varepsilon$ with right-hand side $-\int R_\varepsilon:e(w)+\int f_{\rm o}\cdot w$. Coercivity gives $\|e(V_\varepsilon)\|_{L^2}\le C$. Let $Z_\varepsilon=\chi U+V_\varepsilon$. Then $Z_\varepsilon=0$ on $\partial B_L$ and
\be\lb{sph:eq:dualcorr-weak}
a_\varepsilon(Z_\varepsilon,w)=\mathcal L_\varepsilon(w)\quad\text{for all }w\text{ with zero outer trace}.
\ee
Projection under rotations about the axis and reflection in the midplane preserves the variational equation, so its solution has the chosen angular mode and parity. In the radial case a reflection in a plane through the axis also excludes the azimuthal component.

\smallskip
\noindent\emph{The limiting traction.} On compact portions of the matrix away from the origin, including the outer collar, the correction has a local $H^1$ bound by Korn's inequality on truncated domains as above. Subsequential weak limits give a contact solution $Z_*$, with zero outer trace and zero cavity traction away from the origin. More precisely, smooth pullbacks on each truncated portion allow tests supported away from the origin in \eq{sph:eq:dualcorr-weak}. The inner functional vanishes on those tests for all sufficiently small $\varepsilon$. The difference between $Z_*$ and the limit of $\chi U$ from \eq{sph:eq:dual-memb} has finite strain energy, by lower semicontinuity on an exhaustion away from the origin.

Uniqueness of this limit follows from the following estimate in the chosen symmetries. For $A_t^+=\{t/2<r<4t\}$, at $\varepsilon=0$ we have
\be\lb{sph:eq:cutoff-symm}
\|W/t\|_{L^2(G_0\cap\{x'\in A_t\})}\le C\|e(W)\|_{L^2(G_0\cap\{x'\in A_t^+\})}.
\ee
The estimate also holds at positive gaps when $s_\varepsilon\le t$. For the vertical component use its oddness in $z$ and one-dimensional Poincar\'e. For a radial horizontal component use $e_{\vartheta\vartheta}=W_r/r$. For the $m=2$ horizontal components, apply planar Korn--Poincar\'e on each horizontal section of the enlarged annulus. The angular mode is orthogonal to planar rigid displacements. At every height intersecting $A_t$, the enlarged horizontal section has outer radius $4t$ and inner radius in $[t/2,2t]$. Indeed, the gap is increasing in $r$. After scaling by $t$, Korn--Poincar\'e is uniform for this compact family of annuli. Integrating the section bounds proves \eq{sph:eq:cutoff-symm}.

The difference $W$ of two contact limits has finite strain energy and zero outer trace. Test its homogeneous equation with $\zeta_t^2W$, where $\zeta_t$ vanishes near the contact point, equals one outside $r=2t$ in the neck and has derivative at most $C/t$. This test can be defined to equal $W$ outside the neck since the two cavity faces bound the small neck component. By \eq{sph:eq:cutoff-symm}, the error from differentiating the cutoff tends to zero as $t\to0$. Hence $e(W)=0$, and the outer trace implies $W=0$. Thus the entire family has the same local weak limit.

On the other hand, the homogeneous in-plane field $r^{-1-q}\Phi$ has annular strain norm equal to a positive constant times $t^{-q}$ in the contact neck. Its in-plane strain is nonzero since this exponent does not describe a rigid displacement. The $z^2$ correction in \eq{sph:eq:dual-displ} contributes $O(t^{2-q})$ to the in-plane strain norm. Consequently the limit of $\chi U$ has infinite strain energy, so the finite-energy correction cannot cancel it. In particular, $Z_*\ne0$.

The outer traction of $Z_*$ is nonzero. If both Cauchy data vanished there, extension by zero across the outer sphere would solve the constant-coefficient Lam\'e equation locally across it. Analytic unique continuation and connectedness would give $Z_*=0$.

\smallskip
\noindent\emph{Reciprocity and the localized strain bound.} Choose a smooth outer displacement $\varphi$ in the same symmetry class with nonzero pairing with the limiting traction. For example one may take the traction itself, which is smooth on the outer sphere by boundary regularity. Choose a smooth lift $\widehat\varphi$ supported in a fixed outer collar. By local weak convergence, $a_\varepsilon(Z_\varepsilon,\widehat\varphi)$ converges to this nonzero pairing. Betti reciprocity, using zero outer trace of $Z_\varepsilon$ and \eq{sph:eq:dualcorr-weak} with $w=u_{\rm cav}-\widehat\varphi$, gives
\[
\mathcal L_\varepsilon(u_{\rm cav})=-a_\varepsilon(Z_\varepsilon,\widehat\varphi).
\]
For this fixed load the pairing is bounded away from zero for every sufficiently small gap.

The two terms in the inner functional of~\eq{sph:eq:dual-functional} are supported on $G_\varepsilon\cap\{x'\in A_{s_\varepsilon}\}$. From~\eq{sph:eq:dual-memb} and~\eq{sph:eq:dual-displ}, $\|I_\varepsilon^{\rm in}\|_{L^2}\le Cs_\varepsilon^{-q}$ and $\|f_{\rm i}\|_{L^2}\le Cs_\varepsilon^{-1-q}$. The primal solution has the selected symmetry by uniqueness. Applying \eq{sph:eq:cutoff-symm} to it gives
\[
0<c\le|\mathcal L_\varepsilon(u_{\rm cav})|\le Cs_\varepsilon^{-q}\|e(u_{\rm cav})\|_{L^2(G_\varepsilon(4s_\varepsilon))}.
\]
The volume of $G_\varepsilon(4s_\varepsilon)$ is bounded by $Cs_\varepsilon^4$. Since $\varepsilon=s_\varepsilon^2$ under our normalization, the two lower bounds have powers $\varepsilon^{q/2}=\varepsilon^{1-\beta}$ and $\varepsilon^{-1+q/2}=\varepsilon^{-\beta}$, respectively.
\end{proof}

\begin{proof}[Completion of Theorem~\ref{sph:thm:general-cavity}]
Put $q_*=\min\{q_0,q_2\}$. The upper bound follows from Section~\ref{sec:bending}. Choose $j\in\{0,2\}$ with $q_j=q_*$. Proposition~\ref{sph:thm:cav-lower} gives \eq{main:eq:sharp-lower}, because $1-q_*/2=\beta_\nu$. The fixed dilation restoring $\kappa$ changes only the constants and the radius of the central neck. Finally
\[
q_0^2-q_2^2=\sqrt{\nu^2-30\nu+33}-(5+\nu)
\]
has the sign of $1-5\nu$. This proves the claimed mode selection, including the coincident rates at $\nu=1/5$.
\end{proof}

\subsection{The radial strain and the scalar insulating exponent}

\begin{prop}[The radial strain profile and its scalar component]\lb{sph:prop:scalar-strain}
Take reflection-even meridional data as in Proposition~\ref{sph:thm:radial-sharp}. Put $H(r)=\varepsilon+r^2$ and $d(r)=\sqrt\varepsilon+r$. Let $Y_\varepsilon$ be the regular solution of
\be\lb{sph:eq:scalar-radial-Y}
(rHY_\varepsilon')'+(\nu H'-H/r)Y_\varepsilon=0,\qquad\lim_{r\to0}\frac{Y_\varepsilon(r)}r=\varepsilon^{(p-1)/2},\qquad p=q_0-1.
\ee
Define
\[
\Psi_\varepsilon(x')=Y_\varepsilon(r)\cos\vartheta=Y_\varepsilon(r)\frac{x_1}{r}.
\]
At $r=0$ the functions are defined by smooth continuation.

There is a linear functional $A_\varepsilon$ with $|A_\varepsilon[\varphi]|\le C\|\varphi\|_{C^{2,\gamma}}$ such that, in the orthonormal frame $(e_r,e_\vartheta,e_3)$ and up to both free faces of a fixed small neck,
\be\lb{sph:eq:radial-strain-profile}
e(u_{\rm cav})(x)=A_\varepsilon[\varphi]\operatorname{diag}\!\left(Y_\varepsilon',\frac{Y_\varepsilon}{r},-\eta_{\rm ps}\left(Y_\varepsilon'+\frac{Y_\varepsilon}{r}\right)\right)+O\!\left(\|\varphi\|_{C^{2,\gamma}}d(r)^p\right).
\ee
In particular,
\be\lb{sph:eq:strain-column-scalar}
e(u_{\rm cav})(x)e_1=A_\varepsilon[\varphi](\nabla'\Psi_\varepsilon(x'),0)+O\bigl(\|\varphi\|_{C^{2,\gamma}}d(r)^p\bigr).
\ee

The scalar profile satisfies the exact reduced equation
\be\nn
\operatorname{div}'(H\nabla'\Psi_\varepsilon)+2\nu\Psi_\varepsilon=0.
\ee
For $\lambda=0$, this is the homogeneous thickness-averaged operator of the three-dimensional insulated scalar problem \cite[Sections 2--3]{DLY25}.

There exists a fixed smooth radial loading for which the coefficient is bounded away from zero for all sufficiently small gaps. Under that loading there are constants $c,C>0$ for which
\be\lb{sph:eq:single-strain-lower}
c\,\varepsilon^{-1+\sqrt{(1-\nu)/2}}\le|e_{11}(u_{\rm cav})(0)|=|e_{22}(u_{\rm cav})(0)|\le C\,\varepsilon^{-1+\sqrt{(1-\nu)/2}}.
\ee
For $\nu=0$ the radial strain exponent is $(2-\sqrt2)/2$.
\end{prop}

\begin{proof}
We use the normalization of Section~\ref{sec:membrane}. Let $V$ be the thickness average of $u_r$ and put $f=V/Y_\varepsilon$. The Wronskian estimate in the proof of Proposition~\ref{sph:thm:radial-sharp} gives $|f'|\le C\|\varphi\|_{C^{2,\gamma}}$ and $|f|\le C \|\varphi\|_{C^{2,\gamma}}$ at a fixed radius. Hence $A_\varepsilon[\varphi]=f(0)=\varepsilon^{(1-p)/2}V'(0)$ is a bounded linear functional and $|f(r)-f(0)|\le C\|\varphi\|_{C^{2,\gamma}}r$. The bounds on $Y_\varepsilon'$ and $Y_\varepsilon/r$ in the same proof give
\[
 |V'-A_\varepsilon Y_\varepsilon'|+|V/r-A_\varepsilon Y_\varepsilon/r|\le C\|\varphi\|_{C^{2,\gamma}}d^p.
\]
On the patch of lateral size $R=\max\{r,\sqrt\varepsilon\}$, the radial gradient bound and the $O(R^4)$ patch volume give $\|e(u_{\rm cav})\|_2\le C\|\varphi\|_{C^{2,\gamma}}R^{p+1}$. The reconstruction~\eq{sph:eq:layer-approx} therefore gives a gradient error of $O(\|\varphi\|_{C^{2,\gamma}}R^p)$ up to both faces. Taking symmetric parts yields~\eq{sph:eq:radial-strain-profile}.

The first Cartesian component of $Y_\varepsilon e_r$ is $\Psi_\varepsilon$. Its planar gradient gives the first column of the symmetric tensor $\nabla'(Y_\varepsilon e_r)$, proving~\eq{sph:eq:strain-column-scalar}. In polar coordinates, \eq{sph:eq:scalar-radial-Y} gives
\[
 \operatorname{div}'(H\nabla'\Psi_\varepsilon)
 =\left(\frac{(rHY_\varepsilon')'}r-\frac{HY_\varepsilon}{r^2}\right)\cos\vartheta
 =-\nu\frac{H'}r\Psi_\varepsilon=-2\nu\Psi_\varepsilon.
\]
For the fixed radial load of Proposition~\ref{sph:thm:cav-lower}, the strain supremum on $r<C_0\sqrt\varepsilon$ is at least $c\,\varepsilon^{(p-1)/2}$. The profile there is bounded by $C|A_\varepsilon|\varepsilon^{(p-1)/2}+C\|\varphi\|_{C^{2,\gamma}}\varepsilon^{p/2}$, so $|A_\varepsilon|\ge c_1>0$ for small gaps. Evaluating the profile at the axis, where $Y_\varepsilon'(0)=\lim_{r\to0}Y_\varepsilon/r=\varepsilon^{(p-1)/2}$, proves~\eq{sph:eq:single-strain-lower}.
\end{proof}

\section{Finite contrast for three-dimensional inclusions}\label{sec:spherical-comparisons}
We apply the energy comparisons of Section~\ref{sec:common-variational}. For soft inclusions, the cavity gradient bound controls the energy of a trace extension. For stiff inclusions, a traction estimate for the rigid solution applies at any nondegenerate quadratic contact. In the torsional sector the system reduces to a scalar equation, which also covers balls of unequal radii.

\subsection{Comparison with reflected cavities}
\begin{proof}[Proof of Theorem~\ref{thm:spherical-contrast}\textup{(i)}]
Write $v=u_{\rm fin}-u_{\rm cav}$ in the matrix and $\beta=\beta_\nu\in(0,1)$. Write $B_i$ for either inclusion and $\Gamma_i=\partial B_i$. On $\Gamma_i$, the face and bounds away from the neck of Theorem~\ref{sph:thm:general-cavity} give $|\nabla_\tau u_{\rm cav}|\le C\|\varphi\|_{C^{2,\gamma}}\varepsilon^{-\beta}$. Subtract the surface mean of the trace. Surface Poincar\'e bounds its $H^1(\Gamma_i)$ norm by the same quantity. A bounded right inverse of the trace $H^1(B_i)\to H^{1/2}(\Gamma_i)$ extends it into the inclusion. Both constants are uniform because each inclusion is a translate of a fixed domain. Restoring the mean preserves the gradient bound in $L^2(B_i)$. By \eq{sph:eq:energy-density}, $\bC_1\le d_{\rm s}\bC_0$. The variational comparison \eq{eq:common-soft-lift} gives
\be\lb{eq:soft-balls-strain-energy}
E:=\|e(v)\|_{L^2(M_\varepsilon)}\le C\|\varphi\|_{C^{2,\gamma}}\sqrt{d_{\rm s}}\,\varepsilon^{-\beta}.
\ee

On $r=r_0$, $z=0$, a fixed portion of the matrix connects the endpoints of the radial lines to a fixed open patch of the outer boundary. The anchored Korn inequality and interior regularity give $|\nabla v|\le CE$ there. Write $H(r)=\varepsilon+r^2$. The gap and slope bounds give matrix balls of radius $c_\theta H(r)$ throughout $N_\varepsilon^\theta$. Lemma~\ref{lem:neck-line} bounds the midplane gradient by $CEH^{-7/4}$. Scaled Korn and interior estimates along a vertical segment in $N_\varepsilon^\theta$ give $|\nabla^2v|\le C_\theta H^{-5/2}E$. Its length is at most $CH$, so
\[
 |\nabla v(x)|\le C_\theta E(H^{-7/4}+H^{-3/2})\le C_\theta EH^{-2}.
\]
Together with \eq{eq:soft-balls-strain-energy} and $H\ge\varepsilon$, this proves \eq{eq:soft-balls-weighted}. Adding the cavity bound yields \eq{eq:soft-field-explicit}. For $d_{\rm s}\le M\varepsilon^4$, the error is at most $C_{\theta,M}\|\varphi\|_{C^{2,\gamma}}H^{-\beta}$.

For the lower bound choose the fixed load of Proposition~\ref{sph:thm:cav-lower} for the leading mode. By \eq{sph:eq:3d-lower}, for some fixed $A$,
\[
 \|e(u_{\rm cav})\|_{L^2(S_\varepsilon)}\ge c\,\varepsilon^{1-\beta},\qquad
 S_\varepsilon=G_\varepsilon(A\sqrt\varepsilon).
\]
Its volume is at most $C_A\varepsilon^2$, and the two face strips removed from it in $N_\varepsilon^\theta$ have total volume at most $C_A\theta\varepsilon^2$. The cavity upper bound on those strips is $c\,\varepsilon^{-\beta}$. Choose a fixed $\theta_*\in(0,1/2)$ so small that the $L^2$ norm on the strips is at most $(c/2)\varepsilon^{1-\beta}$.

If $d_{\rm s}\le M\varepsilon^4$, \eq{eq:soft-balls-strain-energy} gives $E\le C_M\|\varphi\|_{C^{2,\gamma}}\varepsilon^{2-\beta}$. Subtracting this error leaves an $L^2$ lower bound of $(c/4)\varepsilon^{1-\beta}$ for $e(u_{\rm fin})$ on $S_\varepsilon\cap N_\varepsilon^{\theta_*}$ for small gaps. Division by the square root of its volume proves the assertion.
\end{proof}

\subsection{Comparison with rigid inclusions}
Let two fixed strictly convex $C^{2,\gamma}$ inclusions in a fixed $C^{2,\gamma}$ domain $\Omega_3$ touch only at the origin with common tangent plane $z=0$. Translate the upper and lower inclusions by $\pm\varepsilon e_3/2$, keeping them a fixed positive distance from the outer boundary. Write their facing graphs as $z_\pm(x')=\pm\varepsilon/2+p_\pm(x')$, where $p_\pm(0)=0$ and $D p_\pm(0)=0$. Assume $D^2(p_+-p_-)(0)$ is positive definite. For $r=|x'|<r_0$ with $r_0$ fixed and small, their gap satisfies
\[
 c(\varepsilon+r^2)\le\delta(x'):=z_+(x')-z_-(x')\le C(\varepsilon+r^2).
\]
For $0<\theta<1/2$ the interior region is
\[
 N_\varepsilon^\theta=\{r<r_0,\ z_-+\theta\delta\le z\le z_+-\theta\delta\}.
\]
As in Section~\ref{sec:main-contrast}, write $d_{\rm h}=\max\{\mu/\mu_1,K/K_1\}$ for the larger inverse modulus ratio.

\begin{prop}[Comparison in a quadratic gap]\lb{sph:thm:hard}
Let $\varphi\in C^{2,\gamma}(\partial\Omega_3;\bR^3)$. Assume both tensors are positive on symmetric matrices and $d_{\rm h}\le1$. Estimates~\eq{eq:main-balls-rigid} and~\eq{eq:hard-field-explicit} hold in $N_\varepsilon^\theta$, and the stiff-side energy bound of Remark~\ref{rem:comparison-energy} holds with $B_L$ replaced by $\Omega_3$. 

For the expansion, assume also that reflection in each coordinate plane preserves $\Omega_3$ and the union of the inclusions. Take odd data $\varphi(-x)=-\varphi(x)$ and facing graphs with $p_\pm(x')=\pm\frac12\kappa r^2+O(r^{2+\gamma})$ for some $\kappa>0$. Let $b[\varphi]$ be the resultant matrix traction on the upper inclusion for the touching minimizer with common rigid displacement \cite[(2.8)--(2.9)]{LX25}, with the normal pointing out of that inclusion. Set $a=(b_1/\mu,b_2/\mu,b_3/(\lambda+2\mu))$, $x'_\parallel=(x_1,x_2,0)$ and $\bar u=(z-z_-)/\delta$. Define the matrix linear in $x'$ by
\[
 T_b(x')=\kappa\,a\otimes x'_\parallel-
 \frac{\kappa(\lambda+\mu)}{\mu(\lambda+2\mu)}
 \big[(b\cdot x'_\parallel)e_3+b_3x'_\parallel\big]\otimes e_3 .
\]
Then the rigid gradient has the expansion
\be\nn
 \nabla u_\infty(x)=\frac{\kappa}{\pi\delta|\log\varepsilon|}
  \big[a\otimes e_3+(1-2\bar u)T_b(x')\big]+\mathcal R_\infty(x),
\ee
where throughout the narrow region
\[
 |\mathcal R_\infty(x)|\le C\|\varphi\|_{C^{2,\gamma}}
 \left(1+\frac{r^{1+\gamma}}{\delta|\log\varepsilon|}
             +\frac{1}{\delta|\log\varepsilon|^2}\right).
\]
For $x\in N_\varepsilon^\theta$, the same displayed profile holds for $\nabla u_{\rm fin}$ with the additional remainder $\nabla(u_{\rm fin}-u_\infty)$ bounded in~\eq{eq:main-balls-rigid}.
Consequently, if the resultant is nonzero and $d_{\rm h}=o(\varepsilon^{3/2})$, then $\varepsilon |\log\varepsilon|\nabla u_{\rm fin}(0)\to (\kappa/\pi)a\otimes e_3$.
\end{prop}

\begin{proof}
Set $v=u_{\rm fin}-u_\infty$. By \eq{sph:eq:energy-density}, $\bC_1\ge d_{\rm h}^{-1}\bC_0$. In the neck the quadratic faces place $z$ within $Cr^2$ of $[-\varepsilon/2,\varepsilon/2]$. Hence the distance to the shortest segment lies between $r$ and $r\sqrt{1+Cr^2}$. The positive definite relative Hessian gives the local convexity condition of \cite[Remark~1.2]{BLL17}. Substituting the distance bounds in \cite[Remark~1.1, (1.12)]{BLL17} gives~\eq{eq:rigid-BLL-exact} up to the matrix-side face traces. The same cited estimate bounds the rigid gradient away from the facing caps by $C\|\varphi\|_{C^2}$. On each translated interface $H^{1/2}\hookrightarrow L^4$ with a fixed constant. On its facing graph the bounded slope gives surface measure at most $Cr\,dr\,d\vartheta$ in polar coordinates. By duality,
\[
\|t_{\infty,i}\|_{H^{-1/2}}\le C\|t_{\infty,i}\|_{L^{4/3}}\le C\|\varphi\|_{C^{2,\gamma}(\partial\Omega_3)}\varepsilon^{-1/4}/|\log\varepsilon|.
\]
Indeed, writing $H=\varepsilon+r^2$, we have $\int_0^{r_0}rH^{-4/3}\,dr\le c\,\varepsilon^{-1/3}$ and $\int_0^{r_0}r(r/H)^{4/3}\,dr\le C$. The energy bound follows from \eq{eq:common-rigid-traction}. On a matrix ball of radius $c_\theta H$, the interior estimate gives $|\nabla v(x)|\le C_\theta H^{-3/2}\|\nabla v\|_2$, which gives~\eq{eq:main-balls-rigid}. Adding the rigid bound~\eq{eq:rigid-BLL-exact} yields~\eq{eq:hard-field-explicit}.

For the expansion, take $E_{ij}=e_i\otimes e_j$ with the third coordinate normal to the gap. In \cite[Theorem~2.4(i), (2.12), p.~761]{LX25}, $\mathbb B_{3,\mathrm I}=a\otimes e_3$ and $x_1\mathbb B_{3,\mathrm {II}}+x_2\mathbb B_{3,\mathrm {III}}=T_b(x')$. The expansions of $p_\pm$ give the gap coefficient $\kappa$ in their (2.7). Expanding the factor $1+O(|\log\varepsilon|^{-1})$ gives the stated additive remainder, since $0\le\bar u\le1$ and $b$ is bounded linearly by the data norm. At the midpoint $T_b(0)=0$ and $\delta=\varepsilon$. More generally, multiplying the error in~\eq{eq:main-balls-rigid} by $(\varepsilon+r^2)|\log\varepsilon|$ bounds it by $C_\theta\|\varphi\|_{C^{2,\gamma}}\sqrt{d_{\rm h}/\varepsilon^{3/2}}$ throughout $N_\varepsilon^\theta$.
The bodies of revolution in $B_L$ are invariant under reflection in each vertical coordinate plane, and midplane reflection exchanges them. Their facing graphs have $p_\pm(x')=\pm\frac12\kappa r^2+O(r^{2+\gamma})$ by the quadratic contact hypothesis. This proves Theorem~\ref{thm:spherical-contrast}\textup{(ii)}, including the midpoint assertion for odd data.
\end{proof}

\subsection{Torsion of spherical cavities and soft inclusions}
Take two balls of fixed radii $R_\pm$ with centers $\pm(R_\pm+\varepsilon/2)e_3$ inside a ball centered at the origin, at a fixed positive distance from its boundary. Use the interior neck $N_\varepsilon^\theta$ defined above. Put $Jx=(-x_2,x_1,0)$ and prescribe $f(r,z)Jx$ on the outer sphere. Assume that the lift $f(|y|,z)$, $y\in\bR^4$, to the five-dimensional outer sphere belongs to $C^{2,\gamma}$. All norms of $f$ below are taken on this lifted outer sphere. Write $u_{\rm cav}$ for the cavity solution and $u_{\rm fin}$ for the finite solution with $k_{\rm sh}\in(0,1]$ and $K_1>0$.

\begin{prop}[Torsional rotation and the soft limit]\lb{sph:thm:torsion}
Let $H=\varepsilon+r^2$ and $\alpha=\sqrt7-2$. The cavity gradient is bounded by $C\|f\|_{C^{2,\gamma}}$ throughout the matrix. There is a scalar $\omega_\varepsilon$ with $\nabla u_{\rm cav}(0)=\omega_\varepsilon J$ such that
\[
 |\nabla u_{\rm cav}-\omega_\varepsilon J|\le C\|f\|_{C^{2,\gamma}}H^{\alpha/2}
\]
throughout the neck, including the free faces. For fixed data, $\omega_\varepsilon$ converges to the contact value of the bounded finite-energy insulated solution in the lifted touching domain. In $N_\varepsilon^\theta$,
\be\lb{sph:eq:soft-comparison}
 |\nabla(u_{\rm fin}-u_{\rm cav})|\le C_\theta \|f\|_{C^{2,\gamma}}\frac{\sqrt{k_{\rm sh}}}{H\sqrt{H+r}}.
\ee
The constants are independent of $K_1$. Consequently $\nabla u_{\rm fin}$ is uniformly bounded there if $k_{\rm sh}=O(\varepsilon^3)$, and $\nabla u_{\rm fin}=\omega_\varepsilon J+o(1)$ uniformly for $r\le L\sqrt\varepsilon$ there if $k_{\rm sh}=o(\varepsilon^3)$, for each fixed $L$.
\end{prop}

\begin{proof}
Lift the balls with their axial centers and radii to $\bR^4\times\bR$. Let $\omega_{\rm cav}$ and $\omega_{\rm fin}$ solve the scalar insulated and transmission problems there with outer trace $f$ and, in the latter case, conductivities $1,k_{\rm sh}$. Uniqueness makes both functions invariant under rotations of $y$. For $r=|y|$ and $u=\omega(r,z)Jx$, direct differentiation gives
\[
 \operatorname{div}u=0,\quad \Delta u=Jx(\omega_{rr}+3r^{-1}\omega_r+\omega_{zz}),\quad
 (\bC_j e(u))n=\mu_j(\partial_n\omega)Jx.
\]
The first two identities hold in each phase. The traction identity holds on each spherical interface since $Jx\cdot n=0$. Hence $u_{\rm cav}=\omega_{\rm cav}Jx$, $u_{\rm fin}=\omega_{\rm fin}Jx$, and $\nabla u=\omega J+(Jx)\otimes\nabla_x\omega$.

The maximum principle bounds $|\omega_{\rm cav}|$ by $\|f\|_{C^{2,\gamma}}$. The five-dimensional insulated estimate \cite[Theorem~1.1]{DLY25} gives $|\nabla\omega_{\rm cav}|\le C\|f\|_{C^{2,\gamma}}H^{(\alpha-1)/2}$ up to the faces, and ordinary boundary regularity applies away from the neck. Put $\omega_\varepsilon=\omega_{\rm cav}(0,0)$. Integrating from the origin along the mid-surface to radius $r$ costs at most $C\|f\|_{C^{2,\gamma}}\int_0^r t^{\alpha-1}\,\rd t$, and the vertical segment costs at most $C\|f\|_{C^{2,\gamma}}H^{(1+\alpha)/2}$. Together with $r|\nabla\omega_{\rm cav}|\le C\|f\|_{C^{2,\gamma}}H^{\alpha/2}$, this proves both cavity gradient bounds.

The tangential gradient bound is square integrable on each lifted cap, since
\[
 \int_0^{r_0}r^3(\varepsilon+r^2)^{\alpha-1}\,\rd r\le C.
\]
Thus the lifted trace has $H^1$ norm at most $C\|f\|_{C^{2,\gamma}}$ on each sphere. Extend the traces harmonically into the balls, with total Dirichlet energy at most $C\|f\|_{C^{2,\gamma}}^2$. Use this extension as a competitor for $\omega_{\rm fin}$. The mixed term in the matrix energy vanishes by the insulated equation, so $W=\omega_{\rm fin}-\omega_{\rm cav}$ satisfies $\|\nabla W\|_{L^2(\widehat M_\varepsilon)}\le C\|f\|_{C^{2,\gamma}}\sqrt{k_{\rm sh}}$ in the lifted matrix. At an interior neck point one can pack at least $c_\theta(1+r/H)^3$ disjoint balls of radius $c_\theta H$ with centers on its rotation orbit. The mean-value bound for $|\nabla W|^2$ has the same left-hand side at every center. Summing it gives
\[
 |\nabla W|\le C_\theta \|f\|_{C^{2,\gamma}}\sqrt{k_{\rm sh}}\,H^{-1}(H+r)^{-3/2}.
\]
At fixed $r_0$ the zero outer trace, Poincar\'e's inequality on a fixed tube, and interior estimates give $|W|\le C\|f\|_{C^{2,\gamma}}\sqrt{k_{\rm sh}}$ on the mid-surface. Integrating inward along that surface uses
\[
 \int_r^{r_0}\frac{\rd t}{H(t)(H(t)+t)^{3/2}}\le\frac{2}{H(r)\sqrt{H(r)+r}}.
\]
The vertical segment has length $O(H)$. Hence
\[
 |W|+r|\nabla W|
 \le C_\theta \|f\|_{C^{2,\gamma}}\frac{\sqrt{k_{\rm sh}}}{H\sqrt{H+r}},
\]
and \eq{sph:eq:soft-comparison} follows from the formula for $\nabla u$.

Finally a fixed outer lift bounds the insulated Dirichlet energy. Local estimates give subsequential convergence away from contact to bounded finite-energy solutions with trace $f$ and zero face flux. For the difference $v$ of two limits, test with $\zeta_\rho^2v$, where $\zeta_\rho$ cuts off $B_\rho$ inside $B_{2\rho}$. The cutoff error is $O(\rho^3\|v\|_\infty^2)$ in dimension five, so the outer trace gives uniqueness. The estimate $|\omega_{\rm cav}-\omega_\varepsilon|\le C\|f\|_{C^{2,\gamma}}H^{\alpha/2}$ then identifies the limit of $\omega_\varepsilon$ by first passing to contact at fixed $r>0$ and then letting $r\downarrow0$.
\end{proof}

\section{The disk traction map and the two-circle interaction}
We compute the Fourier multipliers of the disk traction map, which sends a prescribed boundary displacement to the traction of its Lam\'e extension. These multipliers determine the extension energy in the soft-disk comparison of Section~\ref{sec:soft-disk-proof}. For stiff disks, they give the material coefficients of a boundary integral system on the two circles. In bipolar coordinates we separate its $4\times4$ principal symbol from the remainder terms, a decomposition used in Section~\ref{sec:pre} to estimate the layer densities.

Chao and Gao \cite{CG98} gave an explicit solution for two circular inclusions in an infinite matrix. Here the displacement is prescribed on a fixed outer boundary. In our neck coordinates the rigid estimate of Bao, Li, and Li \cite[Theorem~1.1, (1.9)]{BLL15} gives
\[
 |\nabla u_\infty(x_1,t)|\le
 \frac{C\|\varphi\|_{C^{1,\gamma}(\partial\Omega)}}{\sqrt\varepsilon+|t|}
 \qquad ((x_1,t)\in G_\varepsilon(r_0)).
\]
Their estimate also bounds the matrix gradient uniformly outside $G_\varepsilon(r_0)$. Formula~\eq{finite:eq:mat} recovers the rigid spatial bound as both inverse contrasts tend to zero, and \eq{finite:eq:diff} controls the error for the same outer data. For rigid inclusions, Kang and Yu \cite{KY19} described the planar singular fields, and Dong, Li, Teng and Zhang \cite{DLTZ} obtained higher derivative estimates.

For the cavity comparison we also need the rotation of the displacement on each circle. Proposition~\ref{thm:disk-cavity} identifies the stress resultants that determine the singular part of the full gradient.

\subsection{The disk traction map}
Write $s_\varepsilon=\sqrt\varepsilon$. The background Kelvin tensor is
\be\label{kelvin}
 \Gamma_0(x)=\frac{1}{2\pi\mu(\kappa_0+1)}
 \left(-\kappa_0\log|x|\,I_2+\frac{x\otimes x}{|x|^2}\right),
 \qquad -\Div(\bC_0e(\Gamma_0))=\delta_0 I_2,
\ee
where the two material parameters used in the complex potentials are
\be\label{params}
 \kappa_j=\frac{\lambda_j+3\mu_j}{\lambda_j+\mu_j}
 =1+\frac{2\mu_j}{K_j}>1,\qquad j=0,1 .
\ee

Write $\cT_0u=\bC_0e(u)n$ for the background traction, with $n$ pointing out of the disk. Identify points and vectors with complex numbers. On a disk, a homogeneous background Lam\'e displacement has the representation
\be\lb{km}
2\mu U=\kappa_0 f-z\overline{f'}-\overline g
\ee
for holomorphic $f,g$ and $\kappa_0$ in \eq{params}; see~\cite{M53}. Relative to the outward unit normal $n=e^{\ri\theta}$, the complex traction is
\be\lb{traction}
(\cT_0U)_1+\ri(\cT_0U)_2=-\ri\pa_\theta\big((\kappa_0+1)f-2\mu U\big).
\ee

On a translated unit circle, write a displacement trace as $w(\theta)=\sum_{n\in\bZ}w_n\mathrm{e}^{\ri n\theta}$ in complex notation. Let $\Lambda:H^1(S^1;\bC)\to L^2(S^1;\bC)$ be the real-linear map sending this trace to the background traction of its homogeneous disk extension. Comparing the Fourier coefficients in \eq{km} and \eq{traction} gives
\be\lb{bsw:diskmap}
(\Lambda w)_n=\begin{cases}
2\mu nw_n/\kappa_0,&n\ge2,\\
2\mu(w_1+\overline{w_1})/(\kappa_0-1),&n=1,\\
0,&n=0,\\
-2\mu nw_n,&n\le-1.
\end{cases}
\ee
Indeed \eq{km} gives $\kappa_0 f_n=2\mu w_n$ for $n\ge2$ and $\kappa_0 f_1-\overline{f_1}=2\mu w_1$. Substitution into \eq{traction} proves \eq{bsw:diskmap} for trigonometric polynomials and hence, by density, for $H^1$ traces.

These multipliers give bounded maps $H^q\to H^{q-1}$ for real $q$. Applied to the inclusions, \eq{bsw:diskmap} also shows that vanishing shear and bulk moduli give a vanishing traction map. The $n=-1$ and real $n=1$ modes give the converse. On a circle of radius $r_i$ the multipliers are divided by $r_i$, and the arclength element is $r_i\,\rd\theta$.

For either disk translated and dilated to $B_1$, define $F_w(z)=\sum_{n\ge1}w_nz^n$, where $w_n$ are the complex displacement Fourier coefficients, and let $\mathcal Pw$ be the componentwise Poisson extension. The disk extension is
\be\lb{sildiskformula}
Z(z)=\mathcal Pw(z)+\frac{1-|z|^2}{\kappa_0}\overline{\left(\frac{F_w'(z)-F_w'(0)}{z}\right)}.
\ee
The quotient at zero is removable. In \eq{km} the coefficients of $f$ satisfy $\kappa_0 f_n=2\mu w_n$ for $n\ge2$. On subtracting the Poisson extension of the boundary values, the term $-z\overline{f'}$ contributes $(1-|z|^2)\sum_{n\ge2}n\overline{f_n}\,\overline z^{n-2}$. The constant and first modes give zero correction. This proves \eq{sildiskformula} for polynomials and hence by continuity for the traces used here.

\subsection{Material factors and the boundary system}
By \eq{bsw:diskmap}, the inclusion-to-matrix traction ratio equals $k_{\rm sh}$ on the negative modes and $k_{\rm sh}\kappa_0/\kappa_1$ on modes of order at least two. For the circle multipliers write $k=k_{\rm sh}$. If $k>\max\{1,\kappa_1/\kappa_0\}$, set

\be\nn
A_+=\frac{k+\kappa_1}{k-\kappa_1/\kappa_0},\qquad A_-=\frac{k+\kappa_0^{-1}}{k-1},\qquad \eta_A=\frac12\log(A_+A_-).
\ee
The $\ell$th term of the Neumann series \eq{highseries} carries the material factor $e^{-\ell\eta_A}$. The identity
\[
\frac{\kappa_1}{\kappa_0k}=\frac{1}{\kappa_0}\left(\frac{\mu}{\mu_1}+\frac{2\mu}{\lambda_1+\mu_1}\right)
\]
exhibits both moduli on the positive modes. For $k_{\rm sh}$ and $k_{\rm vol}$ sufficiently large it yields
\be\lb{finite:eq:compare}
c\bigl(k_{\rm sh}^{-1}+k_{\rm vol}^{-1}\bigr)\le \eta_A\le C\bigl(k_{\rm sh}^{-1}+k_{\rm vol}^{-1}\bigr).
\ee
As both inverse ratios tend to zero,
\[
\eta_A=\frac{(\kappa_0+1)\mu}{\kappa_0}\left(\frac{1}{\mu_1}+\frac{1}{\lambda_1+\mu_1}\right)+O\!\left(\left(\frac{\mu}{\mu_1}+\frac{\lambda+\mu}{\lambda_1+\mu_1}\right)^2\right).
\]
The remainder is uniform in the ratio of the two inverse moduli.

Section~\ref{sec:common-variational} gives $\|u_{\rm fin}\|_{H^1(\Omega)}\le C\|\varphi\|_{C^{1,\gamma}}$ whenever $\bC_1\ge\bC_0$. Localize by a cutoff equal to one on $\overline{\Omega_1}$ and supported in $\Omega$. Its transition region stays a fixed distance from the disks. The resulting source $f_{\rm fin}$ has Kelvin potential $H_{\rm fin}=\Gamma_0*f_{\rm fin}$ with uniformly bounded derivatives of each fixed order near the disks. The same construction applies at rigidity.

Write $\Gamma_i=\partial D_i$, $w_i=u_{\rm fin}|_{\Gamma_i}$ and let $\Lambda_j$ denote the disk traction map for $\bC_j$ with the normal pointing out of the disk. Keep the localized field in the matrix and extend $w_i$ into $D_i$ by the homogeneous background equation. The resulting compactly supported field $Z$ satisfies
\[
 -\Div(\bC_0e(Z))=f_{\rm fin}-\sum_i p_i\delta_{\Gamma_i},\qquad p_i=(\Lambda_1-\Lambda_0)w_i,
\]
by the traction jump.
Let $q_i=\Lambda_1 w_i$ be the matrix traction. Both $q_i$ and $\Lambda_0w_i$ annihilate rigid motions, so $p_i$ has zero force and torque. Testing the localized equation against constant displacements gives $\int_{\bR^2}f_{\rm fin}\,\rd x=0$. Define
\[
 (S_i p)(x)=\int_{\Gamma_i}\Gamma_0(x-y)p(y)\,\rd S_y .
\]
The zero forces cancel the nondecaying terms of the Kelvin potentials. Consequently $Z-H_{\rm fin}+\sum_i S_i p_i$ solves the homogeneous background equation on $\bR^2$ and tends to zero at infinity. The interior gradient estimate makes it constant, and decay makes it zero. Taking traces gives

\be\lb{finite:eq:boundary}
p_i=(\Lambda_1-\Lambda_0)w_i,\qquad C_ip_i+\sum_{j=1}^2S_jp_j+R_i=H_{\rm fin}\quad\text{on }\Gamma_i.
\ee

Here $C_i$ inverts $\Lambda_1-\Lambda_0$ on the complement of the rigid modes. The rigid displacement $R_i$ has the same boundary mean as $w_i$ and the same mean curl as $u_{\rm fin}$ in $D_i$.

By \eq{bsw:diskmap}, $C_i$ has nonexceptional multipliers
\[
\frac{\kappa_1}{2\mu(k-\kappa_1/\kappa_0)n}\quad(n\ge2),\qquad\frac{1}{2\mu(k-1)|n|}\quad(n\le-1),
\]
and real first-mode multiplier $1/[2(\lambda_1+\mu_1-\lambda-\mu)]$. All three prefactors are $O(\eta_A)$. At rigidity $C_i=0$ and the density equals the matrix traction $q_i$.

We also have, uniformly,
\be\lb{finite:eq:weak-density}
\sum_i\norm{p_i}_{H^{-1}(\Gamma_i)}\le C\|\varphi\|_{C^{1,\gamma}}.
\ee
Subtracting its value at the facing point from an $H^1$ test trace preserves its pairing with the matrix traction $q_i$ by zero force. Lift the trace with zero value on the other disk and on the outer boundary. For the right disk, write $g(t)$ for the adjusted trace on the facing arc and use $(1/2+x_1/\delta_\varepsilon(t))g(t)$ on each section, with a cutoff near the ends. Its squared gradient integral in the gap is bounded by
\[
C\int_{-r_0}^{r_0}\left(\delta_\varepsilon|g'|^2+\frac{1+|\delta_\varepsilon'|^2}{\delta_\varepsilon}|g|^2\right)\,\rd t\le C\|g\|_{H^1(-r_0,r_0)}^2.
\]
The last inequality uses $\delta_\varepsilon\ge ct^2$ and the one-dimensional Hardy inequality for $g(0)=0$. Reflection gives the left-disk lift, and fixed-size extensions complete both lifts off the gap. Integration by parts in the matrix controls $q_i$ in $H^{-1}$. The uniform $H^1$ bound controls $\Lambda_0w_i$ in $H^{-1/2}$. Finally, the boundary mean and the mean curl are bounded trace functionals on each fixed disk, so the coefficients of $R_i$ are bounded.

\subsection{The Kelvin layer in bipolar coordinates}
Write $p(w)=\sum_{n\in\bZ}p_nw^n$ on a unit circle in the angular coordinate about its center. Force and torque balance mean $p_0=0$ and $\Im p_1=0$. We use the periodic Hilbert transform $\mathrm H(e^{\ri n\theta})=-\ri\operatorname{sgn}(n)e^{\ri n\theta}$, with $\mathrm H1=0$, componentwise for vector fields. The single-layer trace on the circle is
\be\lb{circle2}
\cS p=\frac{\kappa_0}{2\mu(\kappa_0+1)}\sum_{n\ne0}\frac{p_nw^n}{|n|}-\frac{\bar p_1w}{2\mu(\kappa_0+1)}.
\ee
In complex notation, the Kelvin tensor~\eq{kelvin} acts by
\[
\Gamma_0(z-v)p=\frac{-2\kappa_0\log|z-v|p+p+(z-v)(\bar z-\bar v)^{-1}\bar p}{4\pi\mu(\kappa_0+1)}.
\]
On one circle $(z-v)/(\bar z-\bar v)=-zv$, and the Fourier expansion of $-\log|z-v|$ gives~\eq{circle2}.

Put $a=1+\varepsilon/2$, $s_b=\sqrt{a^2-1}$, and $r_{\rm b}=a-s_b$. The second bipolar radius is $r_{\rm b}^{-1}=a+s_b$. Set $\rho=r_{\rm b}^2$ and $\eta=-\log\rho$, so $c\sqrt\varepsilon\le\eta\le C\sqrt\varepsilon$. For $z=e^{\ri\theta}$ the two circle parametrizations are
\be\label{bipolar}
 z_1(z)=s_b\frac{1+r_{\rm b}z}{1-r_{\rm b}z},\qquad
 z_2(z)=s_b\frac{1+r_{\rm b}^{-1}z}{1-r_{\rm b}^{-1}z},\qquad
 B(z)=\frac{z-r_{\rm b}}{1-r_{\rm b}z}.
\ee
The rational function $J(z)=zB'(z)/B(z)$ equals $(1-r_{\rm b}^2)/|1-r_{\rm b}e^{\ri\theta}|^2$ on the unit circle. This is the arclength factor on both disks. Their centered coordinates are $w_1=B(z)$ and $w_2=-1/B(z)$, with reversed orientation on the second circle. Each weighted density $\tau_i(\theta)=J(z)p_i(z_i(z))$ has zero mean. For $|z|=|v|=1$, direct subtraction gives

\[
z_1(z)-z_2(v)=\frac{2s_b(1-\rho z/v)}{(1-r_{\rm b}z)(1-r_{\rm b}/v)},\qquad B(z)-B(v)=\frac{(1-r_{\rm b}^2)(z-v)}{(1-r_{\rm b}z)(1-r_{\rm b}v)}.
\]
Taking logarithms of the moduli gives a convolution on the nonzero modes. The target-dependent term vanishes by zero total force. The source-dependent term contributes only to the zero Fourier mode.

The periodic primitives $Q_i$ of the physical densities $p_i$ satisfy
\[
\partial_\theta(Q_1\circ w_1)=\tau_1,\qquad \partial_\theta(Q_2\circ w_2)=-\tau_2.
\] Composition with $B$ preserves the holomorphic and antiholomorphic parts up to constants, whereas composition with $-1/B$ exchanges them. In the physical angular variable differentiation multiplies mode $n$ by $\ri n$. Thus the multiplier $1/|n|$ is $\ri$ times the periodic primitive on positive modes and $-\ri$ times the primitive on negative modes. Together with the reversed orientation on the second disk, these identities give the compliance entries in \eq{finite:eq:matrix}.

Taking the ratio of the displayed cross-circle difference with its conjugate gives the deviatoric kernel. With $|z|=|v|=1$,
\be\lb{phase}
\frac{z_1(z)-z_2(v)}{\overline{z_1(z)-z_2(v)}}=\frac{B(z)}{B(v)}K(\theta-\phi),\quad K(t)=-\rho+\frac{1-\rho^2}{\rho}\sum_{n\ge1}\rho^ne^{-\ri nt},\quad z=e^{\ri\theta},\ v=e^{\ri\phi}.
\ee
Let $\Pi_+,\Pi_-,\Pi_0$ be the Fourier projections in $L^2(\bS^1;\bC)$ onto the strictly positive, strictly negative and zero modes, respectively. They also act on periodic distributions by the same Fourier multipliers. For a zero-mean function $f$, put $f_\pm:=\Pi_\pm f$. Normalized convolution with the kernel in~\eq{phase} satisfies
\be\lb{commute}
B\big(K*(f/B)\big)=\frac{1-\rho^2}{\rho}\big(1+E(z)\big)f_-(z/\rho)+c\frac{r_{\rm b}B(z)}{(1+r_{\rm b}^2)B(z)+r_{\rm b}},\quad c:=(f/B)_0,
\ee
where $E(z):=B(z)/B(z/\rho)-1$. Indeed,
\[
\begin{aligned}
 \Pi_-(f/B)&=f_-/B+\frac{(1-r_{\rm b}^2)f_+(r_{\rm b})}{z-r_{\rm b}},
 &c&=\frac{1-r_{\rm b}^2}{r_{\rm b}}f_+(r_{\rm b}),\\
 E(z)&=-\frac{(r_{\rm b}+B(z))(1+r_{\rm b}B(z))}{(1+r_{\rm b}^2)B(z)+r_{\rm b}}.
\end{aligned}
\]
The evaluation at $r_{\rm b}$ produces the last term in \eq{commute}. This term is smooth in the physical coordinate $w_1=B(z)$ because $|(1+r_{\rm b}^2)w_1+r_{\rm b}|\ge3/4$. For $f=\bar\tau_2$, $c$ is minus the conjugate of the physical first coefficient of $p_2$. The reversed operator is the complex conjugate of the operator in~\eq{commute} applied to $\bar f$. Its coefficient is $(fB)_0$, which for $f=\bar\tau_1$ is the conjugate of the physical first coefficient of $p_1$.

Write
\[
m_i=\frac1{2\pi}\int_{\partial D_i}p_i(w)\bar w\,\rd S\in\bR,\qquad c_K=\frac1{2\mu(\kappa_0+1)}.
\]
Here $m_i$ is real by the zero torque condition. The finite-rank cross terms in \eq{commute} have coefficients $m_1,m_2$ and are uniformly smooth in the centered circle coordinates. After subtracting these terms, the cross terms that are not convolutions are
\be\lb{error1}
\varphi_1=E(z)(\Pi_-\bar\tau_2)(z/\rho),\qquad\varphi_2=\overline{E(z)}(\Pi_+\bar\tau_1)(\rho z).
\ee

For $n\ge1$ order the Fourier coefficients as $(\tau_{1,n},\tau_{2,n},\overline{\tau_{1,-n}},\overline{\tau_{2,-n}})$. Combining the logarithmic modes and the convolution term in \eq{commute}, the principal symbol of the Kelvin layer, in this ordering at frequency $n$, is $\kappa_0[2\mu(\kappa_0+1)n]^{-1}M_0(\eta n)$, where
\be\lb{matrix2}
M_0(y)=\begin{pmatrix}
1&e^{-y}&0&0\\
e^{-y}&1&b(y)&0\\
0&b(y)&1&e^{-y}\\
0&0&e^{-y}&1
\end{pmatrix},\qquad b(y)=\frac{2\sinh\eta}{\kappa_0\eta}ye^{-y}.
\ee
The nonzero-mode layer equation consists of this symbol, the single-circle first mode from \eq{circle2}, the smooth finite-rank cross terms, and $c_K(1-\rho^2)\Pi(\varphi_1,\varphi_2)/\rho$, where $\Pi=\Pi_++\Pi_-$. In particular the remainder is independent of the inclusion tensor.

Set $g(v)=\sum_{n\ge1}\overline{p_{2,-n}}v^n$. In the right physical coordinate $w=B(z)$ the remainder $\mathcal R_1=(1-\rho^2)\varphi_1/\rho$ is given by
\be\lb{sm:Rphysical}
\mathcal R_1(w)=c(w)g(v(w)),\quad v(w)=\frac{r_{\rm b}w}{(1+r_{\rm b}^2)w+r_{\rm b}},\quad c(w)=\frac{(1+r_{\rm b}^2)(r_{\rm b}+w)^2(1+r_{\rm b}w)^2}{r_{\rm b}w((1+r_{\rm b}^2)w+r_{\rm b})^2}.
\ee
In the exterior bipolar variable, $J(z)g(-1/B(z))$ is holomorphic and vanishes at infinity. The zero of $g$ at the origin cancels the pole of $J$. These are the strictly negative modes of $\bar\tau_2$: the opposite physical modes give an interior function after composition with the reversed coordinate, and its zero cancels the interior pole of $J$. Hence $v(w)=-1/B(z/\rho)$ and $(\Pi_-\bar\tau_2)(z/\rho)=J(z/\rho)g(v(w))$. Substitution in \eq{error1} gives the coefficient $c$. Reflection exchanges the two correction formulas.

With $Y=|1+w|$, $\Delta_w=\eta^2+Y^2$ and $D_w=(1+r_{\rm b}^2)w+r_{\rm b}$, the boundary identities
\be\lb{cross:depth}
1-|v|^2=\frac{(1+r_{\rm b}^2)((1-r_{\rm b})^2+r_{\rm b}Y^2)}{|D_w|^2},\qquad|1-v|^2=\frac{(1-r_{\rm b})^4+r_{\rm b}(1+r_{\rm b}^2-r_{\rm b})Y^2}{|D_w|^2}
\ee
give $c_*\Delta_w\le1-|v|\le C_*\Delta_w$, since $|D_w|\ge3/4$. All derivatives of $v$ are uniformly bounded on the circle, and $|c^{(j)}(w)|\le C_j\Delta_w^{\max\{2-j/2,0\}}$ follows by differentiating its four numerator factors.

Use the bipolar parameterization \eq{bipolar} and the ordering of the four Fourier channels in \eq{matrix2}. For $C_i$, first extend the positive nonexceptional multiplier to the first mode. Correct it on the real first mode by adding the smooth function
\[
\left(\frac{1}{2(\lambda_1+\mu_1-\lambda-\mu)}-\frac{\kappa_1}{2\mu(k-\kappa_1/\kappa_0)}\right)m_iw.
\]
The factor in parentheses is $O(\eta_A)$. The single-circle first term in \eq{circle2} and the two finite-rank cross terms from \eq{commute} also have uniformly bounded smooth norms, by \eq{finite:eq:weak-density}. The remaining self-interaction multipliers add to the diagonal of the Kelvin symbol, giving
\be\lb{finite:eq:matrix}
M_{\rm fin}(y)=\begin{pmatrix}
A_+&e^{-y}&0&0\\
e^{-y}&A_-&b(y)&0\\
0&b(y)&A_-&e^{-y}\\
0&0&e^{-y}&A_+
\end{pmatrix},\qquad y=n\eta,
\ee
where $b(y)$ is defined in \eq{matrix2}. The rigid symbol $M_\infty$ has $A_+=A_-=1$.

\section{Stiff disks and comparison with the rigid solution}\label{sec:pre}
We prove Theorem~\ref{thm:disks}\textup{(i)} by estimating the Fourier projections of the layer densities with a weight that depends on the gap and the material factors. The same estimate applies to the density difference after subtracting the rigid equation. Reconstruction of the layer gives the corresponding bounds for the displacement gradients.

\subsection{The bipolar flow and the physical gradient}
The weighted inverse below uses one bound for the denominator of $\Phi_t=B\circ(e^{-t}B^{-1})$. Let $|w|=1$, $|\arg t|\le\vartheta<\pi/2$, $R=e^{-t}$, $m=\min\{|t|,1\}$ and $Y=|1+w|$. Then
\be\lb{sax:flow-denominator}
 c(\eta+mY)\le |D(w)|\le C(\eta+mY),\qquad D(w)=1-r_{\rm b}^2R+r_{\rm b}(1-R)w .
\ee
Indeed $D(-1)=(1-r_{\rm b})(1+r_{\rm b}R)$ has modulus comparable to $\eta$ uniformly in the sector. For small $|t|$ this follows from $R=1+O(|t|)$, and away from zero $|R|\le e^{-c_\vartheta|t|}$. Since $|D(w)-D(-1)|\le CmY$, this proves the upper bound and the lower bound for $mY\le c_*\eta$ with $c_*$ small. For the remaining points use
\[
 D(w)=(1+r_{\rm b}w)\left(1-r_{\rm b}R\,\frac{r_{\rm b}+w}{1+r_{\rm b}w}\right).
\]
The quotient has modulus one, $|1+r_{\rm b}w|^2=(1-r_{\rm b})^2+r_{\rm b}Y^2$, and $1-r_{\rm b}|R|\ge c_\vartheta m$. Thus $|D(w)|\ge c mY$, giving \eq{sax:flow-denominator}.

\smallskip\noindent\emph{Reconstruction of the physical gradient.}\quad Let $P_+$ and $P_-$ be the positive and negative Fourier projections of $p$ in the centered circle coordinate. Thus $P_+$ is holomorphic inside the circle and $P_-$ is holomorphic outside and vanishes at infinity. For $|w|>1$ the Kelvin layer is
\be\lb{recon:exterior}
\cS p(w)=c_K\Big\{\kappa_0\sum_{n\ge1}\big(\frac{p_{-n}}{nw^n}+\frac{p_n}{n\bar w^n}\big)+\big(\frac{w}{\bar w}-\frac1{\bar w^2}\big)\overline{P_-(w)}-\frac{\bar p_1}{\bar w}\Big\}.
\ee
For $|w|<1$, the corresponding formula is
\[
\cS p(w)=c_K\Big\{\kappa_0\sum_{n\ge1}\frac{p_nw^n+p_{-n}\bar w^n}{n}+\big(\frac1{\bar w^2}-\frac{w}{\bar w}\big)\overline{P_+(w)}-\frac{\bar p_1}{\bar w}\Big\}.
\]
The apparent singularities of the interior expression at zero are removable. In both formulas the term involving the highest derivative of a projection retains the factor $|w|^2-1$ or $1-|w|^2$. For the gradient estimate the extra analytic derivative is therefore controlled by Cauchy's estimate on a circle of radius half the normal distance. The weight is comparable on such circles. A boundary weight $\sigma+|1+w|$ extends to the interior channels by the maximum principle after multiplication by $1+(1-\min\{\sigma,1/2\})w$. For exterior channels apply inversion, using the vanishing of $P_-$ at infinity. Division of an analytic projection by its zero at $w=0$ preserves these bounds by Cauchy's estimate on a fixed inner disk.

\subsection{A weighted inverse for the principal symbol}
\begin{lemma}[Weighted inverse of the stiff transmission symbol]\lb{lem:hardweighted}
Let $0<\eta\le\eta_0$ and $A_\pm\ge1$. Set
\[
D_A=\operatorname{diag}(A_+,A_-,A_-,A_+),\qquad \eta_A=\frac12\log(A_+A_-).
\]
Assume that $\eta_A\le1$ and
\[
\frac{\sinh\eta}{\kappa_0\eta}\sqrt{A_+/A_-}\le q<1
\]
for a fixed $q$. Let $M(y)=M_0(y)+D_A-I_4$. Put $\sigma=\eta+\eta_A$ and $L(w)=\sigma+|1+w|$ on the right circle, with the reflected weight on the left. Suppose each input channel is of the form $f(z)=zB'(z)A(B(z))$, with $A$ holomorphic in the disk and $|A|\le K$ there. Apply the multiplier $M(\eta n)^{-1}$ for $n\ge1$ and divide by $J=zB'/B$ to recover the physical density projections. Each projection is bounded by $CK/L(w)$, with boundary values understood almost everywhere for bounded holomorphic inputs.
\end{lemma}

\begin{proof}
Write $M(y)=D_A^{1/2}(I_4+T(y))D_A^{1/2}$. The three consecutive off-diagonal entries of $T(y)$ are $\rho_Ae^{-y}$, $\rho_A\beta ye^{-y}$, and $\rho_Ae^{-y}$, where
\[
 \rho_A=e^{-\eta_A},\qquad
 \beta=\frac{2\sinh\eta}{\kappa_0\eta}\sqrt{A_+/A_-}.
\]
Let $P(y)$ have the same tridiagonal form with off-diagonal entries $1,\beta y,1$. Reflection splits $P(y)$ into two $2\times2$ blocks, giving
\[
\|P(y)\|=\sqrt{1+(\beta y/2)^2}+\beta y/2\le e^{qy},\qquad y\ge0.
\]
Thus, for every frequency $y=n\eta>0$,
\be\lb{highseries}
M(y)^{-1}=D_A^{-1/2}\sum_{\ell=0}^\infty(-T(y))^\ell D_A^{-1/2},\quad\|T(y)\|\le\rho_Ae^{-(1-q)y}<1.
\ee

Write $T(y)^\ell=e^{-\ell y}\sum_{j=0}^\ell\bC_{\ell,j}y^j$. These coefficient matrices have nonnegative entries. Comparing them with $P(y)^\ell$ at any real $y>0$, and then taking $y=j/(q\ell)$, gives
\be\lb{highcoeff}
\|\bC_{\ell,j}\|\le4\rho_A^\ell(eq\ell/j)^j\quad(1\le j\le\ell),\quad\|\bC_{\ell,0}\|\le\rho_A^\ell.
\ee

We estimate these multipliers in the physical coordinate before summing them. For $\Re t>0$, take $R=e^{-t}$, $w=B(z)$ and $\Phi_t(w)=B(RB^{-1}(w))$. Then
\be\nn
\Phi_t(w)=\frac{(R-r_{\rm b}^2)w+r_{\rm b}(R-1)}{1-r_{\rm b}^2R+r_{\rm b}(1-R)w},\qquad\frac{f(Rz)}{J(z)}=w\Phi_t'(w)A(\Phi_t(w)).
\ee
For $q<\theta_0<1$, $|t-\ell\eta|\le\theta_0\ell\eta$ and $\ell\ge1$, the map $\Phi_t$ sends the closed disk into the open disk. Here $t$ stays in a fixed sector and $(1-\theta_0)\ell\eta\le |t|\le(1+\theta_0)\ell\eta$. Hence \eq{sax:flow-denominator} gives, with $Y=|1+w|$,
\be\nn
 |D(w)|\ge c\big(\eta+\min\{\ell\eta,1\}Y\big),\qquad |R|\le e^{-c\ell\eta}.
\ee
Using $\Phi_t'=R(1-r_{\rm b}^2)^2/D(w)^2$ and $c\eta\le1-r_{\rm b}^2\le C\eta$, we obtain
\be\lb{highkernel}
 |w\Phi_t'(w)A(\Phi_t(w))|\le CKe^{-c\ell\eta}(1+\ell Y)^{-2}.
\ee
When $\ell\eta\ge1$, use $1+\ell Y\le(1+\ell\eta)(1+Y/\eta)$ and absorb $(1+\ell\eta)^2$ into the exponential.

The multiplier $(\eta n)^je^{-\ell\eta n}$ is obtained by $(-\eta)^j\pa_t^j$ at $t=\ell\eta$. Cauchy's formula on the depth circle of radius $\theta_0\ell\eta$, together with the kernel bound~\eq{highkernel} and the coefficient bound~\eq{highcoeff}, introduces the factor
\[
(eq\ell/j)^j\frac{j!}{(\theta_0\ell)^j}\le C\sqrt{j+1}(q/\theta_0)^j.
\]
This factor is summable in $j$. Since $\rho_A^\ell=e^{-\ell\eta_A}$ and $\sigma=\eta+\eta_A$, summing the series in~\eq{highseries} reduces to the bound
\[
\sum_{\ell\ge1}\frac{e^{-c\sigma\ell}}{(1+\ell Y)^2}\le\frac{C}{\sigma+Y}.
\]
For $Y\le\sigma$ sum the geometric series, and for $Y>\sigma$ compare with the integral of $(1+tY)^{-2}$. The term $\ell=0$ is $D_A^{-1}$ applied to $wA(w)$. The maximum principle with the factor $1+(1-\min\{\sigma,1/2\})w$ gives the same weight throughout the disk. Reflection gives the assertion for the other circle.
\end{proof}

\subsection{The complete density estimate}
Choose $\kappa_0^{-1}<q<1$. The factors $A_\pm$ tend uniformly to one as $k_{\rm sh}^{-1}$ and $k_{\rm vol}^{-1}$ tend independently to zero. Thus Lemma~\ref{lem:hardweighted} applies for sufficiently large fixed contrast threshold and sufficiently small gaps. Write $L_{\rm fin}(w)=\eta+\eta_A+|1+w|$ on the right circle and use its reflected weight on the left. At rigidity $A_+=A_-=1$, so $\eta_A=0$ and the right-circle weight reduces to $L_0(w)=\eta+|1+w|$.

Consider \eq{finite:eq:boundary} with prescribed loads $F_i$ in place of $H_{\rm fin}-R_i$. Let $K_F$ bound the first derivatives of the two physical Fourier parts of each load. A uniform $C^3$ bound suffices, by summation of the differentiated Fourier series. For the physical Fourier parts $P_{i,\pm}$ of $p_i$, put
\[
 P=\max_{i,\pm}\|L_{\rm fin,i}P_{i,\pm}\|_\infty,\qquad
 H_p=\sum_i\|p_i\|_{H^{-1}(\Gamma_i)},
\]
where $L_{\rm fin,1}=L_{\rm fin}$ and $L_{\rm fin,2}$ is its reflection. The densities are smooth at each fixed positive gap. On the right circle $\partial_\theta F_+(B(z))=\ri zB'(z)F_+'(B(z))$. Reversal and conjugation give the other three channels. Multiplication by $n$ on the nonzero bipolar modes therefore gives inputs of the form in Lemma~\ref{lem:hardweighted}, with analytic factors bounded by $CK_F$. The finite-rank Kelvin terms and the exceptional first-mode compliance term have differentiated physical projections bounded by $CH_p$.

For the cross-circle remainder $c(w)g(v(w))$ in \eq{sm:Rphysical}, write $\Delta_w=\eta^2+|1+w|^2$. The identities \eq{cross:depth} give
\[
 c\Delta_w\le1-|v(w)|\le C\Delta_w,\qquad
 c L_{\rm fin}(w)\le\eta+\eta_A+|1-v(w)|\le C L_{\rm fin}(w).
\]
The second comparison persists on a circle of radius $c_*\Delta_w$ about $v(w)$ for fixed small $c_*$. The reflected weight bound gives $|g|\le CP/L_{\rm fin}(w)$ on this circle. Cauchy's estimate and the bounds for $c,c',v'$ following \eq{cross:depth} yield
\[
 |\partial_w(c(w)g(v(w)))|
 \le CP(\Delta_w^{3/2}+\Delta_w)/L_{\rm fin}(w).
\]
On $|1+w|\le a_0$ this is at most $Ca_0P$ for $\eta\le a_0$. Off this arc $v(w)$ lies in a fixed compact subdisk, where $H_p$ controls $g$ and its derivatives. The rational map $v$ has no pole on $|w|\ge1$ and maps this exterior into the unit disk by \eq{cross:depth} and the maximum principle. Since $c(w)=O(w)$ at infinity, the composition has at most one growing Laurent term. Its linear and constant coefficients are bounded by $CH_p$, by evaluation of $g,g'$ at $v(\infty)=r_{\rm b}/(1+r_{\rm b}^2)$. Subtracting these terms gives the exterior projection. The linear term gives a smooth interior input. The derivatives of both projections are thus bounded by $Ca_0P+C_{a_0}H_p$. Applying Lemma~\ref{lem:hardweighted} in all four channels gives $P\le Ca_0P+C_{a_0}(K_F+H_p)$, with $C$ independent of $a_0$. Fix $a_0$ so that $Ca_0<1/2$. We obtain
\be\label{finite:eq:complete-weight}
 P\le C(K_F+H_p).
\ee

\subsection{Gradient bounds in the matrix and the disks}
For \eq{finite:eq:boundary}, the cutoff potential and the rigid displacements give $K_F\le C\|\varphi\|_{C^{1,\gamma}}$. Estimate \eq{finite:eq:weak-density} bounds $H_p$ by the same quantity. Thus \eq{finite:eq:complete-weight} controls the density projections for the finite-stiffness problem.

Outside the disks the gradient of the Kelvin layer involves the density projections themselves. The additional analytic derivative in the reconstruction is multiplied by the normal distance to the circle and is controlled by Cauchy's estimate \eq{recon:exterior}. On a neck section the geometric weight is comparable to $s_\varepsilon+\eta_A+|t|$. Replacing $\eta_A$ by the comparable sum in \eq{finite:eq:compare} gives \eq{finite:eq:mat}. On the portion of a fixed collar of the disks outside the neck the facing-point distance is bounded below. The same layer estimates thus give $C\|\varphi\|_{C^{1,\gamma}}$ there. Away from the collars, interior estimates, with the uniform energy bound, give $C\|\varphi\|_{C^{1,\gamma}}$ on the rest of $\Omega_1\cap M_\varepsilon$.

For the rigid displacement in \eq{finite:eq:boundary}, whose coefficients are bounded by $C\|\varphi\|_{C^{1,\gamma}}$, the interior estimate including the disk-side traces reads
\be\lb{finite:eq:inc}
\|\nabla(u_{\rm fin}-R_i)\|_{L^\infty(D_i)}
 \le C\|\varphi\|_{C^{1,\gamma}}\,\frac{k_{\rm sh}^{-1}+k_{\rm vol}^{-1}}{\sqrt\varepsilon+k_{\rm sh}^{-1}+k_{\rm vol}^{-1}},\qquad i=1,2.
\ee
Inside the disks, the tangential derivative of the trace $C_ip_i$ has the additional factor $\eta_A$. The analytic primitives are obtained by division of the coefficients by $|n|$, with the real first mode included separately. Formula \eq{sildiskformula} with parameter $\kappa_1$ gives \eq{finite:eq:inc}: its correction has coefficient $1/\kappa_1\le1$, and the additional analytic derivative is multiplied by the distance to the circle. Using \eq{finite:eq:compare} exposes the same two inverse contrasts in \eq{finite:eq:inc}. The constants are uniform for $\kappa_1>1$ in the stated stiff regime.

\subsection{Comparison with the rigid solution}
Let $D_i^*$ be the touching disks. Following \cite[(2.8)--(2.9)]{LX25}, let $u^*$ minimize the matrix energy with outer trace $\varphi$ among displacements equal to a common rigid motion on the two touching disks. Let $b[\varphi]$ be the resultant matrix traction on $\partial D_1^*$, with the normal pointing out of $D_1^*$. Its normal and tangential components at the neck are $b_1[\varphi]=b[\varphi]\cdot e_1$ and $b_2[\varphi]=b[\varphi]\cdot e_2$, respectively.

\begin{corollary}[Comparison with the rigid solution]
In Theorem~\ref{thm:disks}\textup{(i)}, assume also that both coordinate reflections preserve $\Omega$ and that the outer trace $\varphi\in C^{2,\gamma}(\partial\Omega;\bR^2)$ is odd under $x\mapsto-x$. Then throughout $G_\varepsilon(r_0)$, including both matrix-side traces,
\be\lb{finite:eq:sima-expansion}
\begin{aligned}
\nabla u_{\rm fin}(x_1,t)&=\frac{\sqrt{\varepsilon}}{\pi \delta_\varepsilon(t)}\left(\frac{b_1[\varphi]}{\lambda+2\mu}e_1+\frac{b_2[\varphi]}{\mu}e_2\right)\otimes e_1+\mathcal E^\varepsilon(x_1,t),\\
|\mathcal E^\varepsilon(x_1,t)|&\le C\|\varphi\|_{C^{2,\gamma}}\left(1+\frac{\varepsilon^{1/2+\gamma/8}}{\delta_\varepsilon(t)}+\frac{1}{\sqrt{\varepsilon}+|t|}\min\!\left\{1,\frac{k_{\rm sh}^{-1}+k_{\rm vol}^{-1}}{\sqrt{\varepsilon}}\right\}\right).
\end{aligned}
\ee
The constant is uniform in both inclusion moduli.
\end{corollary}

The rigid expansion of Li and Xu \cite[Theorem~2.1(i)]{LX25} for these symmetries has the displayed leading term. Estimate \eq{finite:eq:diff} adds the finite-stiffness error at the same gap. If $k_{\rm sh}^{-1}+k_{\rm vol}^{-1}=o(\sqrt\varepsilon)$, the remainder in \eq{finite:eq:sima-expansion} is $o(\varepsilon^{-1/2})$ on $G_\varepsilon(r_0)\cap\{|t|\le L\sqrt\varepsilon\}$ for each fixed $L$. A nonzero contact traction component then gives a midpoint gradient of at least $c_\varphi\varepsilon^{-1/2}$. If $k_{\rm sh}^{-1}+k_{\rm vol}^{-1}=O(\varepsilon)$, the error from finite stiffness is bounded throughout $G_\varepsilon(r_0)$.

\subsubsection{Energy and density of the difference}
We prove~\eq{finite:eq:diff} for arbitrary outer data in Theorem~\ref{thm:disks}\textup{(i)}. Let $v=u_{\rm fin}-u_\infty$. Denote the finite-stiffness traces and densities by $w_{{\rm fin},i}=w_i$ and $p_{{\rm fin},i}=p_i$. Let $w_{\infty,i}$ denote the rigid trace and $p_{\infty,i}$ its matrix traction, with the normal pointing out of the disk. The circle formula for $C_ip_i$ and the first-order density estimate give
\[
\abs{\partial_\tau(w_i-R_i)}\le\frac{C\|\varphi\|_{C^{1,\gamma}}\eta_A}{s_\varepsilon+\eta_A+\operatorname{dist}_{\Gamma_i}(\,\cdot\,,P_i)},
\]
where $P_i$ is the facing point and the distance on the circle is arc length. After subtracting the value at $P_i$, the trace has $H^1$ norm at most $C\|\varphi\|_{C^{1,\gamma}}\eta_A/\sqrt{s_\varepsilon+\eta_A}$. Lift these two adjusted traces into the matrix by the trace construction used to prove~\eq{finite:eq:weak-density}. The difference $v$ minus the lift has rigid traces on both disks. By \eq{eq:common-rigid-lift},
\[
\norm{e(v)}_{L^2(M_\varepsilon)}\le C\|\varphi\|_{C^{1,\gamma}}\,\frac{\eta_A}{\sqrt{s_\varepsilon+\eta_A}}.
\]

In the disks the strain of $v$ equals that of the homogeneous inclusion extension of $w_i-R_i$. The $H^{1/2}$-to-$H^1$ bound for the fixed disk extension is uniform in $\kappa_1$. This follows from the Poisson formula with its additional term $(1-|z|^2)/\kappa_1$ times an analytic derivative \eq{sildiskformula}, by integrating the Fourier coefficients. Korn's inequality, using $v=0$ on $\partial\Omega$, now gives
\be\lb{finite:eq:negative-diff}
\norm{v}_{H^1(\Omega)}+\sum_i\norm{p_{{\rm fin},i}-p_{\infty,i}}_{H^{-1}(\Gamma_i)}\le C\|\varphi\|_{C^{1,\gamma}}\,\frac{\eta_A}{\sqrt{s_\varepsilon+\eta_A}}.
\ee
For the density difference use the matrix-side trace lifting and $\Lambda_0(w_{{\rm fin},i}-w_{\infty,i})$, as in \eq{finite:eq:weak-density}. The cutoff potential difference and the rigid and physical moment differences also obey this bound in their respective smooth norms.

\subsubsection{Pointwise comparison}
Subtract the two boundary systems \eq{finite:eq:boundary}, using the rigid operator on the density difference. The additional prescribed displacement is $-C_ip_{{\rm fin},i}$. Its differentiated physical projections are bounded by $C\|\varphi\|_{C^{1,\gamma}}\eta_A/(s_\varepsilon+\eta_A)$, since the nonexceptional differentiated multipliers are $O(\eta_A)$ and the first-mode correction is controlled by the weak density norm. The differences of the cutoff potentials and rigid displacements, as well as $H_p$ for the density difference, satisfy \eq{finite:eq:negative-diff}. All Kelvin remainders act on this density difference. Estimate \eq{finite:eq:complete-weight}, applied with the rigid weight $L_0$, gives a weighted density bound of $C\|\varphi\|_{C^{1,\gamma}}\eta_A/(s_\varepsilon+\eta_A)$. Reconstruction gives
\[
 |\nabla(u_{\rm fin}-u_\infty)(x_1,t)|
 \le C\|\varphi\|_{C^{1,\gamma}}\,\frac{\eta_A}{(s_\varepsilon+\eta_A)(s_\varepsilon+|t|)} .
\]
Since $\frac12\min\{1,\eta_A/s_\varepsilon\}\le \eta_A/(s_\varepsilon+\eta_A)\le\min\{1,\eta_A/s_\varepsilon\}$, \eq{finite:eq:compare} proves \eq{finite:eq:diff}.

\section{The cavity profile and the comparison for soft disks}\label{sec:soft-disk-proof}
The cavity strain determines the symmetric part of the singular gradient. To estimate the energy of a disk extension, we also need the tangential derivative of displacement on each circle, which includes its rotation. We first identify both terms of the cavity gradient from the stress resultants. We then estimate the extension energy and use the transmitted traction to obtain the pointwise comparison in Theorem~\ref{thm:disks}\textup{(ii)}. Throughout the soft comparison $0<k_{\rm sh},k_{\rm vol}<1$. Write $\delta=\delta_\varepsilon$, $H(t)=\varepsilon+t^2$, and $d(t)=s_\varepsilon+|t|$.

Let $\sigma_{\rm cav}=\bC_0e(u_{\rm cav})$ and $I(0)=(-\varepsilon/2,\varepsilon/2)$. The longitudinal and signed shear resultants at the middle section are

\be\lb{hole:eq:forces}
N_\varepsilon=\int_{I(0)}(\sigma_{\rm cav})_{22}(x_1,0)\,\rd x_1,\qquad Q_\varepsilon=-\int_{I(0)}(\sigma_{\rm cav})_{12}(x_1,0)\,\rd x_1.
\ee
Let $J_0=\left(\begin{smallmatrix}0&-1\\1&0\end{smallmatrix}\right)$ be the counterclockwise quarter-turn. Define the profile
\be\nn
\mathcal P_\varepsilon(t)=\frac{N_\varepsilon}{H(t)}\bC_0^{-1}(e_2\otimes e_2)+\frac{3(\lambda+2\mu)Q_\varepsilon}{4\mu(\lambda+\mu)H(t)^2}J_0.
\ee

\begin{prop}[Cavity estimates for the soft comparison]\label{thm:disk-cavity}
For all sufficiently small gaps, the sectional resultants satisfy $|N_\varepsilon|\le C\|\varphi\|_{C^{1,\gamma}}\sqrt\varepsilon$ and $|Q_\varepsilon|\le C\|\varphi\|_{C^{1,\gamma}}\varepsilon^{3/2}$. On each complete section, up to both circles,
\be\lb{hole:eq:profile}
 \nabla u_{\rm cav}(x_1,t)=\mathcal P_\varepsilon(t)
       +\mathcal R_{\rm cav}(x_1,t),\qquad
 |\mathcal R_{\rm cav}(x_1,t)|\le C\|\varphi\|_{C^{1,\gamma}}.
\ee
Consequently
\be\lb{hole:eq:gradient}
 |\nabla u_{\rm cav}(x_1,t)|\le
             \frac{C\|\varphi\|_{C^{1,\gamma}}}{\sqrt\varepsilon+|t|}\qquad ((x_1,t)\in G_\varepsilon(r_0)).
\ee
The gradient is bounded by $C\|\varphi\|_{C^{1,\gamma}}$ on $(\Omega_1\cap M_\varepsilon)\setminus G_\varepsilon(r_0)$.

For a fixed circular outer domain enclosing the disks, there exist a fixed trace $\varphi_*\in C^{1,\gamma}$ and a sequence $\varepsilon_j\downarrow0$ for which
\be\lb{hole:eq:sharp}
|e(u_{\rm cav}(\cdot\,;\varepsilon_j))(0)|\ge c\,\varepsilon_j^{-1/2}.
\ee
\end{prop}

The Airy potential gives the sectional force bounds. Subtracting an interpolant with the same boundary values leaves a clamped biharmonic problem, whose Hessian we estimate on the scale of the local gap. The rotation identity on the center line then gives the sharper bending-moment bound needed for the full gradient profile.

\subsection{Airy data and resultant estimates}
The cavity energy bound from Section~\ref{sec:common-variational} gives
\be\lb{hole:eq:energy}
\norm{\sigma_{\rm cav}}_{L^2(M_\varepsilon)}\le C\|\varphi\|_{C^{1,\gamma}}.
\ee

Away from contact, a finite chain of overlapping fixed-size patches joins each exterior patch to the outer collar. Korn's inequality anchored on that collar, followed by interior or traction-free boundary regularity, bounds the displacement and gradient there by $C\|\varphi\|_{C^{1,\gamma}}$. In particular the rotations at $(0,\pm r_0)$ are bounded.

Extend $\sigma_{\rm cav}$ by zero into the disks. It is symmetric and distributionally divergence-free, and hence has an Airy potential on a simply connected neighborhood $U\Subset\Omega$ of the disks and neck. One may take the union of slightly enlarged disks, which overlap for small gaps. Integrating the two divergence-free rows and then using symmetry gives
\be\nn
 \sigma_{\rm cav}=\operatorname{cof}D^2\Phi
 =\begin{pmatrix}\Phi_{tt}&-\Phi_{x_1t}\\-\Phi_{x_1t}&\Phi_{x_1x_1}\end{pmatrix},
 \qquad \Phi\in H^2_{\rm loc}(U).
\ee
Strain compatibility gives $\Delta^2\Phi=0$ in the matrix. The potential is affine in each disk. After subtraction of one affine function it equals zero on the left and $Nx_1+Qt+M$ on the right. Write $I(t)=(-\delta(t)/2,\delta(t)/2)$. The endpoint traces of $\Phi,\nabla\Phi$ give
\be\lb{hole:eq:moments}
 \int_{I(t)}\Phi_{x_1x_1}\,\rd x_1=N,\qquad
 \int_{I(t)}x_1\Phi_{x_1x_1}\,\rd x_1=-(Qt+M).
\ee
Also $\int_{I(t)}\Phi_{x_1t}\,\rd x_1=Q$, so $N=N_\varepsilon$ and $Q=Q_\varepsilon$ by \eq{hole:eq:forces}.

With $z=x_1/\delta(t)+1/2$ on the neck, put
\[
 A_N=\tfrac12 \delta(t)z^2,\qquad A_M=3z^2-2z^3,\qquad A_Q=t A_M .
\]
Their values and gradients match the zero affine field on the left disk and $x_1,1,t$, respectively, on the right. Use $A_N,A_Q$ as tests on $G_\varepsilon(r_*)$ with fixed $r_*>2r_0$. Extend them affinely into the disks. In a fixed band at either end of the neck the walls have positive separation, so a partition of unity patches the interpolants to smooth extensions of these affine data. An outer cutoff gives compactly supported $H^2$ tests with uniformly bounded Hessians off $G_\varepsilon(r_*)$.

For each such test $A$, $\operatorname{cof}D^2 A$ is divergence-free and vanishes in the disks. Its pairing with the strain of any $H^1$ extension of the cavity displacement is zero. Isotropic compliance and two integrations by parts give
\be\lb{hole:eq:airytest}
 0=\int_U\Phi_{x_1x_1}A_{x_1x_1}\,\rd x_1\,\rd t+\int_U(2\Phi_{x_1x_1}+\Phi_{tt})A_{tt}\,\rd x_1\,\rd t.
\ee
Indeed $\int_U D^2\Phi:D^2A=\int_U\Delta\Phi\,\Delta A$ and $\int_U\Phi_{x_1t}A_{x_1t}=\int_U\Phi_{x_1x_1}A_{tt}$ by compact support. The compliance pairing is $(\lambda+2\mu)/(4\mu(\lambda+\mu))$ times $\int D^2\Phi:D^2A$.

On $G_\varepsilon(r_*)$ the transverse second derivatives of $A_N,A_Q$ are $\delta^{-1}$ and $-12tx_1\delta^{-3}$. By \eq{hole:eq:moments} and the evenness of $\delta$, the leading integrals in \eq{hole:eq:airytest} are $N\int_{-r_*}^{r_*}\delta^{-1}\,\rd t$ and $12Q\int_{-r_*}^{r_*} t^2\delta^{-3}\,\rd t$, with weights comparable to $s_\varepsilon^{-1}$ and $s_\varepsilon^{-3}$. The $tt$ derivatives are bounded by $C$ and $C/(s_\varepsilon+|t|)$, respectively, so both have bounded $L^2$ norm since $|I(t)|=\delta(t)\le C(s_\varepsilon+|t|)^2$. Equation~\eq{hole:eq:energy} controls the error integrals by $C\|\varphi\|_{C^{1,\gamma}}$. For a preliminary bound on $M$, Cauchy--Schwarz in the second identity of \eq{hole:eq:moments} gives
\[
 \int_{-r_*}^{r_*}\frac{(Qt+M)^2}{\delta(t)^3}\,\rd t
 \le \tfrac1{12}\|\Phi_{x_1x_1}\|_{L^2(G_\varepsilon(r_*))}^2 .
\]
The mixed term integrates to zero and $cs_\varepsilon^{-5}\le\int \delta^{-3}\,\rd t\le Cs_\varepsilon^{-5}$. Thus
\be\lb{hole:eq:resultants}
 |N|\le C\|\varphi\|_{C^{1,\gamma}}s_\varepsilon,\qquad |Q|\le C\|\varphi\|_{C^{1,\gamma}}s_\varepsilon^3,\qquad |M|\le C\|\varphi\|_{C^{1,\gamma}}s_\varepsilon^{5/2}.
\ee

\subsection{The clamped remainder}
On $G_\varepsilon(r_*)$ let $A=NA_N+QA_Q+MA_M$ and $W=\Phi-A$. The values and first derivatives of $\Phi$ and $A$ agree on the two faces. Thus $W=\nabla W=0$ on the two facing arcs and $\Delta^2W=f:=-\Delta^2A$. Fix $|t_0|<2r_0$, where $2r_0<r_*$, and set
\[
d_0=d(t_0),\qquad \delta_0=\delta(t_0),\qquad S_0=\frac{|N|}{\delta_0}+\frac{|M|+|Q|d_0}{\delta_0^2}.
\]
In the band $|t-t_0|<cd_0$, the thickness is between $c_1\delta_0$ and $C_1\delta_0$, and $\delta_0$ is between $c_2d_0^2$ and $C_2d_0^2$. Each $x_1$ derivative of the interpolant costs $C/\delta_0$; each $t$ derivative costs $C/d_0$, in the following precise form for the finite orders used below:
\be\lb{hole:eq:interpolant}
|\partial_{x_1}^{\,a}\partial_t^{\,b}A|\le C_{a,b}S_0\delta_0^{2-a}d_0^{-b}.
\ee
This follows by differentiating the three displayed polynomials and using the bounds $|\partial_t^j\delta|\le C_jd_0^{2-j}$ for $j\ge1$. Since $A$ is cubic in $x_1$, $\partial_{x_1}^4A=0$; hence $|f|+\delta_0|\nabla f|\le CS_0/d_0^2$ in that band.

For an interval $J$ in this band, let $E(J)$ be the integral of $|D^2W|^2$ over the corresponding sections. Enlarge $J$ by $C\delta_0$ at each end to form $J^+$, and choose a cutoff $\chi(t)$ equal to one on $J$ and supported in $J^+$, with $|\chi^{(j)}|\le C_j\delta_0^{-j}$ for $j=1,2$. On each section, clamping gives $\|\nabla W\|_2\le C\delta_0\|D^2W\|_2$ and $\|W\|_2\le C\delta_0^2\|D^2W\|_2$. Testing with $\chi^4W$, the cutoff terms are bounded by $\tfrac14\|\chi^2D^2W\|_2^2+CE(J^+\setminus J)$. The forcing contributes at most another quarter of the same squared norm plus the source term below. Rearrangement gives
\be\lb{hole:eq:contraction}
 E(J)\le\theta E(J^+)+C\delta_0^4\int_{\{t\in J^+,\,x_1\in I(t)\}}|f|^2,\qquad 0<\theta<1.
\ee
Starting from $J_{\rm in}=(t_0-C\delta_0,t_0+C\delta_0)$, we can enlarge at least $c/d_0$ times. The area after $j$ steps is $O((j+1)\delta_0^2)$. The terminal energy is $O(\|\varphi\|_{C^{1,\gamma}}^2)$ by \eq{hole:eq:energy} for $\Phi$ and by \eq{hole:eq:interpolant} for $A$, whose squared Hessian norm is at most $C(N^2s_\varepsilon^{-1}+Q^2s_\varepsilon^{-3}+M^2s_\varepsilon^{-5})\le C\|\varphi\|_{C^{1,\gamma}}^2$. Summing the geometric series in \eq{hole:eq:contraction} yields
\[
 E(J_{\rm in})\le Ce^{-c/d_0}\|\varphi\|_{C^{1,\gamma}}^2+C\delta_0^6(S_0/d_0^2)^2.
\]
After rescaling by $\delta_0$, each point has an interior patch or a uniformly smooth patch meeting one clamped wall. The local biharmonic estimate \cite{ADN59}, with the displacement norm controlled by the transverse inequalities, bounds $|D^2W|+\delta_0|D^3W|$ by
\[
 C\delta_0^{-1}E(J_{\rm in})^{1/2}+C\delta_0^2\big(\|f\|_\infty+\delta_0^\gamma[f]_{C^\gamma}\big).
\]
Our bound on $f,\nabla f$ therefore gives, up to both faces,

\be\lb{hole:eq:W}
|D^2W(x_1,t_0)|+\delta_0|D^3W(x_1,t_0)|\le Cd_0^2S_0+C\delta_0^{-1}e^{-c/d_0}\|\varphi\|_{C^{1,\gamma}}.
\ee

\subsection{The strain and rotation on a complete section}
For the cavity displacement, write
\[
 \omega=\tfrac12(\partial_{x_1}u_2-\partial_tu_1),\qquad S(t)=|N|/\delta+(|M|+|Q|d(t))/\delta^2.
\]
Put $a_0=(\lambda+2\mu)/(4\mu(\lambda+\mu))$. The identity $\partial_t\omega=\partial_{x_1}e_{22}-\partial_t e_{12}$ and the compliance formula give, on the center line,
\[
 \partial_t\omega(0,t)
 =a_0(\Phi_{x_1x_1x_1}+\Phi_{x_1tt})(0,t)
 =-12a_0(Qt+M)\delta(t)^{-3}+O(S(t)+\|\varphi\|_{C^{1,\gamma}}).
\]
Here $A_{x_1x_1x_1}=-12(Qt+M)\delta^{-3}$, while \eq{hole:eq:interpolant} bounds $A_{x_1tt}$ and \eq{hole:eq:W} bounds the third derivatives of $W$. By \eq{hole:eq:resultants},
\[
 \int_{-r_0}^{r_0} S(t)\,\rd t
 \le C\big(|N|s_\varepsilon^{-1}+|Q|s_\varepsilon^{-2}+|M|s_\varepsilon^{-3}\big)
 \le C\big(\|\varphi\|_{C^{1,\gamma}}+|M|s_\varepsilon^{-3}\big).
\]
Integrate the rotation identity from $-r_0$ to $r_0$. The endpoint rotations are bounded by $C\|\varphi\|_{C^{1,\gamma}}$, the $Qt$ term integrates to zero, and $\int_{-r_0}^{r_0} \delta^{-3}\,\rd t\ge cs_\varepsilon^{-5}$. It follows that $|M|s_\varepsilon^{-5}\le C(\|\varphi\|_{C^{1,\gamma}}+|M|s_\varepsilon^{-3})$. Absorbing the last term for small $s_\varepsilon$ gives
\be\lb{hole:eq:moment-improved}
 |M|\le C\|\varphi\|_{C^{1,\gamma}}s_\varepsilon^5.
\ee
In particular $\int_{-r_0}^{r_0} S(t)\,\rd t\le C\|\varphi\|_{C^{1,\gamma}}$ and $d(t) S(t)\le C\|\varphi\|_{C^{1,\gamma}}$. Integrating the same identity from $t$ to $r_0$ now yields
\[
 \omega(0,t)=12a_0Q\int_t^{r_0}\frac{q}{\delta(q)^3}\,\rd q+O(\|\varphi\|_{C^{1,\gamma}}).
\]

The Hessian of the interpolant obeys
\[
A_{x_1x_1}=\frac{N}{\delta}-12(Qt+M)x_1\delta^{-3},\qquad|A_{x_1t}|\le Cd(t) S(t),\qquad|A_{tt}|\le Cd(t)^2S(t).
\]
In the first formula the term after $N/\delta$ is bounded by $C(|Q||t|+|M|)/\delta^2\le C\|\varphi\|_{C^{1,\gamma}}$, by \eq{hole:eq:resultants} and \eq{hole:eq:moment-improved}. The bound $d(t) S(t)\le C\|\varphi\|_{C^{1,\gamma}}$ and \eq{hole:eq:W} control the remaining stress terms by $C\|\varphi\|_{C^{1,\gamma}}$. Thus
\[
e(u_{\rm cav})=\frac{N}{\delta(t)}\bC_0^{-1}(e_2\otimes e_2)+O(\|\varphi\|_{C^{1,\gamma}})
\]
uniformly on every complete neck section.

To recover the rotation across the section use $\partial_{x_1}\omega=\partial_{x_1}e_{12}-\partial_te_{11}$. Each third derivative of $A$ here has a $t$ derivative, so \eq{hole:eq:interpolant} bounds it by $CS(t)/d(t)$. The corresponding derivatives of $W$ are bounded by $C(S(t)+\|\varphi\|_{C^{1,\gamma}})$ by \eq{hole:eq:W}. Across a width $\delta(t)\le Cd(t)^2$ this changes the rotation by at most $C\|\varphi\|_{C^{1,\gamma}}$.

Finally, $\delta(t)=\varepsilon+t^2+O(t^4)$. Replacing $\delta$ by $\varepsilon+t^2$ in the strain term costs $O(\|\varphi\|_{C^{1,\gamma}}s_\varepsilon)$. In the rotation integral the error is bounded by
\[
C|Q|\!\int_{-r_0}^{r_0}\frac{|q|^5}{(\varepsilon+q^2)^4}\,\rd q\le C|Q|s_\varepsilon^{-2}\le C\|\varphi\|_{C^{1,\gamma}}s_\varepsilon.
\]
The antiderivative of $q(\varepsilon+q^2)^{-3}$ is $-\tfrac14(\varepsilon+q^2)^{-2}$, and the endpoint term at $r_0$ is bounded. Adding the strain and the skew parts now proves \eq{hole:eq:profile}.

The bounds for $N$ and $Q$ give $|N|/(\varepsilon+t^2)+|Q|/(\varepsilon+t^2)^2\le C\|\varphi\|_{C^{1,\gamma}}/(s_\varepsilon+|t|)$. The expansion therefore gives \eq{hole:eq:gradient} on both faces as well as in the interior.

For the lower bound, take the full-plane solution of Lim and Yu under uniform normal loading. Its midpoint stress has order $s_\varepsilon^{-1}$ and its stress on a fixed annulus enclosing the holes is uniformly bounded \cite[Theorem~2.1 and Corollary~2.2]{LY17}. After subtraction of a rigid motion on the annulus, Korn's inequality and interior regularity bound its trace on a fixed circle in $C^{2,\gamma}$. Along a subsequence the traces converge in $C^{1,\gamma}$ to a single trace $\varphi_*$. At the corresponding gaps, applying \eq{hole:eq:gradient} to the difference of the two boundary data changes the midpoint strain by $o(s_\varepsilon^{-1})$. This gives \eq{hole:eq:sharp}.

\subsection{The energy comparison for soft disks}
By \eq{eq:common-soft-lift} we need only estimate the energy of a disk extension of the cavity trace $w_i=u_{\rm cav}|_{\partial D_i}$. On the facing arc \eq{hole:eq:profile} gives
\[
 |\partial_\tau w_i|\le C\|\varphi\|_{C^{1,\gamma}}\left(1+\frac{s_\varepsilon}{s_\varepsilon^2+t^2}+\frac{s_\varepsilon^3}{(s_\varepsilon^2+t^2)^2}\right).
\]
Both kernels have bounded integrals in $t$ and supremum $O(s_\varepsilon^{-1})$. Since arclength is comparable to $\rd t$ and the far arc has a uniform displacement and gradient bound, integration around the circle gives
\[
 \|w_i\|_\infty+\|\partial_\tau w_i\|_1\le C\|\varphi\|_{C^{1,\gamma}},\qquad \|\partial_\tau w_i\|_2^2\le C\|\varphi\|_{C^{1,\gamma}}^2/s_\varepsilon.
\]
Write $w_n$ for the coefficients of either trace $w=w_i$. The $L^1$ derivative bound gives $|nw_n|\le C\|\varphi\|_{C^{1,\gamma}}$. For $1\le |n|\le s_\varepsilon^{-1}$ this yields $\sum |n||w_n|^2\le C\|\varphi\|_{C^{1,\gamma}}^2(1+|\log s_\varepsilon|)$. Above $s_\varepsilon^{-1}$ the sum is bounded by $Cs_\varepsilon\|\partial_\tau w\|_2^2$, by Parseval. Including the mean gives

\be\nn
 \sum_i\|w_i\|_{H^{1/2}(\partial D_i)}^2\le C\|\varphi\|_{C^{1,\gamma}}^2(1+|\log\varepsilon|).
\ee

Take the homogeneous inclusion extension of $w_i$. Its energy is the real boundary pairing $\langle\Lambda_1w_i,w_i\rangle$. With $K_1=\lambda_1+\mu_1$, the disk formula \eq{bsw:diskmap} gives
\[
 \frac{\langle\Lambda_1w,w\rangle}{2\pi}
 =\frac{2\mu_1}{\kappa_1}\sum_{n\ge2}n|w_n|^2
  +2\mu_1\sum_{n\le-1}|n|\,|w_n|^2+2K_1(\Re w_1)^2 .
\]
Since $\kappa_1>1$, the first two terms are at most $C\mu_1\|w\|_{H^{1/2}}^2$. The real first coefficient is the mean normal displacement and satisfies $|\Re w_1|\le\|w\|_\infty\le C\|\varphi\|_{C^{1,\gamma}}$. Summation over the disks and \eq{eq:common-soft-lift} therefore bound the squared matrix strain error by $C\|\varphi\|_{C^{1,\gamma}}^2((1+|\log\varepsilon|)k_{\rm sh}+k_{\rm vol})$. The outer-collar competitor of Section~\ref{sec:common-variational} also bounds this error by $C\|\varphi\|_{C^{1,\gamma}}^2$, proving \eq{soft:eq:energycompare}. At fixed positive gap, Korn's inequality then gives convergence in $H^1$ in the matrix as both moduli tend to zero.

\subsection{The matrix gradient and the transmitted traction}
We first give an estimate for the background equation in $M_\varepsilon$ with zero outer trace and disk tractions $g_i\in C^{1/2}(\partial D_i)$, oriented by the disk normals. In this paragraph denote a solution by $v$, and set
\[
E=\|e(v)\|_{L^2(M_\varepsilon)},\qquad G=\max_i\|g_i\|_{C^{1/2}(\partial D_i)}.
\]

At $(0,\pm r_*)$ the anchored Korn inequality and interior regularity give $|\nabla v|\le CE$. Along each half of the center line the matrix contains balls of radius $c(\varepsilon+t^2)$. Lemma~\ref{lem:neck-line} with $n=2$ yields $|\nabla v(0,t)|\le CEd(t)^{-5/2}$.

On a full-section patch of length comparable to $\delta(t)$, scaled Korn's inequality gives a rigid motion $q$ with $\|v-q\|_2\le C\delta(t)E$. The scaled patches have uniformly Lipschitz geometry. The interior gradient estimate at the center bounds $|\nabla(v-q)(0,t)|$ by $CE/\delta(t)$. Combining this estimate with the center-line bound gives $|\nabla q|$ at most $CE(d(t)^{-5/2}+\delta(t)^{-1})$. Each point on the section is covered either by an interior patch or by a patch meeting one smooth wall, of size a fixed fraction of $\delta(t)$. The local traction estimate for the Lam\'e system \cite{ADN64} gives
\[
\|\nabla(v-q)\|_{\infty}+\delta(t)^{1/2}[\nabla(v-q)]_{C^{1/2}}\le C\big(\delta(t)^{-2}\|v-q\|_2+G\big)\le C\big(E/\delta(t)+G\big)
\]
on the smaller patches. Subtracting the rigid motion leaves the traction unchanged. On fixed patches away from the neck the bound is $C(E+G)$, using the same anchored Korn argument and the local boundary estimates. Consequently,
\be\lb{soft:eq:backgroundpoint}
|\nabla v(x_1,t)|\le C\big(Ed(t)^{-5/2}+G\big)
\ee
on the complete sections, including the matrix traces.

For a circle field $f$ define $\|f\|_{\varepsilon}=\|f\|_\infty+\varepsilon^{1/2}[f]_{C^{1/2}}$. The section estimates give
\be\lb{soft:eq:scaledtrace}
\max_i\|\partial_\tau v|_{\partial D_i}\|_{\varepsilon}\le C(Es_\varepsilon^{-5/2}+G).
\ee
For pairs of boundary points at arclength distance less than $c\,\varepsilon$ use the local seminorm estimate on a common patch. For larger distances use the supremum bound. Multiplication by the unit tangent preserves the estimate. For the cavity solution, \eq{hole:eq:gradient} and the fixed-patch bound give $\|\nabla u_{\rm cav}\|_\infty\le C\|\varphi\|_{C^{1,\gamma}}/s_\varepsilon$ on the disks' matrix-side traces and nearby patches. Subtract a constant displacement on each patch and apply the homogeneous traction estimate to obtain the same $C\|\varphi\|_{C^{1,\gamma}}/s_\varepsilon$ bound in the scaled trace norm.

Now take $v=u_{\rm fin}-u_{\rm cav}$. Its disk traction is $g_i=\Lambda_1(u_{\rm fin}|_{\partial D_i})$. By \eq{bsw:diskmap} the nonexceptional terms of $\Lambda_1$ are Fourier projections of the tangential derivative with coefficients bounded by $C\mu_1$, and the real first mode has coefficient $2K_1$. These projections are bounded on $C^{1/2}$ by the periodic Hilbert transform estimate. Passing from the scaled trace norm to $C^{1/2}$ costs at most $C/s_\varepsilon$. Thus \eq{soft:eq:scaledtrace} and the cavity trace bound give
\[
 G\le\frac{Cd_{\rm s}}{s_\varepsilon}\big(\|\varphi\|_{C^{1,\gamma}}/s_\varepsilon+Es_\varepsilon^{-5/2}+G\big).
\]

At each positive gap and positive stiffness these norms are finite by interface regularity. Choose $c_0$ so that $Cd_{\rm s}/s_\varepsilon\le1/2$ whenever $d_{\rm s}\le c_0s_\varepsilon$. Absorption, followed by \eq{soft:eq:energycompare}, gives
\[
G\le C\|\varphi\|_{C^{1,\gamma}}\frac{d_{\rm s}}{s_\varepsilon}\big(s_\varepsilon^{-1}+\mathcal A_\varepsilon s_\varepsilon^{-5/2}\big).
\]
Equation~\eq{soft:eq:backgroundpoint} and its fixed-patch counterpart give \eq{soft:eq:grad-difference}. Indeed $d(t)\ge s_\varepsilon$, $d_{\rm s}/s_\varepsilon\le c_0$, and $d_{\rm s} s_\varepsilon^{1/2}\le\mathcal A_\varepsilon$ since $d_{\rm s}<1$ and $d_{\rm s}\le\mathcal A_\varepsilon^2$.

\subsection{Transfer of the cavity rate and profile}
\begin{corollary}[Inheritance of the circular cavity profile]\label{cor:disk-soft-transfer}
Let $\mathcal A_\varepsilon$ be as in \eq{soft:eq:order}. Fix $M>0$ and suppose
\be\lb{soft:eq:smallmoduli}
 \mathcal A_\varepsilon^2\le M\varepsilon^{3/2}.
\ee
For all sufficiently small gaps,
\be\lb{soft:eq:gradient}
 |\nabla u_{\rm fin}(x_1,t)|\le
 \frac{C_M\|\varphi\|_{C^{1,\gamma}}}{\sqrt\varepsilon+|t|}\qquad ((x_1,t)\in G_\varepsilon(r_0)).
\ee
Outside $G_\varepsilon(r_0)$ the matrix gradient is bounded by $C_M\|\varphi\|_{C^{1,\gamma}}$ on $\Omega_1\cap M_\varepsilon$. If $\mathcal A_\varepsilon^2=o(\varepsilon^{3/2})$, then for fixed data and each fixed $L>0$,
\[
 \sup_{G_\varepsilon(r_0)\cap\{|t|\le L\sqrt\varepsilon\}}
 |\nabla u_{\rm fin}-\mathcal P_\varepsilon(t)|=o(\varepsilon^{-1/2}).
\]

For the circular outer domain, fixed trace and gap sequence of Proposition~\ref{thm:disk-cavity}, the midpoint strain satisfies
\[
 |e(u_{\rm fin}(\cdot\,;\varepsilon_j))(0)|\ge c\,\varepsilon_j^{-1/2}\qquad\text{for all sufficiently large }j.
\]
\end{corollary}

\begin{proof}
Under \eq{soft:eq:smallmoduli}, $\mathcal A_\varepsilon\le\sqrt M\,s_\varepsilon^{3/2}$ and $d_{\rm s}\le Ms_\varepsilon^3$. The preceding traction estimate gives $G\le C_M\|\varphi\|_{C^{1,\gamma}}d_{\rm s}/s_\varepsilon^2$, so \eq{soft:eq:backgroundpoint} and \eq{soft:eq:energycompare} yield
\[
 |\nabla(u_{\rm fin}-u_{\rm cav})(x_1,t)|
 \le C_M\|\varphi\|_{C^{1,\gamma}}\big(\mathcal A_\varepsilon d(t)^{-5/2}+d_{\rm s}/s_\varepsilon^2\big).
\]
Since $d(t)\ge s_\varepsilon$, these terms are $O_M(d(t)^{-1})$ and $O_M(s_\varepsilon)$, respectively. Adding \eq{hole:eq:gradient} proves \eq{soft:eq:gradient}, and the fixed-patch bound $C(E+G)$ gives the estimate off the neck. If $\mathcal A_\varepsilon^2=o(s_\varepsilon^3)$, the right-hand side is $o(s_\varepsilon^{-1})$ uniformly on the neck. Equation~\eq{hole:eq:profile} gives the stated profile comparison, with a bounded error if $\mathcal A_\varepsilon^2=O(s_\varepsilon^5)$.

The lower bound persists already under \eq{soft:eq:smallmoduli}. In the matrix ball of radius $\varepsilon/4$ at the origin, subtract a rigid displacement by scaled Korn and apply the interior gradient estimate. This gives
\[
|e(u_{\rm fin}-u_{\rm cav})(0)|\le c\,\varepsilon^{-1}E\le C_M\|\varphi\|_{C^{1,\gamma}}s_\varepsilon^{-1/2}=o(s_\varepsilon^{-1}).
\]
Apply \eq{hole:eq:sharp} for the same trace at the same gaps to conclude.
\end{proof}

\appendix
\section{Local bending regularity through the free faces}\label{sec:layer-proof}
We prove Lemma~\ref{sph:lem:odd-local} by tangential differentiation on the layer. The differential operators in \eq{sph:eq:bending-derivative} preserve the Kirchhoff--Love form of the leading displacement and differentiate along the free faces. This permits commutation in the weak equation without face terms. A local energy iteration controls these derivatives down to the thickness scale. The conormal identity \eq{layer:eq:conormal-recovery} then recovers the normal derivatives. For related estimates on thin elastic layers, see \cite{John65,Berger73,Monneau04,Monneau06}.

\begin{proof}[Proof of Lemma~\ref{sph:lem:odd-local}]
For a lateral subdomain $G\subset D$, write
\[
\|\Phi\|_{h,G}^2=h^{-1}\int_{\mathcal P_h(G)}|\Phi|^2\,\rd y\,\rd Z,\qquad\int_hf=h^{-1}\int_{\mathcal P_h(D)}f\,\rd y\,\rd Z.
\]
We extend laterally supported integrands by zero to $\mathcal P_h(D)$. These norms are uniformly equivalent to the $L^2$ norms of $\Phi(y,ha(y)t)$ on the fixed cylinder. The same normalization is used in $a_h(V,\psi)=\int_h\bC(y)e(V):e(\psi)$.

\smallskip
\noindent\emph{Korn estimates and normalization.} Korn's interpolation inequality \cite[Theorem 3.1]{H18} gives
\be\lb{sph:eq:odd-korn-interpolation}
\|\nabla V\|_{h,G}^2\le C\left(h^{-1}\|V_3\|_{h,G}\|e(V)\|_{h,G}+\|V\|_{h,G}^2+\|e(V)\|_{h,G}^2\right).
\ee
Its thickness assumptions hold since $ha$ and its first derivatives have size $O(h)$. On the flat midplane the normal component in that inequality is precisely $V_3$. For the local iterations one can cover the inner lateral domain by finitely many squares whose concentric enlargements are compactly contained in the outer domain. Use concentric squares with side lengths in a fixed positive interval in each iteration. The geometric constants in the interpolation inequality are uniform on these families. For a fixed $C^2$ lateral domain use interior rectangles and boundary strips between Lipschitz graphs of uniformly positive separation in rotated planar coordinates. Summing the local product terms uses Cauchy--Schwarz and finite overlap.

We will need both the clamped and the normalized forms
\be\lb{sph:eq:odd-clamped-korn}
h\bigl(\|\nabla\psi\|_{h,G}+\|\psi_3\|_{h,G}\bigr)\le C\|e(\psi)\|_{h,G},\qquad \psi=0\ \text{on the lateral boundary},
\ee
\be\lb{sph:eq:odd-free-korn}
\|\nabla(V-R)\|_{h,G}\le Ch^{-1}\|e(V)\|_{h,G},
\ee
for fields of the odd parity in the statement. In the second inequality $R=(-Zc,c\cdot y+c_0)$ is a rigid displacement. It can be chosen so that the flattened third component of $V-R$ and its two lateral derivatives all have zero average on $G\times(-1,1)$.

Here is a compactness proof of both assertions. If either uniform bound failed along $h\to0$, normalize the gradient norm to one. In the free case first impose the three indicated averages. In the clamped case use the zero lateral trace. Fixed-cylinder Poincar\'e bounds the third component by a constant, and oddness gives $\|V'\|_{h,G}\le Ch$. Failure would imply $\|e(V)\|_{h,G}/h\to0$.

The flattened $V_3$ is bounded in $H^1$, and $\partial_tV_3=hae_{33}(V)\to0$. Its limit is a function $w(y)$. The flattened $V'/h$ is also bounded in $H^1$. For its lateral derivatives use planar Korn's second inequality, since its planar strain is bounded by $C(\|e(V)\|_{h,G}/h+\|\partial_ZV'\|_{h,G})$. Its $t$-derivative is $a\partial_ZV'$. The transverse shear identity and parity identify the weak limit as $-at\nabla_yw$. The physical planar strain divided by $h$ consequently tends to $-atD_y^2w$. Thus $w$ is affine. Its zero trace in the clamped case, or its three zero averages in the free case, imply $w=0$. Rellich compactness gives $V_3\to0$ in $L^2$; also $V'\to0$. Equation~\eq{sph:eq:odd-korn-interpolation} then contradicts the normalized gradient norm. Uniformity under variations of $a$ follows by taking a $C^1$ convergent subsequence. For $h$ bounded away from zero ordinary Korn on the compact family of flattened domains gives the same conclusions.

For the solution in the lemma, subtract a vertical constant. Fixed-cylinder Poincar\'e and oddness then give
\[
\|U_3\|_{h,D}\le CK_h,\qquad\|U'\|_{h,D}\le ChK_h.
\]

\smallskip
\noindent\emph{An explicit local energy inequality.} For a lateral cutoff $\chi$ define
\[
\mathcal Q_\chi V=(\chi V'-ZV_3\nabla_y\chi,\ \chi V_3).
\]
Its strain satisfies
\be\lb{sph:eq:bending-cutoff}
\begin{aligned}
e_{\alpha\beta}(\mathcal Q_\chi V)&=\chi e_{\alpha\beta}(V)+(V'\odot\nabla_y\chi)_{\alpha\beta}-Z(\nabla_yV_3\odot\nabla_y\chi)_{\alpha\beta}-ZV_3\partial_{\alpha\beta}\chi,\\
e_{\alpha3}(\mathcal Q_\chi V)&=\chi e_{\alpha3}(V)-\tfrac12Z\partial_ZV_3\partial_\alpha\chi,\qquad e_{33}(\mathcal Q_\chi V)=\chi e_{33}(V).
\end{aligned}
\ee
Suppose $a_h(V,\psi)=\mathcal L(\psi)$, with $|\mathcal L(\psi)|\le M\|e(\psi)\|_{h,G_t}$ for odd tests supported laterally in $G_t$. Let $G_s\Subset G_t$ have cutoff gap $\delta=t-s\in(0,1)$. Choose $\chi=\zeta^4$, with $\zeta=1$ on $G_s$ and $|D^j\zeta|\le C_j\delta^{-j}$, $j\le2$. Put
\[
E_t=\|e(V)\|_{h,G_t},\qquad L_t=\|V'\|_{h,G_t}+h\|V_3\|_{h,G_t}.
\]
In \eq{sph:eq:bending-cutoff} write $e(\mathcal Q_\chi V)=\zeta^4e(V)+B$. The quotient $B/\zeta^2$, extended by zero where necessary, has norm at most
\[
C\left(\delta^{-2}L_t+h\delta^{-1}\|\nabla V\|_{h,G_t}+h\delta^{-1}E_t\right).
\]
Testing with $\mathcal Q_\chi V$, coercivity and Young's inequality therefore give
\[
E_s\le C\left(M+\delta^{-2}L_t+h\delta^{-1}\|\nabla V\|_{h,G_t}+h\delta^{-1}E_t\right).
\]
By \eq{sph:eq:odd-korn-interpolation},
\[
h\|\nabla V\|_{h,G_t}\le C\bigl((L_tE_t)^{1/2}+L_t+hE_t\bigr).
\]
Consequently, for every $0<\theta<1$,
\be\lb{sph:eq:odd-hole-filling}
E_s\le\theta E_t+C_\theta\delta^{-2}(M+L_t),\qquad0<h\le c\theta\delta.
\ee
Indeed, the term $C\delta^{-1}(L_tE_t)^{1/2}$ is at most $\theta E_t/2+C_\theta\delta^{-2}L_t$. The displayed restriction on $h$ absorbs the remaining $Ch\delta^{-1}E_t$. Here $c$ and $C_\theta$ are independent of $h,s,t$. For $V=U$ and $M=0$, the cutoff identity can instead be used directly with the bounds on $U$ and its gradient. It yields
\be\lb{sph:eq:odd-base-energy}
\|e(U)\|_{h,G}\le ChK_h\qquad(G\Subset D).
\ee

\smallskip
\noindent\emph{Source terms after tangential differentiation.} Set $b_i=\partial_i\log a$, $X_i=\partial_i+Zb_i\partial_Z$. These vector fields are tangent to both faces. Define
\be\lb{sph:eq:bending-derivative}
(\mathcal D_iV)_\alpha=X_iV_\alpha-Z(\partial_\alpha b_i)V_3,\qquad(\mathcal D_iV)_3=X_iV_3+b_iV_3.
\ee
It preserves parity and maps $(-Z\nabla w,w)$ to the Kirchhoff--Love field generated by $\partial_iw+b_iw$. Applying $X_i$ componentwise would introduce a shear strain $e_{\alpha3}=-b_i\partial_\alpha w/2$ that need not tend to zero with $h$. The cutoff $\mathcal Q_\chi$ also preserves the Kirchhoff--Love form, replacing $w$ by $\chi w$.

Equip the state space $\bR_{\rm sym}^{3\times3}\oplus\bR^3\oplus\bR$, with the Frobenius inner product on the strain factor and the Euclidean product on the remaining factors. Set
\[
\mathsf JV=(E,P,q)=(e(V),\,Z\nabla V_3,\,ZV_3).
\]
There is an exact identity
\be\lb{sph:eq:bending-state}
\mathsf J(\mathcal D_iV)=X_i(\mathsf JV)+\mathsf L_i\mathsf JV.
\ee
Write $(E_i,P_i,q_i)=\mathsf L_i(E,P,q)$. Its components are
\[
\begin{aligned}
(E_i)_{\alpha\beta}&=(\partial_\alpha b_i)(ZE_{\beta3}-P_\beta)+(\partial_\beta b_i)(ZE_{\alpha3}-P_\alpha)-q\,\partial_{\alpha\beta}b_i,\\
(E_i)_{\alpha3}&=b_iE_{\alpha3},\qquad(E_i)_{33}=2b_iE_{33},\\
(P_i)_\alpha&=(\partial_\alpha b_i)(ZP_3+q),\qquad(P_i)_3=b_iP_3,\qquad q_i=0.
\end{aligned}
\]
Every entry is polynomial in $Z$, with smooth coefficients depending on $a$. All coefficient bounds required below are uniform in $h$.

Let $\mathsf B_0(y)$ act as $\bC(y)$ on the strain factor and vanish on the other two factors. Thus $a_h(V,\psi)=\int_h(\mathsf JV)^*\mathsf B_0\mathsf J\psi$. Set
\[
\mathsf K_i=\mathsf L_i^*\mathsf B_0+\mathsf B_0\mathsf L_i-b_i\mathsf B_0-X_i\mathsf B_0,\qquad\mathfrak C_iB=X_iB+b_iB-\mathsf L_i^*B-B\mathsf L_i.
\]
The last term in $\mathsf K_i$ accounts for material variation. Since $\operatorname{div}X_i=b_i$ and $X_i$ is tangent to the faces, integration by parts gives
\be\lb{sph:eq:bending-commuted}
a_h(\mathcal D_iV,\psi)=-a_h(V,\mathcal D_i\psi)+\int_h(\mathsf JV)^*\mathsf K_i\mathsf J\psi.
\ee

For an ordered word $I$ in the indices $1,2$, put $U_I=\mathcal D_IU$, with $U_\varnothing=U$. Suppose the equation for $U_I$ has been written as
\[
a_h(U_I,\psi)=\sum_{|J|<|I|}\int_h(\mathsf JU_J)^*B_{I,J}\mathsf J\psi.
\]
For fixed $h>0$, local conormal regularity on one-face patches gives $U\in H^{N+1}$ up to either face, as detailed below for the terminal estimate. The identities for weak solutions then follow by density in the tests. The empty word has zero right-hand side. For a new leading index $i$, combining~\eq{sph:eq:bending-commuted} with~\eq{sph:eq:bending-state} and integrating by parts gives the exact recursion
\be\nn
a_h(U_{iI},\psi)=\int_h(\mathsf JU_I)^*\mathsf K_i\mathsf J\psi+\sum_{|J|<|I|}\int_h\bigl[(\mathsf JU_{iJ})^*B_{I,J}+(\mathsf JU_J)^*(\mathfrak C_iB_{I,J})\bigr]\mathsf J\psi.
\ee
All states on the right have order at most $|I|$. The recursion constructs the matrices $B_{iI,J}$, with identical words combined by addition. The coefficient bounds can be written as
\be\lb{sph:eq:bending-coefficient-bounds}
\sup_{\mathcal P_h(G)}|\partial_y^\alpha\partial_Z^\beta B_{I,J}|\le C_N,\qquad
1\le |I|\le N,\quad |J|<|I|,\quad |\alpha|+\beta\le N-|I|.
\ee
Indeed $\mathsf B_0$ is independent of $Z$, and $X_i$ preserves polynomial dependence on $Z$; each use of $\mathfrak C_i$ takes at most one additional lateral derivative of its coefficient. The initial matrix $\mathsf K_i$ uses derivatives of $a$ through order three and of $\bC$ through order one. For $|I|=j$ the recursion differentiates $\mathsf K_i$ at most $j-1$ times. A further $N-j$ derivatives in \eq{sph:eq:bending-coefficient-bounds} use at most $3+(j-1)+(N-j)=N+2$ derivatives of $a$ and $N$ of $\bC$. Products with $\mathsf L_i$ obey the same count. No inverse power of $h$ enters this recursion. By \eq{sph:eq:odd-clamped-korn}, $\|\mathsf J\psi\|_{h,G}\le C\|e(\psi)\|_{h,G}$. Consequently the commuted source satisfies the explicit bound
\be\lb{sph:eq:bending-source-norm}
|a_h(U_I,\psi)|\le C_N\sum_{|J|<|I|}\|\mathsf JU_J\|_{h,G}\|e(\psi)\|_{h,G}.
\ee

\smallskip
\noindent\emph{Derivative induction and the final thickness scale.} For $V=U_I$, we prove on nested fixed lateral domains that
\be\lb{sph:eq:odd-tangential-energy}
\|e(V)\|_h+h\|\nabla V\|_h+\|V'\|_h+h\|V_3\|_h\le C_jhK_h,\qquad|I|=j\le N.
\ee
The case $j=0$ follows from the normalization bounds and the strain estimate~\eq{sph:eq:odd-base-energy}. Suppose it holds through order $k-1$. Planar Korn on slices, with a fixed lateral cutoff, gives
\[
\|\nabla_y(V'\circ \Phi_h)\|_{L^2}\le C\bigl(\|e(V)\|_h+h\|\nabla V\|_h+\|V'\|_h\bigr),\qquad \Phi_h(y,t)=(y,ha(y)t).
\]
The difference of the flattened and physical planar strains is bounded by $Ch|\partial_ZV'|$. For $W=\mathcal D_iV$ of order $k$, this estimate and \eq{sph:eq:bending-derivative} show $\|W'\|_h+h\|W_3\|_h\le ChK_h$. By~\eq{sph:eq:bending-source-norm}, the source dual norm satisfies $M\le ChK_h$, since only states of orders at most $k-1$ contribute. Explicitly,
\[
\|\mathsf JU_J\|_h\le C\bigl(\|e(U_J)\|_h+h\|\nabla U_J\|_h+h\|(U_J)_3\|_h\bigr)\le ChK_h\qquad(|J|<k).
\]

Apply \eq{sph:eq:odd-hole-filling} to $W$. Choose expanding domains with successive cutoff gaps $\delta_n=\delta_0\,2^{-n}$, whose sum stays inside a fixed larger domain. Fix
\[
0<\theta<\min\{1/8,2^{-(k+3)}\}.
\]
For $h\le c\theta\delta_0$, let $n_h=\lfloor\log_2(c\theta\delta_0/h)\rfloor$, the last index for which $h\le c\theta\delta_{n_h}$. Iteration up to this index yields
\[
E_0(W)\le C_\theta hK_h\delta_0^{-2}\sum_{n=0}^{n_h}(4\theta)^n+\theta^{n_h+1}E_{n_h+1}(W).
\]
The sum is bounded by $(1-4\theta)^{-1}$. With $\alpha=-\log_2\theta>k+3$, the choice of the final index gives $\theta^{n_h+1}\le(h/(c\theta\delta_0))^\alpha$. Here $\delta_0$ is fixed before $h$ tends to zero.

For the terminal norm, cover the relevant layer by balls or one-face patches of diameter $c_0h$, with bounded overlap. Interior and homogeneous conormal estimates for the original solution $U$, after subtraction of local constant displacements, give $\|\nabla^{\ell+1}U\|_h\le C_\ell h^{-\ell}K_h$. Here $0\le\ell\le k\le N$. On each rescaled one-face patch the graph is bounded in $C^{k+1}$ and the flattened variational coefficients in $C^k$, uniformly in $h$. The local $H^{k+1}$ estimate uses at most $k$ tangential differences of the weak equation. The invertible normal principal block recovers derivatives containing two or more normal differentiations by the equation. Each patch meets at most one face. The frozen traction problem satisfies the complementing condition of \cite{ADN64}. Indeed, a decaying half-space solution at a nonzero tangential frequency, tested against its conjugate over one tangential period, has zero strain and hence vanishes. Ellipticity, the coefficient bounds and compactness of the frozen tensors and normals give uniform constants in these local estimates. The coefficients in $\mathcal D_I$ and their derivatives are uniformly bounded; the normalized third displacement is bounded by $CK_h$. Hence $E_{n_h+1}(W)\le C_kh^{-k}K_h$.

The terminal contribution is at most $Ch^3K_h$. This proves the strain bound in \eq{sph:eq:odd-tangential-energy}. The already controlled displacement components and \eq{sph:eq:odd-korn-interpolation} then give $\|\nabla W\|_h\le CK_h$, completing the induction. All smallness conditions on $h$ can be chosen once for the finitely many orders $k\le N$.

\smallskip
\noindent\emph{Recovery on the two faces.} To recover full gradients, write $a_i=(-Z\nabla_yb_i,b_i)$, so that $\mathcal D_iV=X_iV+a_iV_3$. For each physical derivative $\partial_k$, including $\partial_Z$, the exact rule is
\[
X_i(\partial_kV)=\partial_k(\mathcal D_iV)-\partial_k(Zb_i)\partial_ZV-a_i\partial_kV_3-(\partial_ka_i)V_3.
\]
At fixed $t$, $\partial_{y_i}(V\circ \Phi_h)=(X_iV)\circ \Phi_h$. The rule above and $X_iV_3=(\mathcal D_iV)_3-b_iV_3$ express $j$ flattened lateral derivatives of $\nabla U$ in terms of gradients and third displacement components of words of length at most $j$. For $j\le N$ the coefficients here use derivatives of $a$ through order $j+2$. For the stress, \eq{sph:eq:bending-state} and \eq{sph:eq:odd-tangential-energy} bound each differentiated strain by states of size $O(hK_h)$. Consequently, on each inner cylinder,
\[
 \|\sigma\circ\Phi_h\|_{L_t^2H_y^N}\le C_NhK_h,\qquad
 \|\nabla U\circ\Phi_h\|_{L_t^2H_y^N}\le C_NK_h,\qquad \sigma=\bC(y)e(U).
\]

Apply \eq{layer:eq:conormal-balance} with $\tau=h$.

With $N\ge3$ the right-hand side has $L_t^2H_y^2$ norm at most $Ch^2K_h$. Integration in $t$ and planar $H^2\hookrightarrow L^\infty$ on inner subdomains give $\|\mathcal T\|_\infty\le Ch^2K_h$. The same argument applied to
\[
 \partial_t\nabla_y W_3=h\nabla_y(a\, e_{33}(U)\circ\Phi_h),\qquad
 \partial_t\nabla_y W'=h\nabla_y(a\,\partial_ZU'\circ\Phi_h)
\]
gives $\|\nabla_yW_3-\nabla_yw_h\|_\infty\le Ch^2K_h$ and $\|\nabla_yW'\|_\infty\le ChK_h$. In the first estimate subtract the thickness mean, and in the second use $W'(y,0)=0$. Finally, averaging
\[
 \nabla_y W_3=(\nabla' U_3)\circ\Phi_h+ht\,\nabla_y a\,(\partial_Z U_3)\circ\Phi_h
\]
gives $\|\nabla_y w_h\|_{H^N(D_1)}\le C_NK_h$. The vertical normalization supplies the zeroth order bound. At $N=3$, $H^4$ on an inner planar domain embeds into $C^{2,\alpha}$ for every $0<\alpha<1$, and at $N=4$, $H^5$ embeds into $C^3$.

Apply \eq{layer:eq:conormal-recovery} with $n=(-ht\nabla_y a,1)$ and $p_0=(-\nabla_yw_h,0)$. The tangential estimates give $\operatorname{sym}(G_y+p_0\otimes n)=O(hK_h)$, since the leading terms form a skew matrix. Hence

\[
\mathsf A_{\bC}(n)(\partial_ZU\circ\Phi_h-p_0)=\mathcal T-\bigl(\bC(y)\operatorname{sym}(G_y+p_0\otimes n)\bigr)n
\]
therefore gives $\partial_ZU'=-\nabla_yw_h+O(hK_h)$ and $\partial_ZU_3=O(hK_h)$. Physical lateral derivatives differ from flattened ones by $ht(\partial_\alpha a)\partial_ZU$. Thus the physical gradient is $\mathcal R(\nabla_yw_h)+O(hK_h)$. Taking its symmetric part proves the two pointwise estimates. All estimates hold on the closed thickness interval.
\end{proof}

\medskip\noindent\textbf{Acknowledgments.}
J. Zhang thanks Dr. Yihao Gong for helpful discussions and suggestions. H. Li was partially supported by the National Natural Science Foundation of China (Grant No.~12471191) and the Beijing Natural Science Foundation (Grant No.~1242006).

\section*{Declaration of AI use}
Generative AI tools were used for English-language editing and to assist in checking portions of the mathematical exposition. Responsibility for the final statements, proofs, and citations rests with the author.

\end{document}